\documentclass{article}

\usepackage[T1]{fontenc}
\usepackage[utf8]{inputenc}
\usepackage{csquotes}
\usepackage[british]{babel}

\usepackage{mathtools}

\usepackage{amsmath}
\usepackage{amssymb}
\usepackage{amsfonts}
\usepackage{dsfont} 
\usepackage{amsthm}
\usepackage{thmtools}
\usepackage{tabularx}

\usepackage[font=small, singlelinecheck=false]{caption}

\usepackage[left=3.5cm,right=3.5cm]{geometry} 

\usepackage{xcolor} 

\usepackage{multirow} 
\usepackage{adjustbox} 
\usepackage{placeins} 

\usepackage[style = numeric, backend=bibtex]{biblatex}
\usepackage{enumitem} 
\usepackage{array} 

\NewDocumentCommand{\R}{}{\mathbb{R}}

\usepackage{hyperref}
\usepackage[capitalise]{cleveref}
\numberwithin{equation}{section}

\usepackage{authblk}

\title{Central limit theorems describing isolation by distance under various forms of power-law dispersal}
\date{January 2025}
\author[1]{Rapha\"el Forien}
\author[2]{Bastian Wiederhold}
\affil[1]{INRAE, BioSP, 84914, Avignon, France}
\affil[2]{Department of Statistics, University of Oxford, OX1 3LB, Oxford, UK}
\DeclareMathOperator{\blangle}{\big\langle}
\DeclareMathOperator{\brangle}{\big\rangle}
\DeclareMathOperator{\Vol}{Vol}

\DeclareMathOperator{\tail}{a}
\DeclareMathOperator{\power}{b}
\DeclareMathOperator{\impact}{c}
\DeclareMathOperator{\dimension}{d}

\theoremstyle{definition}
\newtheorem{definition}{Definition}
\numberwithin{definition}{section}

\newtheorem{remark}[definition]{Remark}

\theoremstyle{plain}
\newtheorem{theorem}[definition]{Theorem}
\newtheorem{proposition}[definition]{Proposition}
\newtheorem{corollary}[definition]{Corollary}

\newtheorem{lemma}[definition]{Lemma}

\begin{document}


\maketitle


\begin{abstract}
In this paper, we uncover new asymptotic isolation by distance patterns occurring under long-range dispersal of offspring. We extend a recent work of the first author, in which this information was obtained from forwards-in-time dynamics using a novel stochastic partial differential equations approach for spatial $\Lambda$-Fleming-Viot models. The latter were introduced by Barton, Etheridge and V\'eber as a framework to model the evolution of the genetic composition of a spatially structured population. Reproduction takes place through extinction-recolonisation events driven by a Poisson point process. During an event, in certain ball-shaped areas, a parent is sampled and a proportion of the population is replaced. We generalize the previous approach of the first author by allowing the area from which a parent is sampled during events to differ from the area in which offspring are dispersed, and the radii of these regions follow power-law distributions. In particular, while in previous works the motion of ancestral lineages and coalescence behaviour were closely linked, we demonstrate that local and non-local coalescence is possible for ancestral lineages governed by both fractional and standard Laplacians.
\end{abstract}


\section{Introduction} \label{sec:introduction}
\paragraph{Probability of identity by descent} In most natural settings, the physical expansion of the population is too broad to be visited by a single individual during its lifetime and several generations are necessary to bridge large spatial distances. Consequently, in the presence of neutral mutations, individuals further apart in space carry, on average, more different genetic information; an effect which is known as \textit{isolation by distance}. In \textit{infinitely-many-alleles models}, mutations are deemed unlikely to affect the same nucleotide twice, such that mutation events create new types never seen before in the population. Under this assumption, relatedness of individuals can be measured by the \textit{probability of identity by descent} $F(x)$, that is, the probability that two individuals sampled from two locations at separation $x$ are of the same type because no mutation occured since their most recent common ancestor.\par
Historically, the displacement distance between offspring and parents has been assumed to be of finite variance.  Wright and Mal\'ecot (\parencite{Wright}, \parencite{Malecot}) were the first to observe that if $ G_t $ is a Gaussian transition kernel associated with Brownian motion, $ F(x) $ is approximately given by the \textit{Wright-Mal\'ecot formula} (ignoring model specific parameters)
 \begin{equation} \label{introductionWMF}
 	F(x) \approx \int_0^\infty e^{-2 \mu t} G_{2t} (x) dt,
 \end{equation}
 where $\mu > 0$ is the rate at which neutral mutations arise in a given \textit{ancestral lineage or line of descent}.
 An ancestral lineage refers to the sequence of ancestors of an individual sampled from a population. As we consider spatial processes we effectively refer to the spatial random walk described by the positions of the ancestors as an ancestral lineage. If individuals in the population mutate at a certain rate $\mu$, with the same rate we see mutation along the ancestral lines in reversed time.
 
The formula reflects the idea of probability of identity by descent. 
Since $ \mu$ is the neutral mutation rate, $ e^{-2 \mu t} $ is the probability that no mutation occurred on the line of descent of either individual since the time $ t $ of their most recent common ancestor. Mutation rates estimated from biological data are often very small. If we wish to detect mutations, we therefore need to consider large time scales. Over large time and space scales, in turn, the central limit theorem tells us that jump processes of ancestral lineages are approximated by Brownian motions. 
Hence, the two ancestral lineages are at the same place and coalesce with a rate proportional to $ G_{2t} (x) $ after time $ t $. The universality of the central limit theorem in approximating jump processes and ancestral lineages explains why the Wright-Mal\'ecot formula has been proven for a large variety of different spatial models such as the stepping stone model or the Gaussian version of the spatial $\Lambda$-Fleming-Viot model (\parencite{CD02}, \parencite{neutralevolution2002}, \parencite{2010companionpaper}).

\paragraph{Long-range dispersal in biology} Over the past decades, the assumption of a finite-variance distance relationship between parents and offspring has been challenged by many biological examples such as butterfly populations (\parencite{powerlawbutterflies}), expansion of aquatic species (\parencite{poweraquaticspecies}) or the dispersal of plant seeds (\parencite{powerlawplantseeds}) and the spread of pollen (\parencite{powerlawpollen}). In plants, long-distance dispersal of seeds has even been proposed to be an instance of the Levy foraging hypothesis (\parencite{powerlawseedlevyflight}), which supposes that long jumps are more efficient in the search for sparse resources. \par

Studying the effect of events covering an exceptionally large area on genealogies was one of the motivations for the introduction of the spatial $\Lambda$-Fleming-Viot model (SLFV) and already \cite{anewmodel} explored deviations from Kingman's coalescent. However, only recently, the papers \parencite{forien} and Smith and Weissman \parencite{weissman} have started to classify regimes for probability of identity by descent patterns under long-range reproduction. 
Exceptionally large events can even persist in the limit and lead to ancestral lineages described by stable processes. Additionally, coalescence of lineages may not only happen \textit{locally} when lineages meet, but \textit{globally}, when lineages are still separated in space.
The first instance of a long-range Wright-Mal\'ecot formula was established in \parencite{forien}, in which a long-range behaviour of ancestral lineages always occurred together with global coalescence. 
This was complemented and confirmed by the simulation-based approach of \parencite{weissman}. Exploring dimensions one and two, Smith and Weissman discovered that the same stable behaviour of ancestral lineages can appear together with both local and global coalescence. 
The present paper includes and considerably extends those previous classifications by generalizing the new technique of \parencite{forien}.

\paragraph{Novel central limit theorem technique} 

The classical formula \eqref{introductionWMF} contains an innate mathematical difficulty: Brownian motions and thus ancestral lineages do not meet in dimension two or higher. 
Wright and Mal\'ecot as well as \parencite{neutralevolution2002} and \parencite{2010companionpaper} circumvented this problem by introducing a parameter for the local scale. 
Only recently, the paper \parencite{forien} provided a mathematically rigorous solution to this obstacle. 
The key is the observation that even in higher dimensions, in all prelimiting models a non-zero probability of coalescence remains, vanishing only in the limit. 
By scaling the probability of identity by descent as the model converges to the limit, \parencite[Theorem 2.4]{forien} has been able to recover the classical Wright-Mal\'ecot formula regardless of dimension.

A central limit theorem type result involves taking a limit of the scaled population process $\boldsymbol{\rho}_t^N \to \lambda$, corresponding to a law of large numbers, and considering the difference $\boldsymbol{\rho}_t^N - \lambda$ around the resulting deterministic limiting object $\lambda$. If the differences $\boldsymbol{\rho}_t^N - \lambda$ are scaled appropriately, in our case by a sequence we dubbed $\sqrt{N \eta_N}$, the fluctuations
\begin{equation} \label{introductionclt}
Z^N = \sqrt{N \eta_N} ( \boldsymbol{\rho}^N - \lambda)
\end{equation}
will converge to a stochastic partial differential equation (SPDE), which one can regard as a central limit theorem.
Strikingly, the fluctuations of the model in \cite{forien} and this work have a straightforward connection to the probability of identity by descent. The scaling of the fluctuations $\sqrt{N \eta_N}$ required for the central limit theorem will be precisely the scaling required for the convergence of the probability of identity by descent. Additionally, various other parts of the SPDE will directly reflect the Wright-Malecot formula: the operator, such as the Laplacian, corresponds to the motion of lineages and the spatial correlation of the noise of the SPDE will reflect the coalescence behavior. We believe this to be a wider phenomena, in a follow up work \cite{wmfvaryingsize} will illustrate this correspondence for a different class of SPDEs.

\paragraph{Spatial $\Lambda$-Fleming-Viot model - forwards and backwards}
The arguments of \parencite[Theorem 2.4]{forien} are based on a central limit theorem for the SLFV with mutation. The SLFV process was first introduced in \parencite{Eth08}, while the mutation mechanism was envisioned by \parencite{veber_spatial_2015}. The process is based on the idea that reproduction occurs through a sequence of spatial events, which replace a proportion of the population inside the affected area.

More precisely, individuals carry types described by real numbers from $[0,1]$ and are distributed over a geographic space $\mathbb{R}^{\dimension}$. Between reproduction events, individuals mutate by assuming new uniform labels from the interval $[0,1]$ with a certain mutation rate $\mu$, corresponding to the infinitely-many-alleles model. In each reproduction event, a ball $B(x,r)$ centered at $x \in \mathbb{R}^{\dimension}$ and with radius $r >0$ is affected. If there is no individual in the area, nothing happens. If there are individuals, we sample uniformly at random one of them to be the parent, kill each individual in the ball independently of each other with probability $u \in (0,1]$ and add new individuals of the parental type sampled from a Poisson process with intensity proportional to $u \mathds{1}_{B(x,r)}$, in such a way that the mean density of individuals is kept constant (see \cite{berestycki2009survival}). The proportion $u$ is known as the impact of an event. Location of events, as well as their parameters are sampled according to a Poisson process. 

The SLFV is then obtained as the limit of the above process when the spatial density of individuals tends to infinity.
The limit is described by a process $(\rho_t, t \geq 0)$ such that, for each $t \geq 0$, $\rho_t$ is a measurable map from $\R^{\dimension}$ to the space of probability measures on the set of genetic types (see \Cref{SLFVdefinition} below).
At each reproduction event, we update the map $\rho_t$ inside the ball affected by the event by replacing a proportion $u$ of the mass by a Dirac mass on the parental type, which is selected at random from the distribution of genetic types in the ball just before the event.
On the other hand, the mutations in the limit result in an exponential decay of $\rho_t$ towards the constant map which assigns the uniform distribution on $[0,1]$ to each point $x \in \R^{\dimension}$.
Even though the density of individuals is infinite in this process, the population is still subject to genetic drift as the allele frequencies fluctuate randomly at each reproduction event.
The strength of this genetic drift can be tuned by changing the impact parameter $u$, which affects the amount by which allele frequencies change at each reproduction event.

In population genetics one is often interested in the distribution of the genealogy of a sample of individuals.
Here, a sample of individuals is understood as a sample of genetic types from the distribution of types at some specified locations.
Using the reversibility of the Poisson point process driving the reproduction events, \cite{anewmodel} have shown that the SLFV admits a genealogical dual which can be described as a system of coalescing random walks.
Following the genealogy of a sample backwards-in-time, at each reproduction event, we mark each lineage present in the ball affected by the event independently of each other with probability $u$, and all marked lineages coalesce and move to a new location selected uniformly at random from the same ball.
This backwards-in-time process is connected to the forwards-in-time SLFV through a mathematical relationship called duality (see \Cref{sec:proof_existence} for a statement in the context of the present paper).

How are we able to consider a law of large numbers if the process is already an infinite-population limit? The SLFV retains a notion of finite-population size through the size of the impacts $u$. Smaller impacts correspond to smaller stochastic changes in the composition, so a smaller genetic drift, which can be seen as considering larger population sizes.

Throughout this work, we will switch between the forwards- and backwards-in-time perspectives on the process. The main results will be obtained forwards-in-time and through the Wright-Malecot formula we will give interpretations in terms of the backwards-in-time behavior of lineages. Under the considered limits, the ancestral lineages describing the position of the ancestors will converge to Brownian motions or stable Lévy processes with surprising coalescing behavior, described at the end of the introduction. We expect that these dynamics, which we obtain through the central limit theorem, could be obtained by working directly with the backwards-in-time genealogical dual of the process under the appropriate scalings. Proving this explicitly goes beyond the scope of this work; the crux of the central limit theorem technique is that one can work forwards-in-time.
	
\paragraph{Formulae for short and long-range dispersal} \label{paragraphshortlong}	Crucially, the central limit theorem approach yields a generalisation of the Wright-Mal\'ecot formula to populations with long-range reproduction events, that is to say when the distribution of the radius of balls in which events take place has a heavy tail $ r^{- ( 1+ d+ \alpha)} dr $ with $\alpha < 2$ (if $\alpha \geq 2$ we degenerate to the classical case, as this work will show). 
Under a suitable spatial and temporal rescaling, ancestral lineages converge in the long-range case to stable processes instead of Brownian motions (see e.g. \parencite{forien_central_2017}). 
In the new approach, the probability of identity by descent is evaluated for a pair of individuals sampled independently at two locations chosen at random from two probability densities $\phi, \psi$. 
Let $ P_t^N (\phi, \psi) $ denote this probability for a population evolving according to the $ N $-th rescaled SLFV model $ (\boldsymbol{\rho}_t^N, t \geq 0) $. 
It is shown in \cite{forien} that
	\begin{equation*}
		\lim_{t \to \infty} \lim_{N \to \infty} N \eta_N P_t^N (\phi, \psi) = \int_{(\mathbb{R}^{\dimension})^2}  F ( \vert x -y \vert) \phi (x) \psi (y) dx dy,
	\end{equation*}
	for some function $ F: \mathbb{R}_+ \rightarrow \mathbb{R}_+ $ which depends on the dispersal regime. 
As mentioned before, the scaling factor $N \eta_N$ from the central limit theorem \eqref{introductionclt} is precisely the scaling factor of the probability of identity by descent to recover a non-trivial formula, i.e. captures how quickly the probability of identity by descent converges to zero as the limit is approached.
	In the short-range case, where the event radius is a constant $ R > 0 $,
	\begin{equation} \label{introductionshortrange}
		F(\vert x - y \vert) = \int_0^\infty e^{-2 \mu t} G_{2t} (\vert x - y\vert) dt.
	\end{equation}
This is the classical Wright-Mal\'ecot formula, see \cite{neutralevolution2002}. 
In the long-range case, where the event size distribution has a heavy tail $ r^{- ( 1+ d+ \alpha)} dr$, ancestral lineages become stable processes with parameter $ \alpha \in (0,2) $. 
In this case, it was shown that
	\begin{equation} \label{introductionlongrangeformula}
		F (\vert x -y\vert) = \int_0^\infty \int_{(\mathbb{R}^{\dimension})^2} e^{- 2 \mu t} G_t^\alpha ( x- z_1) G_t^\alpha (y - z_2) K (z_1, z_2) dz_1 dz_2 dt,
	\end{equation}
	where 
	\begin{equation} \label{kquantity}
	K (z_1, z_2) = \frac{C}{\Vert z_1 - z_2 \Vert^{d + \alpha}}.
	\end{equation}
 Reflecting the convergence of ancestral lineages, the Brownian transition kernel $ G_t $ is replaced by a stable transition kernel $ G_t^\alpha $. The formula involves a new term $ K(z_1, z_2) $, which appears as the spatial correlations of the noise of the SPDE arising in the central limit theorem. The correlations $ K(z_1, z_2) $ encode the rate of coalescence for ancestral lineages located at positions $ z_1 $ and $ z_2 $. 
	
\paragraph{New phenomena} In this paper, we will expand the technique of \parencite{forien} to more general models. 
The broader regimes will unveil some new phenomena, a detailed interpretation can be found in \Cref{subsec:interpretation}.

In \parencite{forien} there was a strong link between the behaviour of ancestral lineages and asymptotic coalescence, while our regimes show that this is no necessity. The short-range Wright-Mal\'ecot formula \eqref{introductionshortrange} can be interpreted as having a coalescence/spatial correlation term $K(x_1, x_2) = C \delta_{x_1} (x_2)$ for a constant $C > 0$. This indicates that coalescence happens only locally and ancestral lineages need to be at the same position $x_1 = x_2$ in order to coalesce asymptotically. In the long-range regime \eqref{introductionlongrangeformula}, ancestral lineages became stable processes and the coalescence/correlation term $ K (x_1, x_2) $ of \eqref{kquantity} indicated that coalescence happens globally over large distances. In this work, we will see mixtures of both behaviours: Brownian motions with global coalescence and stable ancestral lineages which still only coalesce locally.
More generally, in most cases, local coalescence or a range of different powers $\alpha$ appearing in \eqref{kquantity} is possible for a fixed either stable or Brownian behaviour of ancestral lineages (see \Cref{centrallimittheorem} and the corresponding formulae in \Cref{wmf}).

Furthermore, we will demonstrate that the coalescence/spatial correlation behaviour shows a dimension-dependent phase transition. We will identify a critical threshold for the tail behaviour below which global coalescence occurs and above which only local coalescence takes place. 
The critical case will exhibit local coalescence and requires a special scaling regime.

Our more general SLFV models involve two different radii $r_1$ and $r_2$ in reproduction events (\Cref{SLFVdefinition}). Inside a ball of radius $r_1$ a parental type will be chosen. Thereafter, a proportion of the population will be replaced by the parental type inside a ball of radius $r_2$. 
Our new Wright-Mal\'ecot formulae show that the joint distribution of $r_1$ and $r_2$ determines the asymptotic behaviour of lineages while only the $r_2$ marginal determines the coalescence regime.
Moreover, if both tails have independent heavy-tailed distributions, one of the tails dominates the parameter of the resulting stable lineages (see \eqref{one_tail_regime}, \eqref{two_tail_regime} and Table \ref{table:parameters}).\par
For potential inference of long-range dispersal patterns, this new variety of formulae has far-reaching consequences, which we hope to cover in future work.
The formulae illustrate that long-range reproduction leaves a clear trace in isolation by distance patterns.
However, as we obtain many situations which lead to similar or even the exact same formulae, inference of more precise model parameters might be hindered.\par
The Wright-Mal\'ecot formulae have been obtained entirely forwards-in-time, raising the hope that the techniques can be applied to other complex spatial models with difficult backwards-in-time behavior and untractable dual processes. 
Indeed, in the related work \cite{wmfvaryingsize}, we employ the strategy of this paper to a model, in which the dual process is too complex to describe due to varying population size.

\paragraph{Outline} In \Cref{sec:definition}, we introduce the SLFV models considered in this work. The central limit theorem and the corresponding Wright-Mal\'ecot formulae are stated in \Cref{sec:mainresults}. In \Cref{sec:proof_outline}, we give an overview of the proof strategy of \parencite{forien} and the adaptations needed in the present case. \Cref{sec:proof_lemmas} provides proofs of all steps requiring significant changes from \parencite{forien} to cover the new regimes. \Cref{sec:proof_existence} proves the existence and uniqueness of the SLFV models of \Cref{sec:definition}.

\paragraph*{Acknowledgements} The authors wish to thank Alison Etheridge for valuable suggestions and feedback. Furthermore, the authors are grateful to the Hausdorff Institute for Mathematics Bonn for its hospitality in the course of the 2022 Junior Trimester Program "Stochastic modelling in life sciences". This program was funded by the Deutsche Forschungsgemeinschaft (DFG, German Research Foundation) under Germany's Excellence Strategy - EXC-2047/1 - 390685813.
R.F. was supported in part by the Chaire Modélisation Mathématique et Biodiversité of Veolia-\'Ecole Polytechnique-Muséum National d'Histoire Naturelle-Fondation X.
B.W. was supported by the Engineering and Physical Sciences Research Council Grant [EP/V520202/1].

\section{Definition of the model} \label{sec:definition}
This paper makes use of spatial $ \Lambda $-Fleming-Viot (SLFV) processes, a framework introduced in (\parencite{Eth08}, \parencite{anewmodel}) in which the evolution of a population is driven by a family of reproduction events.
These events affect balls of varying radius and in each event, a fraction of the population inside the ball is replaced by the offspring of one or a few individuals. 
The SLFV process presented here differs in two ways from the most commonly studied version (\parencite{anewmodel}): a mutation mechanism, introduced in \parencite{veber_spatial_2015}, and different areas for the origin of parents and dispersal of offspring. We first define the state space $ \Xi $ of the process $ (\rho_t, t \geq 0) $, then we clarify its evolution. \par

Let $ \mathcal{M} \big([0,1]\big) $ denote the set of probability measures on the interval $ [0,1] $. The state space $ \Xi $ is defined as the set of maps
\begin{equation*}
	\rho : \mathbb{R}^{\dimension} \rightarrow \mathcal{M} \big([0,1] \big),
\end{equation*}
where we identify two maps $ \rho_1, \rho_2 $, if $ \{ x \in \mathbb{R}^{\dimension} : \rho_1 (x) \neq \rho_2 (x)\} $ is a Lebesgue null set. This space is in one-to-one correspondence with the space of measures on $\mathbb{R}^{\dimension} \times [0,1]$ with Lebesgue measure as a spatial marginal and probability measures as the marginals on $[0,1]$.
At time $ t $, the measure $ \rho_t(x) $ describes the probability distribution of the type of an individual sampled from position $ x $. We will frequently switch between both perspectives on the state space throughout this work.
To introduce a suitable topology on the space $\Xi$, let $ (\phi_n,n \geq 1), \phi_n : \mathbb{R}^{\dimension} \times [0,1] \rightarrow \mathbb{R} $ be a sequence of continuous functions with compact support on $\mathbb{R}^{\dimension}$, which separate elements of $ \Xi $. For $\rho_1, \rho_2 \in \Xi$ we define the metric $ d $ by
\begin{equation} \label{equationvaguemetric}
	d (\rho_1, \rho_2) = \sum_{n = 1}^\infty \frac{1}{2^n} \big\vert \langle \rho_1, \phi_n \rangle - \langle \rho_2 , \phi_n \rangle \big\vert,
\end{equation}
which induces the vague topology on $\Xi$. 
It is useful, for later proofs, to define $p$-norms as in \parencite{forien} only on the spatial argument of the test functions
\begin{equation} \label{equationlpnorm}
	\Vert \phi \Vert_p := \Bigg(\int_{\mathbb{R}^{\dimension}} \sup_{k \in [0,1]} \vert \phi (x, k) \vert^p dx \Bigg)^{1/p}.
\end{equation}
We can assume without loss of generality that for multi-indices $ \kappa \in \mathbb{N}^d $
\begin{equation} \label{eq:derivativesbounded}
	\sup_{n \in \mathbb{N}} \Big\{ \Vert \partial_\kappa \phi_n \Vert_p : \vert \kappa \vert \leq 2 , p \in \{ 1, \infty \} \Big\} < \infty.
\end{equation}
The next definition describes the evolution of the spatial $ \Lambda $-Fleming-Viot process. It extends \parencite[Definition 1.1]{forien} by introducing balls of different sizes for the selection of parents and dispersal of offspring and allowing the impact of reproduction events to be random. 

\begin{definition} \label{SLFVdefinition}
	Suppose we are given a positive constant $\mu$ and a $\sigma$-finite measure $ \nu $ on $ (0,1] \times (0, \infty)^2 $, which encodes the event impact and the event radii and satisfies
		\begin{align} \label{existence_condition}
		    \int_{(0,1] \times (0,\infty)^2} u r_2^d \, \nu(du, dr_1, dr_2) < \infty.
		\end{align}
	Then let $ \Pi $ be a Poisson point process on $ \mathbb{R}_+ \times \mathbb{R}^{\dimension} \times (0,1] \times (0, \infty)^2 $ with intensity $ dt \otimes dx \otimes \nu(du, dr_1 , dr_2) $.
	The spatial $ \Lambda $-Fleming-Viot process with mutation (SLFV) is defined as a stochastic process $ \rho = (\rho_t, t \geq 0) $ taking values in $ \Xi $ which evolves through reproduction and mutation. 
	\begin{enumerate}
		\item Reproduction: For every point $ (t,x,u,r_1, r_2) \in \Pi $, first choose a spatial position $ y $ uniformly at random from within the parent search area $ B(x,r_1) $ and a parental type $ k_0 \in [0,1] $ according to $ \rho_{t^-} (y, dk) $. Then, replace a proportion $ u $ of the population inside the dispersal area $ B(x,r_2) $, i.e. set
		\begin{equation}
			\forall z \in B(x,r_2), \enspace \rho_t (z, dk) = (1 - u) \rho_{t^-} (z, dk) + u \delta_{k_0} (dk).
		\end{equation}
		\item
		Mutation: If a point $ x \in \mathbb{R}^{\dimension}$ is hit by events at times $ t_1 $ and $ t_2 $ with no events in between, we set for $ t_1 \leq s < t_2 $,
		\begin{equation} \label{equationmutation}
			\rho_s (x, dk) = e^{- \mu (s -t_1)} \rho_{t_1} (x, dk) + \left(1 - e^{- \mu (s- t_1 )}\right) dk.
		\end{equation} 
	\end{enumerate}
\end{definition}

\begin{figure}
    \centering
    \includegraphics[width=\textwidth]{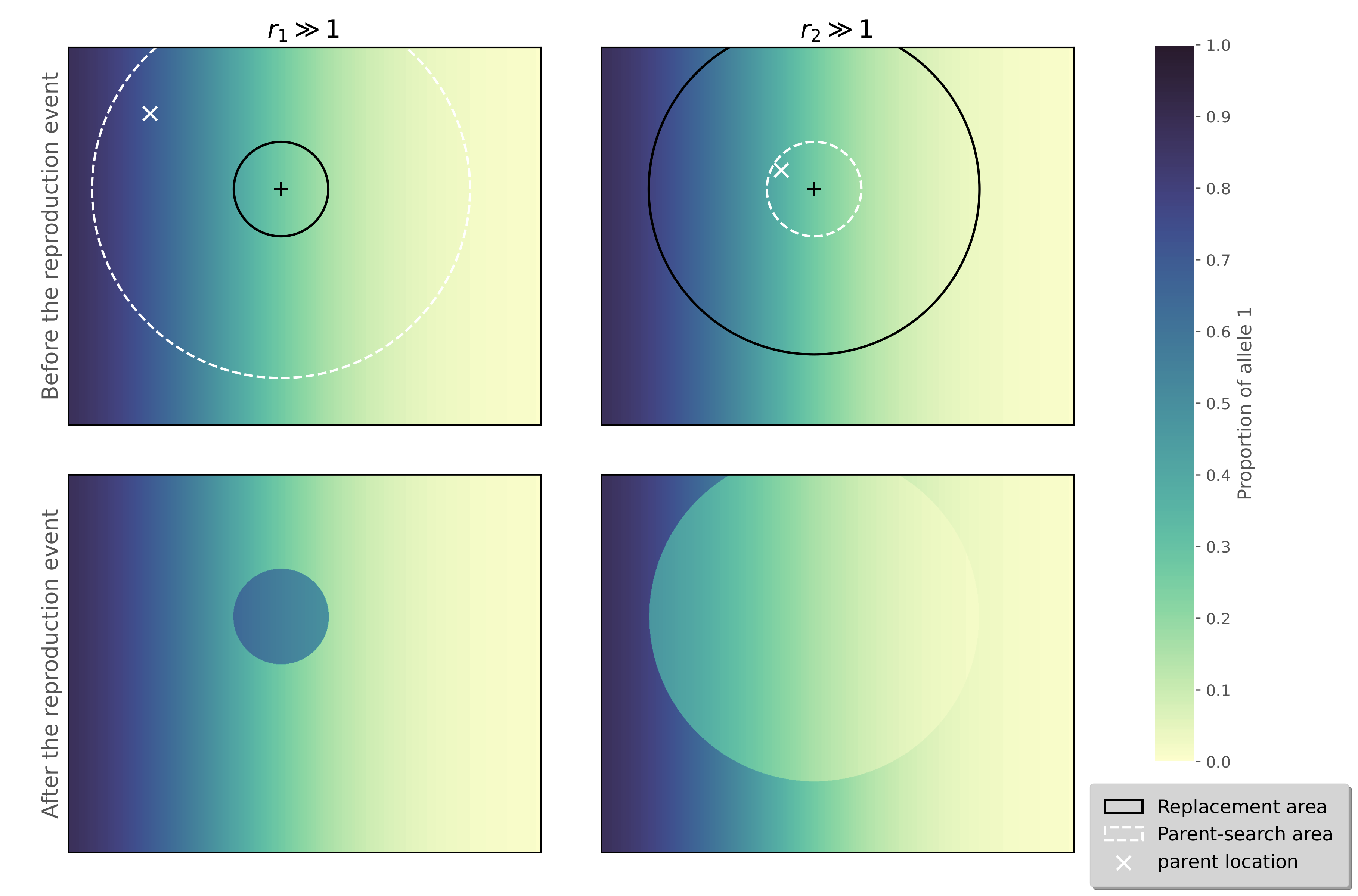}
    \caption{Example of the effect of a single reproduction event on the frequency of one allele in a spatial $\Lambda$-Fleming-Viot model with two alleles ($0$ and $1$). The characteristics of the reproduction event are $(t, x, u, r_1, r_2)$ with $u = 0.4$, and we plot $z \mapsto \rho_{t^-}(z, \lbrace 1 \rbrace)$ on the top row and $z \mapsto \rho_t(z, \lbrace 1 \rbrace)$ on the bottom row. On the left, the parent-search radius $ r_1 $ is larger than the replacement radius $ r_2 $ and the parent is selected at a point with high proportion of allele 1, even though this allele is relatively rare in the replacement area. We then plot the most likely outcome: that $k_0 = 1$, and a fraction $u$ of the population in the ball of radius $r_2$ is replaced by new offspring carrying allele 1. On the right, the parent-search radius $ r_1 $ is small and the replacement radius $r_2$ is relatively large. The parent is then selected near the centre of the event, where allele 1 is less abundant. We again plot the most likely outcome: that $k_0 = 0$, and a fraction $u$ of the population in a large area is replaced by new offspring carrying allele 0. In both cases, the average distance between parent and offspring is large, but the spatial correlations induced on allele frequencies are very different. Note also that, when $r_1 > r_2$, the region occupied by allele 1 may become disconnected.}
    \label{fig:first_event}
\end{figure}

In Figure~\ref{fig:first_event}, we show the effect of a reproduction event as described above on the frequency of a single allele in two cases: if the parent search radius $r_1$ is larger than the replacement radius $r_2$ and if the former is smaller than the latter.

The mutation mechanism is specifically tailored to capture probabilities of identity by descent. In infinite-alleles models, if an individual acquires a mutation, a new type is created. In other words, identity in state is almost surely equivalent to identity by descent. This can be captured by a mutation kernel with no atoms; our choice of the uniform distribution is arbitrary. In the large population limit, uniform mutation at rate $ \mu $ then corresponds to an exponential decay towards a uniform distribution. In \cite{forien}, more general mutation mechanisms are also considered.

\cref{SLFVdefinition} allows us to define the infinitesimal generator of the SLFV with mutations on a suitable set of functions.
For $\rho \in \Xi$, $x \in \mathbb R^d$, $r \in (0,\infty)$, $u \in (0,1]$ and $k_0 \in [0,1]$, define
\begin{align*}
    \Theta_{x,r,u,k_0}(\rho)(z,dk) := \begin{cases}
    (1-u) \rho(z,dk) + u \delta_{k_0}(dk) & \text{ if } z \in B(x,r), \\
    \rho(z,dk) & \text{ if } z \notin B(x,r).
    \end{cases}
\end{align*}
Let $D(\mathcal{G})$ denote the set of functions $F_{f,\phi} : \Xi \to \R$ of the form $F_{f,\phi}(\rho) = f(\langle \rho, \phi \rangle)$, where $\phi : \R^{\dimension} \times [0,1] \to \R$ is continuous and compactly supported, and $f : \R \to \R$ is continuously differentiable.
The infinitesimal generator of the SLFV of \cref{SLFVdefinition} is the operator $\mathcal{G}$, defined on $D(\mathcal{G})$ as
\begin{align*}
    &\mathcal{G} F_{f,\phi}(\rho) := \mu f'(\langle \rho, \phi \rangle) \langle \lambda - \rho, \phi \rangle \\ 
    & \hspace{2cm} + \int_{\mathbb R^d} \int_{(0,1] \times (0,\infty)^2} \frac{1}{|B(x,r_1)|} \int_{B(x,r_1)} \int_{[0,1]} \left( F_{f,\phi}(\Theta_{x,r_2,u,k}(\rho)) - F_{f,\phi}(\rho) \right) \\
    & \hspace{10cm} \rho(y,dk) dy\, \nu(du, dr_1, dr_2)\, dx,
\end{align*}
where $|B(x,r_1)|$ denotes the Lebesgue measure of $B(x,r_1)$ and $\lambda$ denotes the constant map which associates Lebesgue measure on $[0,1]$ to each $x \in \R^{\dimension}$, i.e. $\lambda(x,dk)=dk$.
We then say that a $\Xi$-valued process solves the martingale problem associated to $(\mathcal{G}, D(\mathcal{G}))$ if, for any $F \in D(\mathcal{G})$,
\begin{align*}
    F(\rho_t)-F(\rho_0)-\int_0^t \mathcal{G}F(\rho_s) ds
\end{align*}
is an $\mathcal{F}_t$-local martingale, where $(\mathcal{F}_t, t \geq 0)$ is the natural filtration associated to the process $(\rho_t, t \geq 0)$.

\begin{theorem} \label{thm:existence}
Under condition \eqref{existence_condition}, there exists a unique Markov process taking values in $\Xi$ which solves the martingale problem associated to $(\mathcal{G}, D(\mathcal{G}))$.
\end{theorem}

We prove \cref{thm:existence} in \cref{sec:proof_existence}, extending the arguments used in the proof of Theorem~1.2 in \cite{etheridge_rescaling_2020}.
Hereafter, the unique solution to the martingale problem associated to $(\mathcal{G}, D(\mathcal{G}))$ is called the SLFV with mutations (or simply the SLFV), with impact and radii intensity measure $\nu$ and mutation rate $\mu$.

\section{Main results} \label{sec:mainresults}

This paper generalizes the results of \parencite{forien} to a spectrum of regimes exhibiting new phenomena. 
Its centrepiece is a central limit theorem for the fluctuations of the process, from which we deduce the asymptotic behaviour of pairs of ancestral lineages.

\subsection{Choice of parameters}

\begin{figure}
\includegraphics[width=\textwidth]{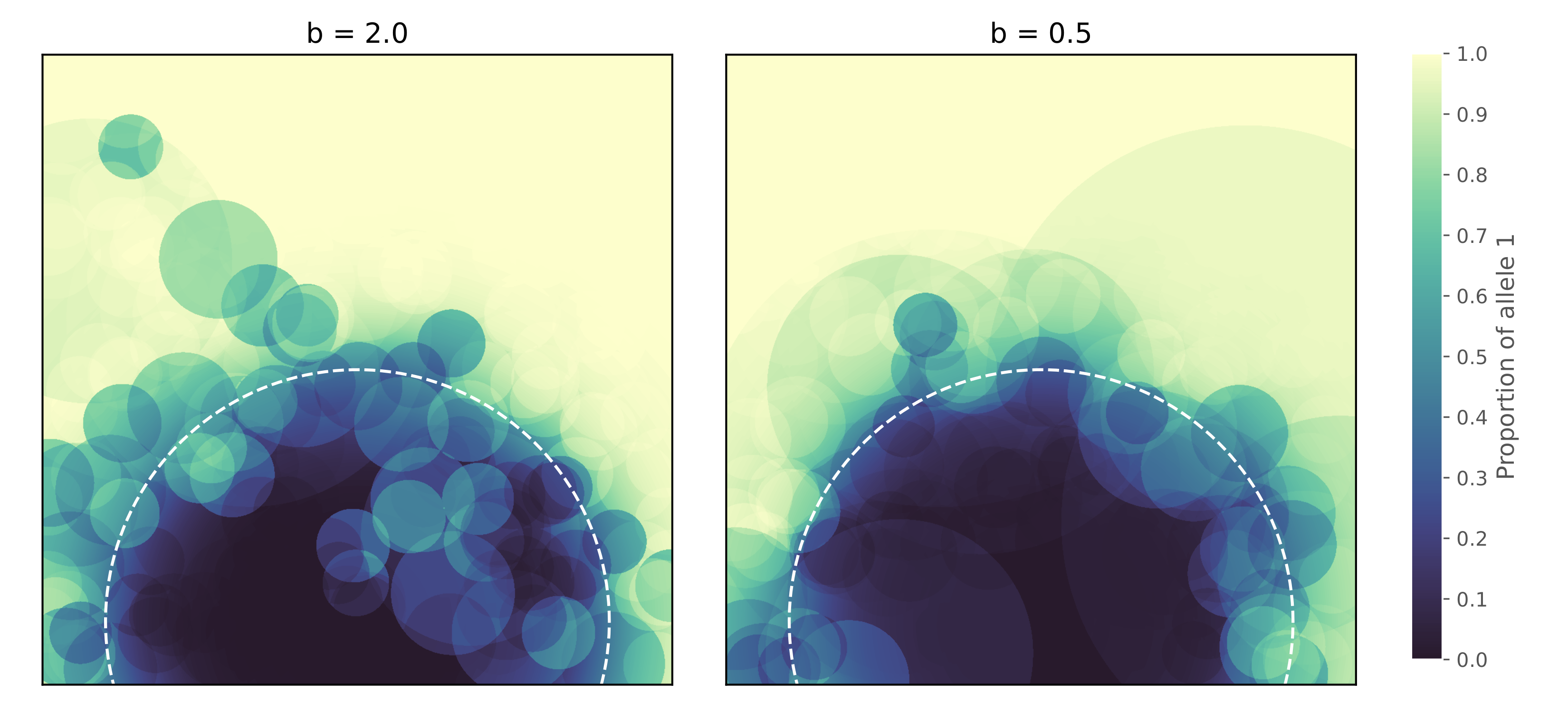}
\caption{Two realisations of the spatial $\Lambda$-Fleming-Viot model with two alleles (0 and 1) in two dimensions. Initially, the population inside the white circle is monomorphic of type 0, while the population outside the circle is monomorphic of type 1. The plots show the proportion of the type 1 allele after some time. The parameters of the SLFV are chosen according to the one-tail regime, one where $\power = 2$ (on the left) and one where $\power = 0.5$ (on the right). In both simulations, the parameter $\tail$ is chosen so that $\alpha = 1.3$, and $\impact = 1$.} \label{fig:two_regimes}
\end{figure}

We consider two possible choices for the intensity measure $\nu$, which are interpreted in detail in view of the main results in Subsection \ref{subsec:interpretation}. In the \textit{one-tail regime}, the premise is that the parental search radius $r_1$ is always a power of the replacement radius $r_2$. The intensity measure $\nu$ governing the reproduction events is chosen as
\begin{equation} \label{one_tail_regime}
\nu_N (du, dr_1 , dr_2) = \mathds{1}_{r_2 \geq 1}\, \delta_{u_N r_2^{- \impact}} (du) \, \delta_{r_2^{\power}} (dr_1)\, \frac{dr_2}{r_2^{1 + \dimension + \tail - \impact}},
\end{equation}
where $u_N \in (0,1]$, $\tail \in (0, \infty) $ and $ \power, \impact \in [0, \infty)$ are parameters. 
Here, ancestral lineages will converge to symmetric $\alpha$-stable processes with parameter
\begin{align} \label{def:alpha_one_tail}
    \alpha = \tail \land \hspace{0.5mm} \frac{\tail}{\power} \land 2 .
\end{align}
The parameter $\impact$ appears in the tail distribution of the radius $r_2$ and the impact $u = u_N r_2^{- \impact}$. When we consider the behaviour of lineages, both appearances annihilate each other, so that the parameter of the stable processes \eqref{def:alpha_one_tail} does not depend on $\impact$ (see Subsection \ref{subsec:interpretation} for details). Introducing the parameter $\impact$ allows us to modulate the impact of large events. For example, if $\impact = \dimension$ the offspring mass $uV_{r_2}$ produced is the same for every reproduction event regardless of the radius $r_2$.
Examples of realisations of the SLFV in the one-tail regime are displayed in Figure~\ref{fig:two_regimes}. 

In contrast, in the \textit{two-tails regime}, the parental search radius $r_1$ and the replacement radius $r_2$ are sampled independently of each other accordingly to heavy-tailed distributions. The intensity measure $\nu$ has the form
\begin{equation} \label{two_tail_regime}
\nu_N (du, dr_1, dr_2) = \mathds{1}_{r_1 \geq 1} \, \mathds{1}_{r_2 \geq 1} \, \delta_{u_N r_1^{-\impact_1} r_2^{-\impact_2}} (du) \, \frac{dr_1 dr_2}{r_1^{1 + \tail_1 - \impact_1} r_2^{1 + \dimension + \tail_2 - \impact_2}},
\end{equation}
where $u_N \in (0,1]$, $\tail_1 , \tail_2 \in (0, \infty)$ and $ \impact_1, \impact_2 \in [0, \infty)$ are parameters. In this case, limiting ancestral lineages will be symmetric $\alpha$-stable processes of parameter
\begin{align} \label{def:alpha_two_tails}
    \alpha = \tail_1 \land \tail_2 \land \, 2.
\end{align}
It can then easily be checked that both these measures satisfy \eqref{existence_condition}. 
Indeed,
\begin{align*}
    &\int_{(0,1] \times (0, \infty)^2} u r_2^{\dimension} \, \nu_N(du, dr_1, dr_2)= \begin{cases}
    u_N \int_1^{+\infty} \frac{dr_2}{r_2^{1+\tail}} & \text{ in the one-tail regime}, \\
    u_N \int_1^{+\infty} \frac{dr_1}{r_1^{1+\tail_1}} \int_1^{+\infty} \frac{dr_2}{r_2^{1+\tail_2}} & \text{ in the two-tails regime.}
    \end{cases}
\end{align*}
In both regimes, the right hand side is finite. 

\begin{remark}
The paper \parencite{forien} considers the one-tail regime with parameters $b = 1 $ and $ c = 0$, which amounts to considering only a single radius $r_1 = r_2$ and deterministic impact. Our results augment the corresponding results of \parencite{forien} to the above broader parameter setting and an interpretation thereof is given in \Cref{subsec:interpretation}.
\end{remark}

\subsection{Choice of scaling} \label{subsec:choiceofscaling}
Our goal is to derive a central limit theorem for the process. To derive a law of large numbers in \Cref{theo:convergencetolebesgue}, we apply the following scaling, identified in \parencite{forien}, to the process. The scaling will be explained in \Cref{rmkscaling}. 
We first choose a sequence $ (\delta_N, N \geq 1) $ of positive real numbers which converges to zero as $N \to \infty$.
Let $(\rho^N_t, t \geq 0)$ be the SLFV of \cref{SLFVdefinition}, with
 	\begin{equation*}
 		u_N = \frac{u_0}{N}, \hspace{1.5cm} \mu_N = \delta^\alpha_N \frac{\mu}{N},
 	\end{equation*}
 	where $\alpha$ is given by \eqref{def:alpha_one_tail} in the one-tail regime and by \eqref{def:alpha_two_tails} in the two-tails regime.
The rescaled process is then defined as
 	\begin{equation*}
 		\boldsymbol{\rho}_t^N (x, dk) := \rho^N_{Nt/\delta^\alpha_N } \left(\frac{x}{\delta_N} ,  dk\right).
 	\end{equation*}
The sequence of rescaled processes $ (\boldsymbol{\rho}_t^N, t \geq 0) $ converges to a deterministic limit given by Lebesgue measure $\lambda$. 
The convergence will be made precise in \Cref{theo:convergencetolebesgue}. This allows us to consider the fluctuations
\begin{equation} \label{scalingoffluctuations}
Z_t^N := \sqrt{N \eta_N} (\boldsymbol{\rho}_t^N - \lambda),
\end{equation}
where $\eta_N$ is a sequence of real numbers capturing the magnitude of fluctuations. 
In the one-tail regime, we set
\begingroup
\allowdisplaybreaks
\begin{align*}
&\eta_N := 
\begin{dcases}
    \delta_N^{\alpha - \tail - \impact} & \text{ if } \tail + \impact < \dimension, \\
   \delta_N^{\alpha - \dimension} \log (1/\delta_N)^{-1} & \text{ if } \tail + \impact = \dimension,\\
    \delta_N^{\alpha - d}  & \text{ if } \tail + \impact > \dimension,
\end{dcases}\\
\intertext{while in the two-tail regime,}
&\eta_N := 
\begin{dcases}
    \delta_N^{\alpha - \tail_2 - \impact_2} & \text{ if } \tail_2 + \impact_2 < \dimension, \\
   \delta_N^{\alpha - \dimension} \log (1/\delta_N)^{-1} & \text{ if } \tail_2 + \impact_2 = \dimension,\\
    \delta_N^{\alpha - d}  & \text{ if } \tail_2 + \impact_2 > \dimension.
\end{dcases}
\end{align*}
\endgroup

\begin{remark} \label{rmkscaling}
We wish to clarify how our choice of scaling aligns with our aim to show a central limit theorem for the process. To obtain a deterministic law of large numbers, we require the impact $u_N$ to converge to zero and would like to consider large spatial scales. Without loss of generality, we can assume $ u_N = u_0/N $ and scale space by $1/ \delta_N$ for a sequence of positive reals converging to zero. 

Ancestral lineages need to be among the impacted individuals in order to jump, which causes the speed of an ancestral lineage to be proportional to the average impact of events (\parencite{anewmodel}). If we sample individuals at rescaled distance $D$, which in unrescaled spatial units is proportional to $1/\delta_N$, it will therefore take roughly a time $ N D^\alpha \approx N /\delta_N^\alpha$ to cover this distance. We scale time by $N /\delta_N^\alpha$, so that the time to cover a distance $D$ is of order $O(1)$ and lineages are capable of movement as we approach the limit. To be able to detect mutation and genetic events in the population, the expected number of mutations along the lineages needs to be of order $O(1)$ and as $1/ \mu \sim N /\delta_N^\alpha$ we scale the mutations by the inverse.
\end{remark}

\subsection{Law of large numbers}
With the above definitions at hand, we can formulate the law of large numbers extending \parencite[Proposition 3.6]{forien}. Its proof can be found at the end of \Cref{subsubsec:tightnessstochasticintegrals}.

\begin{restatable}{theorem}{theoremconvergencetolebesgue} \label{theo:convergencetolebesgue}
	Let $ (\delta_N, N \geq 1) $ be a sequence of positive numbers converging to zero as $ N \to \infty$ and such that $ N \eta_N \to \infty $ as $ N \to \infty $. Let $ \boldsymbol{\rho}^N $ be the SLFV started from $\rho_0 \equiv \lambda$. Then, for each fixed $T > 0$,
	\begin{equation*}
		\lim_{N \to \infty} \mathbb{E} \Bigg[ \sup_{t \in [0,T] } d (\boldsymbol{\rho}_t^N, \lambda) \Bigg] = 0,
	\end{equation*}
where $ d $ is the metric defined in \eqref{equationvaguemetric}.
\end{restatable}

\subsection{Central limit theorem} \label{sectioncentrallimittheorem}

Prior to stating the theorem, we need to introduce a few notations.
Let $V_r$ denote the volume of a ball $B(x,r)$ of radius $r > 0$ centered at $x \in \mathbb{R}^{\dimension}$ in $\R^{\dimension}$ and let $V_r(x,y)$ denote the volume of the intersection of $B(x,r)$ and $B(y,r)$  for $x, y \in \R^{\dimension}$.
Let $\alpha$, $\beta$, $\gamma$, $\sigma$ and $\zeta$ be given by Table~\ref{table:parameters}.
Define an operator $\mathcal{D}_\alpha$ acting on twice-continuously differentiable functions $ \phi : \mathbb{R}^{\dimension} \times [0,1] \to \R $ by
\begin{equation} \label{alphastableprocessgenerator}
\mathcal{D}_\alpha \phi (x,k) := \begin{cases} \zeta \int_0^\infty \left( \overline{\phi} (x,k,r) - \phi (x,k)\right) \dfrac{dr}{r^{1+ \alpha}} & \text{ if } \alpha < 2, \\[1em]
\dfrac{\sigma^2}{2} \Delta \phi(x,k) & \text{ if } \alpha = 2,
\end{cases}
\end{equation}
where $\overline{\phi} (x,k,r)$ is defined as the average of $y \mapsto \phi(y,k)$ over the ball $B(x,r)$, i.e.
\begin{equation} \label{singleaverage}
	\overline{\phi} (x, k, r) = \frac{1}{V_r} \int_{B(x,r)} \phi (y,k ) dy,
\end{equation}
and the Laplacian operator $\Delta$ acts only on the spatial coordinates, i.e.
\begin{equation*}
    \Delta \phi(x,k) = \sum_{i=1}^d \frac{\partial^2 \phi}{\partial x_i^2}(x,k).
\end{equation*}
In the case $\alpha < 2$, $\mathcal{D}_\alpha \phi(x,k)$ is equal to
\begin{equation*}
\frac{\zeta}{V_1(d+ \alpha)} \lim\limits_{r \to 0} \int_{\mathbb{R}^{\dimension} \setminus B(x,r)} \frac{\phi (z,k) - \phi (x,k)}{\| z- x \|^{\dimension + \alpha}} dz,
\end{equation*}
which in turn, if multiplied by $ \zeta^{-1} V_1 (\dimension + \alpha) 2^\alpha \pi^{-\dimension/2} \Gamma (\frac{\dimension + \alpha}{2}) \vert \Gamma ( - \frac{\alpha}{2}) \vert^{-1} $, is the fractional Laplacian applied to $x \mapsto \phi(x,k)$ (see \parencite{fractionalLaplacian}). 
Furthermore, define a measure $K_\beta$ on $(\R^{\dimension})^2$ by
\begin{equation} \label{eq:correlation}
K_\beta(dx_1, dx_2) := \begin{cases}
\delta_{x_2}(dx_1) dx_2 & \text{ if } \beta \geq \dimension, \\[1em]
\dfrac{1}{\| x_1 - x_2 \|^{\beta}} dx_1 dx_2 & \text{ if } \beta < \dimension.
\end{cases}
\end{equation}
\def\spacingtabular{8pt}
\setlength{\extrarowheight}{4pt}
\begin{table}
    \centering
    \begin{tabular}{c|c|c}
         & One tail regime & Two tail regime  \\
         [\spacingtabular]
         \hline
         $\alpha$ & $\tail \wedge \frac{\tail}{\power} \wedge 2$ & $\tail_1 \wedge \tail_2 \wedge 2$\\[\spacingtabular]
         \hline
         
         $\beta$ & $\tail + \impact$ & $\tail_2 + \impact_2$ \\[\spacingtabular]
         \hline
         
         $\gamma$ & $\begin{array}{cc}
         \dfrac{u_0^2 V_1^2}{\tail + \impact - \dimension} & \text{ if } \beta > \dimension \\
         u_0^2 V_1^2 & \text{ if } \beta = \dimension \\
         u_0^2 C_{d,\beta}^{(2)} & \text{ if } \beta < \dimension
         \end{array}$ & $\begin{array}{cc}
         \dfrac{u_0^2 V_1^2}{(\tail_1 + \impact_1) (\tail_2 + \impact_2 - \dimension)} & \text{ if } \beta > \dimension \\[0.5em]
         \dfrac{u_0^2 V_1^2}{\tail_1 + \impact_1} & \text{ if } \beta = \dimension \\[0.5em]
         \dfrac{u_0^2 C_{d,\beta}^{(2)}}{\tail_1 + \impact_1} & \text{ if } \beta < \dimension
         \end{array}$ \\[38pt]
         \hline
         
        $\sigma^2$ & $\dfrac{u_0 V_1}{\dimension + 2} \left( \dfrac{1}{\tail - 2 \power} + \dfrac{1}{\tail - 2} \right)$ & $\dfrac{u_0 V_1}{\dimension + 2} \left( \dfrac{1}{(\tail_1-2)\tail_2} + \dfrac{1}{\tail_1 (\tail_2-2)} \right)$ \\[\spacingtabular]
         \hline
         
         $\zeta$ & $\begin{array}{cc} \dfrac{u_0 V_1}{1 \lor \power} & \text{ if } \power \neq 1 \\[0.5em]
         \dfrac{C_{d,\alpha}^{(1)}}{V_1 (\dimension + \alpha)} & \text{ if } \power = 1 \end{array}$ & $\begin{array}{cc}
         \dfrac{u_0 V_1}{\tail_1 \vee \tail_2} & \text{ if } \tail_1 \neq \tail_2 \\
         u_0 V_1 \left( \dfrac{1}{\tail_1} + \dfrac{1}{\tail_2} \right) & \text{ if } \tail_1 = \tail_2
         \end{array}$ 
    \end{tabular}
    \caption[short caption]{Expressions for the different parameters of the limiting fluctuations. The constants $C^{(1)}_{d,\alpha}$ and $C^{(2)}_{d,\beta}$ are given by
    \[
    C^{(1)}_{d,\alpha} := \int_0^\infty \frac{V_r(0,e_1)}{V_r^2} \frac{dr}{r^{1+\alpha}}, \hspace{1cm} C^{(2)}_{d,\beta} := \int_0^\infty V_r(0,e_1) \frac{dr}{r^{1+\dimension+\beta}},
    \]
    where $e_1$ denotes an arbitrary point of the unit sphere in $\mathbb{R}^{\dimension} $. } 
\label{table:parameters}
\end{table}

Our next result shows the convergence of the fluctuations \eqref{scalingoffluctuations} to the solution of a stochastic partial differential equation (SPDE). 
A majority of this work will be dedicated to proving this result. 
The SPDE is of a common form, and mirrors the corresponding result of \parencite[Theorem 2.1, Theorem 2.3]{forien}. However, the one- and two tail regimes allow for a plethora of new parameter combinations specified in \Cref{table:parameters}. We denote the space of tempered distributions acting on $\mathcal{S}(\mathbb{R}^{\dimension} \times [0,1])$ by $\mathcal{S}' (\R^{\dimension} \times [0,1])$. The proof of the result is given in \Cref{subsec:proof-clt} and differences to \cite{forien} are detailed in \Cref{subsec:proof_differences}.
\begin{theorem} \label{centrallimittheorem}
Let $ (Z_t^N, t \geq 0) $ be the fluctuations of an SLFV process with mutation as defined in \eqref{scalingoffluctuations}.
Suppose that the assumptions of \cref{theo:convergencetolebesgue} are satisfied and that, in addition,
\begin{align}
    &\lim_{N \to \infty} \sqrt{\frac{\eta_N}{N}} \delta_N^{{\dimension} \wedge {\impact}} = 0, \quad \text{ in the one-tail regime,}  \label{assumption_jumps_one_tail}\\
    &\lim_{N \to \infty} \sqrt{\frac{\eta_N}{N}} \delta_N^{{\dimension} \wedge {\impact}_2} = 0, \quad \text{ in the two-tail regime.}  \label{assumption_jumps_two_tail}
\end{align}
Then the sequence $ ( Z_t^N, t \geq 0 ) , N = 1,2,... $ converges in distribution in $ D ( \mathbb{R}_+, \mathcal{S}' (\mathbb{R}^{\dimension} \times [0,1])) $ to a process $ (Z_t, t \geq 0) $ solving the mild form of the equation
\begin{equation} \label{SPDE}
\left\lbrace
\begin{aligned}
& d Z_t = [\mathcal{D}_\alpha Z_t - \mu Z_t] dt + \sqrt{\gamma}\, dW(t), \\
&Z_0 = 0.
\end{aligned}
\right.
\end{equation}
Here, $ (W(t), t \geq 0) $ denotes a Wiener process on $ \mathbb{R}^{\dimension} \times [0,1] $, uncorrelated in time, with covariation measure
\begin{equation*}
\mathcal{Q}(dx_1 dk_1 dx_2 dk_2) = K_\beta (dx_1, dx_2) \big(dk_1 \delta_{k_1} (dk_2) - dk_1 dk_2 \big).
\end{equation*}
\end{theorem}

\begin{remark}
The results of \parencite{forien} appear in \Cref{centrallimittheorem} as a special case of the one-tail regime when $\power = 1, \impact = 0$ and $\tail + \impact < \dimension$. 
\end{remark}

\begin{remark}
The assumptions \eqref{assumption_jumps_one_tail} and \eqref{assumption_jumps_two_tail} on the convergence speed of the sequence $(\delta_N, N \geq 1)$ follow in many situations already from the assumption that $\delta_N \to 0$ as $ N \to \infty$, e.g. in the case $\impact \geq \dimension$ and $\tail + \impact \neq \dimension$ (so that there the logarithmic term does not occur in $\eta_N$). The additional assumptions \eqref{assumption_jumps_one_tail} and \eqref{assumption_jumps_two_tail} are only used in the proof of \Cref{lem:martingalejumps}.
\end{remark}

\subsection{Wright-Mal\'ecot Formulae} \label{subsec:WMF}

From the central limit theorem, we can derive expressions for the probability of identity for two individuals sampled from a population evolving according to our SLFV model (see \Cref{subsec:WMFproof}). The connection was established in \parencite[Theorem 2.4]{forien}, and the main idea can be described as follows:
Let $\phi, \psi : \mathbb{R}^d \rightarrow \mathbb{R}$ be two smooth and compactly supported probability densities. Recall from the introduction that $ P_t^N (\phi, \psi) $ denotes the probability that two individuals sampled independently according to $ \phi, \psi $ from the $ N $-th rescaled SLFV model $ (\boldsymbol{\rho}_t^N, t \geq 0) $ are of the same type. 
We can write
\begin{equation*}
P_t^N (\phi , \psi) = \mathbb{E} \big[ \blangle \boldsymbol{\rho}_t^N \otimes \boldsymbol{\rho}_t^N , (\phi \otimes \psi) \mathds{1}_{\Delta} \brangle \big]
\end{equation*}
where $ \mathds{1}_\Delta : [0,1]^2 \rightarrow \mathbb{R} $ is defined as the indicator of the diagonal $\Delta = \lbrace (k_1, k_2) \in [0,1]^2 : k_1 = k_2 \rbrace$.
It was observed in \parencite[p.33f]{forien} that
\begin{equation*}
P_t^N (\phi, \psi) = (N \eta_N)^{-1} \mathbb{E} \big[ \blangle Z_t^N \otimes Z_t^N , (\phi \otimes \psi ) \mathds{1}_\Delta \brangle \big],
\end{equation*}
since $ Z_t^N := \sqrt{N \eta_N } (\boldsymbol{\rho}_t^N - \lambda) $ and the Lebesgue measure $ \lambda $ does not charge the diagonal $\Delta$. 
This implies
\begin{equation*}
	\lim_{N \to \infty} N \eta_N P_t^N (\phi , \psi) = \mathbb{E} \big[ \blangle Z_t \otimes Z_t , (\phi \otimes \psi) \mathds{1}_\Delta \brangle \big],
\end{equation*}
and connects $ P_t^N(\phi, \psi) $ to the fluctuations enabling us to recover Wright-Mal\'ecot formulae.

\begin{theorem}
\label{wmf}
Let $ G_t^\alpha $ be the fundamental solution associated to $\partial_t - \mathcal{D}_\alpha $. Then the Wright-Mal\'ecot formula is given by
\begin{equation} \label{equationwmftestfunctions}
\lim_{t \to \infty} \lim_{N \to \infty} N \eta_N P_t^N (\phi, \psi) =  \int_{(\mathbb{R}^{\dimension})^2} F_{\dimension, \alpha, \beta} \Big( \vert x - y  \vert \Big) \phi (x) \psi (y) dx dy,
\end{equation}
where the function $ F_{\dimension, \alpha, \beta} : \mathbb{R}_+ \rightarrow \mathbb{R}_+ $ is defined as
\begin{equation} \label{eq:functionofwmf}
F_{\dimension, \alpha, \beta} \Big(\vert x - y \vert \Big) := \gamma \int_0^\infty \int_{(\mathbb{R}^{\dimension})^2} e^{- 2 \mu t} G_t^\alpha (x - z_1) G_t^\alpha (y - z_2) K_\beta (dz_1, dz_2) dt.
\end{equation}
\end{theorem}
\Cref{subsec:WMFproof} is devoted to the proof of this theorem.
\subsection{Interpretation} \label{subsec:interpretation}

\begin{figure}
    \centering
    \includegraphics[width=\textwidth]{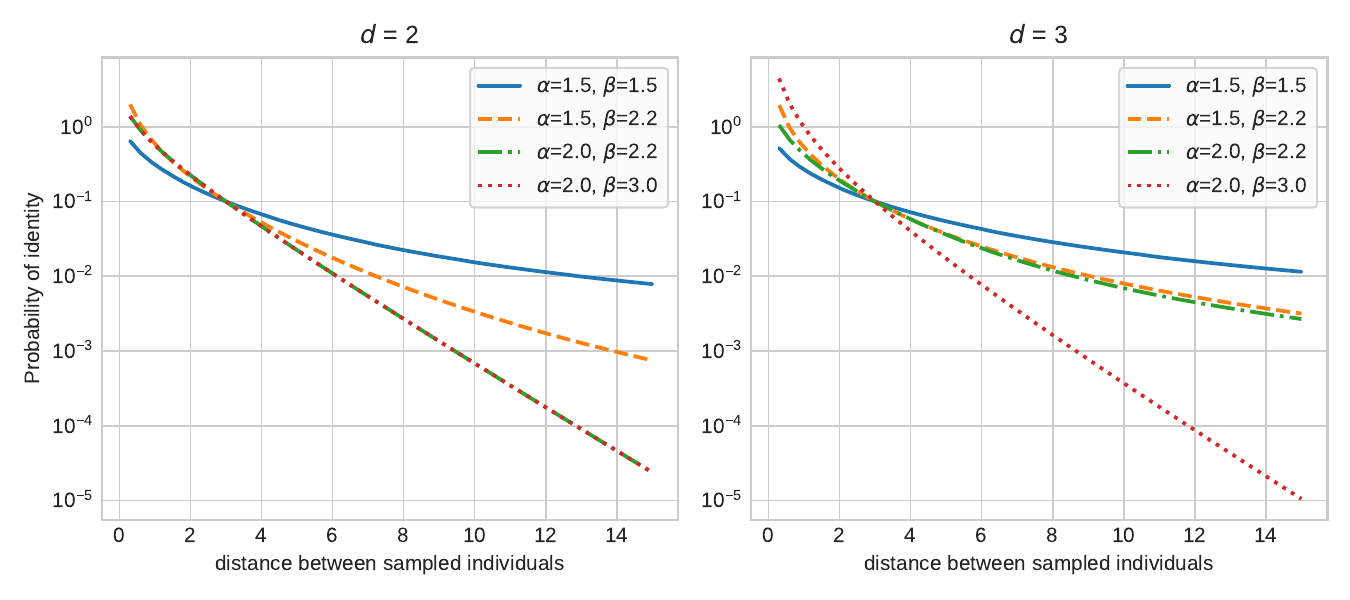}
    \caption{Graph of the function $r \mapsto F_{\dimension,\alpha, \beta}(r)$ defined in \eqref{eq:functionofwmf} for several choices of $\alpha$ and $\beta$. The function is computed after setting $\zeta = 1$, $\sigma^2= 1$, $\mu = 0.2$, and $\gamma$ is chosen in order to normalize the functions at $r = 3$. We do not specify $u_0, \tail, \power, \impact$ and as such the functions are only defined up to multiplicative constants. We thus normalize by our choice of $\gamma$; the aim of the figure is to show how alpha and beta affect the way in which the probability of identity decreases as a function of the distance between the samples.
    On the left, $\dimension = 2$, and four curves are plotted: the blue curve corresponds to long range dispersal ($\alpha = 1.5 < 2$) and long range coalescence ($\beta =1.5 < d$), the orange dashed curve corresponds to long range dispersal ($\alpha < 2$) with short range coalescence ($\beta = 2.2 \geq \dimension$), the green and red curves correspond to short range dispersal ($\alpha = 2$) and short range coalescence ($\beta \geq \dimension$). The last two curves are superposed because of the normalization. On the right plot, four curves with the same parameters are plotted in dimension $\dimension = 3$: the blue and orange curves correspond to long range dispersal ($\alpha = 1.5 < 2$) and long range coalescence ($\beta < \dimension$), the green curve corresponds to short range dispersal ($\alpha = 2$) and long range coalescence ($\beta = 2.2 < \dimension$), and the red curve corresponds to short range dispersal ($\alpha = 2$) and short range coalescence ($\beta = 3 \geq \dimension$). It can be noted that having either long range dispersal or long range coalescence significantly slows down the decrease of the probability of identity with geographic distance, and that this effect is more pronounced for small values of $\alpha$ and $\beta$.}
    \label{fig:wright_malecot}
\end{figure}

For a given set of parameters of the one- or two-tail regime, \Cref{table:parameters} provides the values of $\alpha$ and $\beta$, from which we deduce the form of $\mathcal{D}_\alpha $, $K_\beta (dx_1, dx_2)$ and the corresponding Wright-Mal\'ecot formula \eqref{equationwmftestfunctions}. 
The generator $\mathcal{D}_\alpha, \alpha \in (0, 2]$ reflects either a stable process or Brownian behaviour of ancestral lineages. The correlations $K_\beta(dx_1, dx_2)$ are either of the form $\delta_{x_1} (x_2)$, indicating local coalescence, i.e. asymptotically ancestral lineages need to meet to coalesce, or are of the form
\begin{equation*}
K_\beta (dx_1,dx_2) = \frac{1}{\Vert x_1 - x_2 \Vert^\beta} dx_1 dx_2.
\end{equation*}
This quantity reflects the rate of long range coalescence after noting that
\begin{align*}
\gamma K_\beta (x_1, x_2) = u_0^2 C_{d, \beta}^{(2)} K_\beta (x_1,x_2) 
&= u_0^2 \int_{\vert x_1 - x_2 \vert/2}^\infty V_r (x_1, x_2) \frac{dr}{r^{1 + \dimension + \tail + \impact}} dx_1 dx_2.
\end{align*}
Ancestral lineages at positions $x_1$ and $x_2$ coalesce precisely if they are affected by a reproduction event of radius bigger than $\vert x_1 - x_2 \vert/2$, whose centre falls into the intersection of $B(x_1,r)$ and $B(x_2,r)$ of volume $V_r (x_1, x_2)$. 
This happens at a rate proportional to the squared impact $u_0^2 /r^{2 \impact}$, as both ancestral lineages need to be affected by the event. 
Integrating with respect to $r$ we arrive at the above expression. Our results reveal two main new phenomena.

\paragraph{Decorrelation of ancestral motion and coalescence behaviour}
		Most notably, a phase transition occurs at $ \beta = \dimension $. If $ a + c =: \beta  \geq \dimension $ or $\tail_2 + \impact_2 =: \beta \geq \dimension$ no spatial correlations persist and coalescence happens locally, while if $\beta < \dimension$ the spatial correlations are of non-trivial form $K_\beta (dx_1, dx_2) = \Vert x_1 - x_2 \Vert^{- \beta} dx_1 dx_2$ and coalescence happens at long range. 
  Crucially, the correlations/coalescence behaviour $ K_\beta (x_1, x_2) $ and the parameter $ \alpha $ of stable ancestral lineages are not as closely linked as suggested by the results of \parencite{forien}. Fractional Brownian lineages can now appear with only local coalescence, e.g. by choosing in the one-tail regime in dimension $\dimension = 2$,
		\begin{equation*}
		   \tail = 1.5, \enspace \power = 1, \enspace \impact = 0.7.
		\end{equation*}
		By \Cref{centrallimittheorem}, ancestral lineages are fractional ($\alpha = \tail = 1.5$) and coalescence happens locally as $\beta = \tail + \impact = 2.2 > 2 = \dimension$ and therefore $K(x_1, x_2) = u_0^2 V_1^2 \delta_{x_1} (x_2)$. Even more surprisingly, in dimension $\dimension \geq 3$ we can obtain Brownian lineages, while still maintaining long range coalescence. For example, we can choose in the one-tail regime in $\dimension = 3$,
		\begin{equation*}
		    \tail = 2, \enspace \power = 0.5, \enspace \impact = 0.2.
		\end{equation*}
		Again, by \Cref{centrallimittheorem}, ancestral lineages are Brownian $\alpha = \frac{2}{1 \lor 0.5} \land 2 = 2$ and because $\beta = \tail + \impact = 2.2 < 3 = \dimension$ the spatial correlations/coalescence rate is 
		\[
		K(dx_1, dx_2) = u_0^2 C_{d, 2}^{(2)} \frac{ dx_2 dx_2}{\Vert x_1 - x_2 \Vert^\beta}.
		\]
		Both of the above cases can be found in \Cref{fig:wright_malecot}.
\paragraph{Asymmetric influence of $r_1$ and $r_2$} 
    The operator $\mathcal{D}_\alpha$ is a fractional or standard Laplacian depending on the value of $\alpha$. We can note that $\alpha$ is governed by the tails of the two marginal distributions of $r_1$ and $r_2$.
    More precisely, the radius with the heaviest tail dominates and determines the dispersal regime of the ancestral lineages (with a correction by a factor $r_2^d$ in the tail of $r_2$ coming from the rate at which a reproduction event affects a given spatial location, see Section~\ref{subsec:duality} below).

    On the other hand, the coalescence regime (i.e. whether $K_\beta$ is concentrated on the diagonal or not) is only determined by the tail of the replacement radius $r_2$.
    Long range coalescence can only originate from reproduction events with a very large replacement radius.
    This is why $\power$ does not appear in the formula for $\beta$ in the one-tail regime, and neither do $\tail_1$ and $\impact_1$ in the two-tail regime.

    The parameter $\impact$ (or $\impact_1$ and $\impact_2$ in the two-tail regime) describes how the impact of a reproduction event depends on its size.
    For instance, if $\impact = 0$ in the one-tail regime, the impact is constant, and the mass of offspring produced during an event is proportional to $r_2^d$, while if $\impact = \dimension$, the mass of offspring produced is independent of the size of the reproduction event.
    The parametrisation of the two regimes in the definition of $\nu_N$ (see \eqref{one_tail_regime} and \eqref{two_tail_regime}) is such that these parameters do not affect the value of $\alpha$, but they do impact $\beta$.
    Hence the parameter $\impact$ (or $\impact_2$) is another way of decoupling the dispersal from the coalescence regime.

\section{Proof of the main results} \label{sec:proof_outline}
At the heart of the present work lies \Cref{centrallimittheorem}, a central limit theorem for the SLFV process. First, we represent the SLFV process as a semimartingale driven by worthy martingale measures in \Cref{subsec:semimartingaledecomposition}. Worthy martingale measures were used by \parencite{walsh_introduction_1986} to generalize classical stochastic integrals to include space. In turn, this theory allows us to write the sequence of SLFV processes as a sequence of stochastic integrals. Second, in \Cref{subsec:proof-clt}, we proceed by proving tightness of the sequence and uniqueness of limit points, which proves \Cref{centrallimittheorem}. The Wright-Mal\'ecot formulae will follow almost immediately from the corresponding central limit theorems (\Cref{subsec:WMFproof}). \par
The overall strategy of the proofs follows \cite{forien}. Non-trivial adaptations to cover the new regimes are necessary in several intermediate lemmas, which we introduce in the overall context in this section.
We will summarize these adaptations in \Cref{subsec:proof_differences}, before proving the lemmas in \Cref{sec:proof_lemmas}.

\subsection{Semimartingale decomposition} \label{subsec:semimartingaledecomposition}
We characterize the SLFV via a martingale problem adapting the calculations of \parencite[Proposition 3.3]{forien} to include two different radii and varying impact. Let \[E^p := L^p (\mathbb{R}^{\dimension} \times [0,1]) \cap L^\infty (\mathbb{R}^{\dimension} \times [0,1]) \] be the space of uniformly bounded $ L^p $-functions. For $\phi \in E^p$, the averaged function $\overline{\phi}$ was defined in \eqref{singleaverage}. Accordingly, the double averaged function is for radii $ r_1, r_2 > 0 $ is defined as
\begin{equation*}
	\overline{\overline{\phi}} (x, k, r_1, r_2) := \frac{1}{V_{r_1} V_{r_2}} \int_{B(x,r_1)} \int_{B(y,r_2 )} \phi (z,k) dz dy.
\end{equation*}
For $ r_1, r_2 > 0 $ and $\rho \in \Xi$, we define $ [\rho]_{r_1, r_2} : \mathbb{R}^{\dimension} \times \mathbb{R}^{\dimension} \rightarrow \mathcal{M} ([0,1]) \cup \{ 0 \} $ by
\begin{equation*}
\begin{aligned}
	&[\rho]_{r_1, r_2} (x_1,x_2, dk)\\
	 &= \begin{dcases}
		\frac{1}{V_{r_2} (x_1, x_2)} \int_{B(x_1,r_2) \cap B(x_2,r_2)} \frac{1}{V_{r_1}} \int_{B(y,r_1)} \rho (z, dk) dz dy & \text{ if } \vert x_1 - x_2 \vert < 2 r_2,\\
		0 & \text{ otherwise;}
	\end{dcases}
\end{aligned}
\end{equation*}
where once more $ V_r(x_1, x_2) := \Vol (B(x_1,r) \cap B(x_2, r))$ for $ r > 0 $. Furthermore, we define the map 
\begin{equation} \label{eq:Gamma}
	\begin{aligned}
		&\Gamma^{\nu} (\rho) (x_1, x_2 , dk_1, dk_2)\\
		&=\int_{(0,1] \times (0, \infty)^2} u^2 V_{r_2} (x_1,x_2) \Big[ [\rho]_{r_1, r_2} (x_1,x_2, dk_1) \delta_{k_1} (dk_2) - \rho (x_1, dk_1) [\rho]_{r_1, r_2} (x_1, x_2, dk_2) \\ 
		& \hspace{1cm} - [\rho]_{r_1, r_2} (x_1, x_2, dk_1) \rho(x_2, dk_2) + \rho (x_1 ,dk_1) \rho (x_2, dk_2) \Big] \nu (du, dr_1, dr_2).
	\end{aligned}
\end{equation}
The symbols $[\rho]_{r_1,r_2}$ and $\Gamma^\nu$ are used to express the predictable variation process in the semimartingale below. More precisely, the integral $\frac{1}{V_{r_1}} \int_{B(y,r_1)} \rho(z,dk) dz$ describes the distribution of parental types for an event at center $y \in \mathbb{R}^d$ and of parental search radius $r_1$. In turn, $[\rho]_{r_1, r_2}$ is the distribution of the parental type in an event given that both $x_1$ and $x_2$ are in the replacement ball of the event and $\Gamma^\nu$ measures the rate at which the covariance between $\rho_t (x_1, dk)$ and $\rho_t (x_2,dk)$ increases over time due to these reproduction events.

\begin{lemma}
\label{martingaleproblem}
	Let $ (\rho_t, t \geq 0) $ be the SLFV with mutations. For any test function $ \phi \in E^1 $,
	\begin{multline} \label{equationinitialmartingaleproblem}
			\mathcal{M}_t(\phi) := \blangle \rho_t , \phi\brangle -\blangle \rho_0 , \phi\brangle - \int_0^t \Bigg\{ \mu\blangle \lambda - \rho_s , \phi\brangle \\
			+ \int_{(0,1] \times (0, \infty)^2} u V_{r_2} \blangle \rho_s , \overline{\overline{\phi}} (\cdot, r_1, r_2) - \phi\brangle \nu (du, dr_1, dr_2) \Bigg\} ds
	\end{multline}
	is a square-integrable martingale adapted to the natural filtration $ (\mathcal{F}_t,t \geq 0) $ with predictable quadratic variation process
	\begin{equation*}
		t \mapsto \int_0^t\blangle \Gamma^{ \nu} (\rho_s) , \phi \otimes \phi\brangle ds.
	\end{equation*}
\end{lemma}

\begin{proof}
The proof is a straightforward adaptation of the calculations of \parencite[p.18f]{forien} to reproduction events with two radii $r_1$ and $r_2$. The finite variation term is a consequence of 
\begin{align*}
\lim_{\delta t \to 0} &\frac{1}{\delta t} \mathbb{E} \left[ \langle \rho_{t-\delta t} , \phi \rangle - \langle \rho_t , \phi \rangle \vert \rho_t = \rho \right] \\
&= \int_{\mathbb{R}^{\dimension}} \int_{(0,1] \times (0, \infty)^2} \frac{1}{V_{r_1}} \int_{B(x, r_1) \times [0,1]} \langle u \mathds{1}_{B(x, r_2)} (\delta_{k_0} - \rho ), \phi \rangle \rho(y,dk_0) dy \nu (du, dr_1, dr_2) dx \\
&\hspace{10cm} + \mu \langle \lambda - \rho , \phi \rangle
\end{align*}
and the exchange of the order of integration. The four terms of \eqref{eq:Gamma} originate from a similar calculation applied to
\begin{align*}
\lim_{\delta t \to 0} &\frac{1}{\delta t} \mathbb{E} \left[ (\langle \rho_{t-\delta t} , \phi \rangle - \langle \rho_t , \phi \rangle )^2 \vert \rho_t = \rho \right] \\
&= \int_{\mathbb{R}^{\dimension}} \int_{(0,1] \times (0, \infty)^2} \frac{1}{V_{r_1}} \int_{B(x, r_1) \times [0,1]} \langle u \mathds{1}_{B(x, r_2)} (\delta_{k_0} - \rho ), \phi \rangle^2 \rho(y,dk_0) dy \nu (du, dr_1, dr_2) dx. \qedhere
\end{align*} 
\end{proof}

As a next step, we will consider the martingale problem under the rescaling introduced in \Cref{subsec:choiceofscaling}. For this purpose, we define the rescaled measure $ \nu_N^\alpha $ 
by
\begin{equation*}
\begin{aligned}
	&\int_{(0,N] \times (0, \infty)^2} f(u, r_1, r_2) \nu_N^\alpha (du, dr_1, dr_2) \\
    &:= \int_{(0,1]\times(0, \infty)^2} f(N u, \delta_N r_1, \delta_N r_2) \delta_N^{-( \alpha +d)} \nu_N (du, dr_1, dr_2),
    \end{aligned}
\end{equation*}
and the operator $ \mathcal{L}^N $ acting on $E^1$ by
\begin{equation} \label{equationrescaledgenerator}
\begin{aligned}
	\mathcal{L}^N \phi (x,k) 
    &:= \int_{(0,N]\times(0, \infty)^2} u V_{r_2} \left(\overline{\overline{\phi}} (x, k, r_1, r_2) - \phi (x,k) \right) \nu_N^\alpha (du, dr_1, dr_2). \notag
\end{aligned}
\end{equation}
The idea to determine the rescaled martingale problem is to use the equation
\begin{equation*}
	\blangle \boldsymbol{\rho}_t^N, \phi\brangle =\blangle \rho_{\frac{N t}{\delta^\alpha_N}} , \phi_N\brangle,
\end{equation*}
where $ \phi_N (x,k) := \delta_N^{\dimension} \phi ( \delta_N x, k) $, i.e. rewrite $ \blangle \boldsymbol{\rho}_t^N , \phi \brangle $ as a version of $ (\rho_t, t \geq 0) $ with accelerated time integrated against test functions with the correct spatial scaling.
\begin{lemma}[extends {\parencite[p.19f]{forien}}] \label{rescaledmartingaleproblem}
Let $(\boldsymbol{\rho}_t^N, t \geq 0) $ be the rescaled SLFV model as described in \Cref{subsec:choiceofscaling}. Then for any $ \phi \in E^1 $
\begin{equation*}
	\langle \boldsymbol{\rho}_t^N , \phi\brangle =\blangle \boldsymbol{\rho}_0, \phi\brangle + \int_0^t \Bigg\{ \mu\blangle \lambda - \boldsymbol{\rho}_s^N, \phi\brangle +\blangle \boldsymbol{\rho}_s^N, \mathcal{L}^N \phi \brangle \Bigg\} ds + \boldsymbol{\mathcal{M}}_t^N (\phi),
\end{equation*}
where $ \boldsymbol{\mathcal{M}}_t^N (\phi) $ is a square-integrable martingale with predictable quadratic variation given by
\begin{equation*}
	\Big\langle \boldsymbol{\mathcal{M}}_t^N(\phi) \Big\rangle_t 
	= \frac{1}{N} \int_0^t \Big\langle \Gamma^{\nu_N^\alpha} (\boldsymbol{\rho}_s^N) , \phi \otimes \phi \Big\rangle ds.
\end{equation*}
\end{lemma}
Hereafter, we apply the scaling of the fluctuations
\begin{equation}
Z_t^N := \sqrt{N \eta_N} (\boldsymbol{\rho}_t^N - \lambda )
\end{equation}
to the martingale problem. The resulting martingale problem falls within the martingale measures framework developed by Walsh in \cite{walsh_introduction_1986}. Martingale measures generalise the notation of martingales and stochastic integrals to involve a spatial component (\Cref{def:martingalemeasure}). Additionally, for a suitable class called worthy (see \Cref{def:worthy}), powerful theorems about the convergence of stochastic integrals exist (see \Cref{theo:mainconvergencetheorem}, \Cref{theo:convergenceoffinitedimensionaldistributions}). In this framework, we can rewrite the fluctuations as a sequence of stochastic integrals, and applications of the convergence theorems will lead to \Cref{centrallimittheorem}.
For functions $f: [0,T] \times \mathbb{R}^{\dimension} \times [0,1] \rightarrow \mathbb{R}$ and a worthy martingale measure $M$, we will use the abbreviation
\begin{equation*}
M_t (f) := \int_{[0,t] \times \mathbb{R}^{\dimension} \times [0,1]} f (s,x,k) M (ds, dx, dk).
\end{equation*}
Rewriting the semimartingale decomposition in \Cref{rescaledmartingaleproblem} with this new terminology and applying the scaling of the fluctuations \eqref{scalingoffluctuations} results in the following representation.
\begin{proposition}[extends {\parencite[Proposition 3.1]{forien}}] \label{prop:fluctuationsrescaled}
Assume that $\rho_0^N = \lambda$. For any $ \phi \in E^1 $, the fluctuations $Z_t^N := \sqrt{N \eta_N} (\lambda - \boldsymbol{\rho}_t^N )$ can be written as
	\begin{equation*}
		\big\langle Z_t^N, \phi\brangle = \int_0^t\blangle Z_s^N, \mathcal{L}^N ( \phi) - \mu \phi\brangle ds + M_t^N ( \phi),
	\end{equation*}
where $ M^N $ is a worthy martingale measure on $\R^{\dimension} \times [0,1]$ with
	\begin{equation*}
		\big\langle M^N (\phi)\brangle_t = \eta_N \int_0^t\blangle \Gamma^{\nu_N^\alpha} (\boldsymbol{\rho}_s^N), \phi \otimes \phi\brangle ds.
	\end{equation*}
\end{proposition}
\begin{proof}
The calculations and argumentation of \parencite[p.19f]{forien} readily carry over to the new SLFV models. A dominating measure for the martingale measures is given by
\begin{equation} \label{dominatingmeasuresdefinition}
D_N (dt, dx_1, dk_1, dx_2, dk_2 ) := \eta_N\vert \Gamma \vert^{\nu_\alpha^N} (\boldsymbol{\rho}_{t^-}^N) (x_1, x_2, dk_1, dk_2) dx_1 dx_2 dt,
\end{equation}
where we define
\begin{equation*}
	\begin{aligned}
		&\vert \Gamma \vert^{\nu} (\rho) (x_1, x_2 , dk_1, dk_2)\\
		&:=\int_{(0,1]\times (0, \infty)^2} u^2 V_{r_2} (x_1,x_2) \Big[ [\rho]_{r_1, r_2} (x_1,x_2, dk_1) \delta_{k_1} (dk_2) + \rho (x_1, dk_1) [\rho]_{r_1, r_2} (x_1, x_2, dk_2) \\ 
		& \hspace{1cm} + [\rho]_{r_1, r_2} (x_1, x_2, dk_1) \rho(x_2, dk_2) + \rho (x_1 ,dk_1) \rho (x_2, dk_2) \Big] \nu (du, dr_1, dr_2). 
	\end{aligned}
\end{equation*}
\end{proof}
The operator
\begin{equation*}
	\mathcal{L}^N \phi (x,k) = \int_{(0,1] \times (0, \infty)^2} u V_{r_2} \Big(\overline{\overline{\phi}} (x, k , r_1 , r_2) - \phi (x,k) \Big) \nu_N^\alpha (du, dr_1, dr_2)
\end{equation*}
can be interpreted as the generator of a continuous-time jump process $(X_t^N, t \geq 0)$ describing the spatial motion of an ancestral lineage in the process with behaviour as follows. 
At rate
\begin{equation*}
\int_{(0,1]\times(0, \infty)^2} u V_{r_2} \nu_N^\alpha (du, dr_1, dr_2)
\end{equation*}
the process jumps, if located at position $x$, to position $x + (Y_1 R_1 + Y_2 R_2)$, where $ (R_1, R_2)$ is a random variable taking values in $(0,\infty)^2$ with distribution
\begin{equation*}
\frac{\int_{(0,1]} u V_{R_2} \nu_N^\alpha (du, dR_1, dR_2)}{\int_{(0,1]\times(0, \infty)^2} u V_{r_2} \nu_N^\alpha (du, dr_1, dr_2)},
\end{equation*}
and $Y_1, Y_2$ are two independent uniform random variables on $B(0,1)$.
The operator $\mathcal{L}^N - \mu$ generates the strongly-continuous semigroup $P_t^N$ acting on $E^1$, defined as
\begin{equation} \label{eq:definitionofthesemigroups}
P_t^N \phi (x,k) := e^{- \mu t} \mathbb{E}_x \left[ \phi ( X_t^N, k ) \right].
\end{equation}
Now, the semimartingale decomposition in \Cref{prop:fluctuationsrescaled} enables us to represent $\big\langle Z_t^N, \phi \big\rangle$ as a stochastic integral with respect to a martingale measure.
Indeed, by \parencite[Theorem 5.1]{walsh_introduction_1986}, we can write
\begin{equation} \label{eq:stochasticintegralrepresentation}
\langle Z_t^N, \phi \rangle = \int_{[0,t] \times \mathbb{R}^{\dimension} \times [0,1]} P_{t-s}^N \phi (x,k) M^N (ds, dx, dk).
\end{equation}
This representation, together with a convergence theorem which can be found in \cite{walsh_introduction_1986}, will play a crucial role in the proof of the central limit theorem.

\subsection{Central limit theorem} \label{subsec:proof-clt}

In this subsection, we prove the convergence of $(Z_N, N \geq 1)$, i.e. the central limit theorem for our SLFV models, following the overall strategy of \parencite[Section 3]{forien}. 
Classically, tightness of the sequence will show that any subsequence admits limit points, and we shall conclude by proving uniqueness of the finite-dimensional distributions of these limit points. Afterwards, in \Cref{subsec:proof_differences} we highlight the main points where additional arguments are needed compared to \cite{forien}. In turn, the proofs of these main points are detailed in \Cref{sec:proof_lemmas}. 

\subsubsection{Tightness of the stochastic integrals} \label{subsubsec:tightnessstochasticintegrals}
Here, we will apply Mitoma's theorem:
\begin{theorem}[{\parencite[Theorem 6.13]{walsh_introduction_1986}}] \label{theoremmitoma}
Let $ X = ((X_t^N)_{t \in [0,1]}, N \geq 1)$ be a sequence of stochastic processes in $ D ([0,1 ], \mathcal{S}'(\mathbb{R}^{\dimension})) $. Then $ X $ is tight if and only if for each Schwartz function $ \phi \in \mathcal{S} (\mathbb{R}^{\dimension}) $, the evaluated sequence $ ((X_t^N (\phi))_{t \in [0,1]}, N \geq 1)$ is tight in $ D([0,1], \mathbb{R}) $.
\end{theorem}
This allows us to reduce tightness of $(Z^N, N \geq 0)$ to tightness of $\big( \big\langle Z_t^N, \phi \big\rangle_{t \geq 0} , N \geq 1 \big) \subset D (\mathbb{R}_+, \mathbb{R})$ for $\phi \in \mathcal{S} (\mathbb{R}^{\dimension} \times [0,1])$. \Cref{theo:mainconvergencetheorem} will show tightness of the latter sequences, when written in the stochastic integral form of \eqref{eq:stochasticintegralrepresentation}. The vague convergence $\boldsymbol{\rho}_t^N \to \lambda$ will be an additional byproduct of a bound in this theorem and we give a short proof below.
To apply \Cref{theo:mainconvergencetheorem} to the sequence of stochastic integrals \eqref{eq:stochasticintegralrepresentation} we need to verify three conditions. The first one requires the dominating measures of the integrators $M^N$ to be uniformly bounded, whereas the other two conditions describe sufficient properties of the integrands $P_{t-s}^N \phi$.

\begin{restatable}[extends {\parencite[Lemma 3.4]{forien}}]{Lemma}{boundonmeasures} \label{lem:boundondominatingmeasures}
	Condition 1 of \Cref{theo:mainconvergencetheorem}: Let $ D^N $ be the dominating measure from \eqref{dominatingmeasuresdefinition}. Then in the one-tail and two-tail regimes, we can find a constant $ C > 0 $ such that for all $ {N \geq 1}, 0 \leq s  \leq t $ and $ \phi , \psi \in \mathcal{S} ( \mathbb{R}^{\dimension} \times [0,1]) $,
	\begin{equation} \label{boundondominatingmeasures}
		\blangle D^N, \mathds{1}_{[s,t]} \phi \otimes \psi \brangle \leq C \vert t - s \vert \Big( \Vert \phi \Vert_1 \Vert \psi \Vert_1 + \Vert \phi \Vert_2 \Vert \psi \Vert_2\Big).
	\end{equation}
\end{restatable}
\begin{restatable}[extends {\parencite[Lemma 3.5]{forien}}]{Lemma}{Boundonintegrands}
\label{lem:semigroupproperties}
	Conditions 2 - 3 of \Cref{theo:mainconvergencetheorem} + convergence of the semigroups: Let $P_t^N$ and $P_t$ be the semigroups associated with $\mathcal{L}^N - \mu$ and $\mathcal{D}_\alpha - \mu$ according to \eqref{eq:definitionofthesemigroups}. 
 \begin{enumerate}
 \item 
If $\phi \in \mathcal{S} ( \mathbb{R}^{\dimension} \times [0,1] )$, then $P_t^N \phi \in \mathcal{S} (\mathbb{R}^{\dimension} \times [0,1])$ and $t \mapsto P_t^N \phi$ is continuous.
     \item 
     Then for all $ t \geq 0, N \geq 1, p \geq 1 $ and $ \phi \in \mathcal{S} ( \mathbb{R}^{\dimension} \times [0,1]) $
	\begin{equation*}
		\big\Vert P_t^N \phi \big\Vert_p \leq e^{- \mu t} \Vert \phi \Vert_p
	\end{equation*}
	and for any $ \kappa \in \mathbb{N}^d $, the derivatives are bounded by
	\begin{equation*}
		\big\Vert \partial_\kappa P_t^N \phi \big\Vert_p \leq e^{- \mu t} \Vert \partial_\kappa \phi \Vert_p.
	\end{equation*}
 \item
	Furthermore, there exists constants $ C > 0$ and $ \theta > 0 $ such that for all $ N \geq 1 $ and $ t \geq 0 $
	\begin{equation*}
		\big\Vert P_t^N \phi - P_t \phi \big\Vert_p \leq C \delta_N^{\theta} t e^{- \mu t} \max_{\vert \kappa \vert \leq 4} \Vert \partial_\kappa \phi \Vert_p.
	\end{equation*}
 \end{enumerate}
\end{restatable}
 The convergence of the semigroups will be required to identify limit points in the next subsection. 
The proof of \Cref{lem:semigroupproperties} will require the boundedness and convergence of the corresponding generators.

\begin{restatable}[extends {\parencite[Prosition A.3]{forien_central_2017}}]{lemma}{convergenceofgenerators} \label{lem:convergenceofgenerators}
Let $ \mathcal{D}_\alpha $ be the generator \eqref{alphastableprocessgenerator} and $ \phi : \mathbb{R}^{\dimension} \rightarrow \mathbb{R} $ be a twice continuously-differentiable function with $ \max_{\vert \kappa \vert \leq 4} \Vert \partial_\kappa \phi \Vert_{p} < \infty $ for $ 1 \leq p \leq \infty $. Then there exists constants $\theta, C, C' > 0$ such that
\begin{equation}  \label{equationboundongenerator}
\big\Vert \mathcal{L}^N \phi \big\Vert_p \leq C \max_{\vert \kappa \vert \leq 2 } \Vert \partial_\kappa \phi \Vert_p
\end{equation}
and 
\begin{equation} \label{equationboundondifferenceofgenerators}
\big\Vert \mathcal{L}^N \phi - \mathcal{D}_\alpha \phi \big\Vert_p \leq C' \delta_N^\theta \max_{\vert \kappa \vert \leq 4 } \Vert \partial_\kappa \phi \Vert_p.
\end{equation}
\end{restatable}

To cover our new regimes, the proofs of \Cref{lem:boundondominatingmeasures} and \Cref{lem:convergenceofgenerators} require significant changes and will therefore be detailed in \Cref{subsec:bound_dominating_measures} and \Cref{subsec:convergence_of_generators}. The remainder of this subsection aligns with the arguments of \cite[p.22f]{forien}.
\begin{proof}[Proof of \Cref{lem:semigroupproperties}]
Recall the definition of $P_t^N \phi$ in \eqref{eq:definitionofthesemigroups}. The random walk is invariant by translation
\begin{equation*}
\mathbb{E}_x \left[ \phi (X_t^N , k ) \right] = \mathbb{E}_0 \left[ \phi (x + X_t^N , k ) \right],
\end{equation*}
which by Fubini's theorem and $\int_{\mathbb{R}^{\dimension}} \phi (x + 
 y) dx = \int_{\mathbb{R}^{\dimension}} \phi (x) dx$ for $ y\in \mathbb{R}^{\dimension}$ implies
\begin{equation*}
    \int_{\mathbb{R}^{\dimension}} \mathbb{E}_x \left[ \phi (X_t^N , k) \right] dx = \int_{\mathbb{R}^{\dimension}} \phi (x,k) dx.
\end{equation*}
This allows us to write, by Jensen's inequality,
\begin{equation*}
\begin{aligned}
\Vert P_t^N \phi \Vert_p &\leq \Bigg( \int_{\mathbb{R}^{\dimension}} \sup_{k \in [0,1]} \big\vert e^{- \mu t} \mathbb{E}_x \left[ \phi (X^N, k ) \right] \big\vert^p dx \Bigg)^{1/p} \\
&\leq e^{- \mu t } \Bigg( \int_{\mathbb{R}^{\dimension}} \mathbb{E}_x \left[ \sup_{k \in [0,1]} \big\vert \phi (X_t^N, k ) \big\vert^p \right] dx \Bigg)^{1/p}\\
&\leq e^{- \mu t } \Vert \phi \Vert_p.
\end{aligned}
\end{equation*}
The bounds on the derivatives and that the $P_t^N \phi$ are Schwartz functions follows from similar calculations.
To prove property 3., let us consider the difference
\begin{equation*}
g_t^N = P_t^N \phi - P_t \phi
\end{equation*}
whose time derivative is
\begin{equation*}
\partial_t g_t^N = ( u \mathcal{L}^{N, \alpha} - \mu ) g_t^N + u (\mathcal{L}^N  - \mathcal{D}_\alpha) P_t \phi.
\end{equation*}
As the difference $g_0^N$ is zero, this yields the representation
\begin{equation*}
g_t^N = u \int_0^t P_{t-s}^{N} ( \mathcal{L}^N - \mathcal{D}_\alpha ) P_s \phi ds.
\end{equation*}
Taking the $L^p$-norm on both sides, we get
\begin{equation*}
\Vert g_t^N \Vert_p \leq u \int_0^t e^{- \mu (t-s)} \Big\Vert (\mathcal{L}^N - \mathcal{D}_\alpha ) P_s \phi \Big\Vert_p ds,
\end{equation*}
and we can apply \eqref{equationboundondifferenceofgenerators} in \Cref{lem:convergenceofgenerators}.
\end{proof}

We are now in a position to apply \Cref{theo:mainconvergencetheorem} to obtain the tightness of the stochastic integrals and the convergence to Lebesgue measure.
\begin{corollary}
Under the assumptions of \Cref{theo:mainconvergencetheorem}, the sequence $ \big(\blangle Z_t^N, \phi \brangle_{t \geq 0} , N \geq 1\big) \subset D(\mathbb{R}_+, \mathbb{R})$ is tight for any $ \phi \in \mathcal{S} (\mathbb{R}^{\dimension} \times [0,1]) $.
\end{corollary}
\begin{proof}
With minor adaptations (see \cite[Proof of Proposition 3.6]{forien}), Lemmas \ref{lem:boundondominatingmeasures} and \ref{lem:semigroupproperties} prove conditions 1 - 3 of \Cref{theo:mainconvergencetheorem}, which yields the tightness of the stochastic integrals \eqref{eq:stochasticintegralrepresentation}.
\end{proof}

\begin{proof}[Proof of \Cref{theo:convergencetolebesgue}]
In addition to the tightness, \Cref{theo:mainconvergencetheorem} shows that there exists a constant $ C^* > 0 $ such that for all $T >0,  N \geq 1 $ and $ \phi \in \mathcal{S}(\mathbb{R}^{\dimension} \times [0,1]) $,
\begin{equation} \label{equationconvergencetolebesgueinequality}
\mathbb{E} \Bigg[ \sup_{t \in [0, T]} \Big\vert \blangle Z_t^N, \phi \brangle \Big\vert^2 \Bigg]^{1/2} \leq C^* \max_{p \in \{ 1, 2\}} \max_{\vert \kappa \vert \leq 2} \Vert \partial_\kappa \phi \Vert_p.
\end{equation}
This allows us to bound 
\begin{equation*}
\mathbb{E} \Bigg[ \sup_{t \in [0, T] } d (\boldsymbol{\rho}_t^N, \lambda) \Bigg] \leq \frac{C^*}{\sqrt{N \eta_N}} \sum_{n = 1}^\infty \frac{1}{2^n} \max_{q \in \{ 1, 2\} } \max_{\vert \kappa \vert \leq 2 } \Vert \partial_\kappa \phi_n \Vert_p,
\end{equation*}
where $ d $ is the vague metric on $ \Xi $ defined in \eqref{equationvaguemetric} and $ \phi_n $ the corresponding sequence of test functions. Using equation \eqref{eq:derivativesbounded}, the right-hand side is bounded by $\frac{C}{\sqrt{N \eta_N}}$ for some $C > 0$, proving the result. 
\end{proof}



\subsubsection{Uniqueness of finite-dimensional distributions}
We need to prove the convergence of the integrands and integrators of \eqref{eq:stochasticintegralrepresentation} to show convergence of the finite-dimensional distributions. In the last subsection, we have seen the $L^p$-convergence of the integrand $P_t^N \rightarrow P_t$.
The convergence in law of the martingale measures $M^N$ will follow from \cref{theo:martingaleconvergence} applied to $M^N (\phi)$ for test functions $\phi \in \mathcal{S} ( \mathbb{R}^{\dimension} \times [0,1])$ with the limit given by
\begin{equation*}
M(\phi) = \sqrt{\gamma} \int_{\mathbb{R}^{\dimension} \times [0,1]} \phi (x,k) W(dxdk)
\end{equation*}
for the Wiener process $W$ of \Cref{centrallimittheorem}.
With both convergence of integrands and integrators at hand, by an adaptation \Cref{theo:convergenceoffinitedimensionaldistributions} of \parencite[Theorem 7.12]{walsh_introduction_1986} this suffices for the convergence of the finite dimensional distributions
\begin{equation*}
\Big( \blangle Z_{t_i}^N, \phi_i \brangle \Big)_{i = 1,...,n} \xrightarrow[N \to \infty]{d} \Bigg( \int_{[0,t_i] \times \mathbb{R}^{\dimension} \times [0,1]} P_{t_i-s} \phi_i (x,k) M (ds, dx, dk) \Bigg)_{i = 1,...,n}
\end{equation*}
 for $ t_1, ..., t_n \geq 0 $ and functions $ \phi_1, ..., \phi_n \in \mathcal{S} (\mathbb{R}^{\dimension} \times [0,1]) $. In turn, by \parencite[Theorem 6.15]{walsh_introduction_1986} we will conclude that the sequence $ ((Z_t^N)_{t \geq 0}, N \geq 1)$ converges in law in $ D(\mathbb{R}_+, \mathcal{S}' (\mathbb{R}^{\dimension} \times [0,1])) $ to the solution $ (Z_t, t \geq 0) $ of \eqref{SPDE} defined by
 \begin{equation}
\blangle Z_t, \phi \brangle = \int_{[0,t] \times \mathbb{R}^{\dimension} \times [0,1]} P_{t-s} \phi (x,k) M (ds, dx, dk).
 \end{equation}
 

We separate the prerequisites of the martingale convergence theorem into two lemmas:
\begin{restatable}[{extends \parencite[Lemma 3.10 Part 1]{forien}}]{Lemma}{lem:martingalejumps} \label{lem:martingalejumps}
	For any Schwartz function $ \phi \in \mathcal{S} ( \mathbb{R}^{\dimension} \times [0,1]) $ the jumps of the martingale measures
	\begin{equation*}
		\sup_{t \geq 0} \Big\vert M_t^N (\phi) - M_{t^-}^N (\phi) \Big\vert \xrightarrow[N \to \infty]{} 0
	\end{equation*}
	vanish almost surely.
\end{restatable}
\begin{restatable}[{extends \parencite[Lemma 3.10 Part 2]{forien}}]{Lemma}{convergenceofquadraticvariation} \label{lem:convergenceofquadraticvariation}
	Let $ \mathcal{Q} $ be defined as
	\begin{equation} \label{eq:equationQ}
		\mathcal{Q}(dx_1 dk_1 dx_2 dk_2) = \gamma K_\beta (dx_1, dx_2) \big(dk_1 \delta_{k_1} (dk_2) - dk_1 dk_2 \big),
	\end{equation}
	where $\gamma$ and $ K_\beta(x_1, x_2) $ are defined as in \Cref{table:parameters} and \eqref{eq:correlation}. Then for any fixed $t \geq 0 $ and for any Schwartz function $ \phi \in \mathcal{S} ( \mathbb{R}^{\dimension} \times [0,1]) $ the covariation process converges in probability
	\begin{equation*}
		\blangle M^N (\phi) \brangle_t \xrightarrow[N \to \infty]{} t \blangle \mathcal{Q}, \phi \otimes \phi \brangle.
	\end{equation*}
\end{restatable}
We postpone the proofs of \Cref{lem:martingalejumps} and \Cref{lem:convergenceofquadraticvariation} to \Cref{subsec:convergenceofmartingalemeasures}, as both require adaptations to hold for the generalized regimes. Following \cite{forien}, the limiting martingale measure of $ (M^N, N \geq 1) $ is then characterised by:
\begin{lemma}[{extends \parencite[Lemma 3.9]{forien}}] \label{lemmaconvergencemartingalemeasures}
Suppose $ M $ is a continuous martingale measure on $ \mathbb{R}_+ \times \mathbb{R}^{\dimension} \times [0,1] $ with covariation
\begin{equation*}
	dt \delta_{t} (ds) \gamma K_\beta(dx_1, dx_2) \big(dk_1 \delta_{k_1} (dk_2) - dk_1 dk_2 \big).
\end{equation*}
Then the sequence $ (M^N, N \geq 1) $ converges in distribution to $ M $ in $ D (\mathbb{R}_+, \mathcal{S}' (\mathbb{R}^{\dimension} \times [0,1])) $.
\end{lemma}
\begin{proof}
Let $ \phi \in \mathcal{S} (\mathbb{R}^{\dimension} \times [0,1]) $ be a Schwartz function and consider the sequence $ ((M_t^N(\phi)_{t \geq 0}, N \geq 1) $. \Cref{lem:martingalejumps} and \Cref{lem:convergenceofquadraticvariation} amount to the two conditions of \Cref{theo:martingaleconvergence}, yielding the convergence of $ (M_t^N (\phi), t \geq 0) $ to $ (M_t(\phi), t \geq 0) $ in distribution in $ D(\mathbb{R}_+, \mathbb{R}) $. \\
Now let $ \phi_1, ..., \phi_k \in \mathcal{S} (\mathbb{R}^{\dimension} \times [0,1]) $ be Schwartz functions. Since both $ \blangle M^N (\phi_i) \brangle_t $ and $ \blangle M^N (\phi_j) \brangle_t $ converge, the covariation $ \blangle M^N (\phi_i), M^N (\phi_j) \brangle_t $ converges by polarization. Applying \Cref{theo:martingaleconvergence} again, we get convergence of finite-dimensional distributions $ (M_t^N(\phi_1),..., M_t^N (\phi_k), t \geq 0) $. Consequently, $ (M^N, N \geq 1) $ satisfies the prerequisites of \parencite[Theorem 6.15]{walsh_introduction_1986} and the convergence of $ (M^N, N \geq 1) $ to $ M $ follows.
\end{proof}
 

\subsubsection{Summary of significant changes required to the proofs of \parencite{forien}} \label{subsec:proof_differences}
In order to unify the arguments, we will generalize the various lemmas of the proof of \parencite{forien} to the one- and two-tail regimes. While the general form of the lemmas has been preserved, the proofs require non-trivial adaptations, which we provide in \Cref{sec:proof_lemmas}. In summary:
\begin{enumerate}
    \item 
    \Cref{lem:boundondominatingmeasures} deals with the bound on the dominating measures. Here, we are able to use the same bounds in new combinations to cover the new regimes.
        \item
   Prior calculations cannot be transfered to \Cref{lem:convergenceofgenerators}, because in the one- and two-tail regime one of the averages of $\mathcal{L}^N$ vanishes asymptotically. 
    \item
    For general sequences $(\delta_N, N \geq 1)$ converging to zero, the proof of \Cref{lem:martingalejumps} does not carry over. It is necessary to impose the new conditions \eqref{assumption_jumps_one_tail} and  \eqref{assumption_jumps_two_tail} on the convergence speed of $\delta_N$ relative to $N$ in the one- respectively two-tail regime.
    \item
    \Cref{lem:convergenceofquadraticvariation} analyses the quadratic variation of the martingale measures. If ${\beta < \dimension}$, convergence becomes much harder to justify, which led us to develop the auxiliary \Cref{lem:boundondominatingmeasures2}. The case $\beta > \dimension$ has been adapted to two different radii and the case $\beta = \dimension$ involves new calculations to account for the novel logarithmic scaling $\eta_N = \delta_N^{\alpha - \dimension} \log^{-1} (1/\delta_N)$.
\end{enumerate}

\subsection{Derivation of the Wright-Mal\'ecot formulae} \label{subsec:WMFproof}
Let us explain how one derives the Wright-Mal\'ecot formula from the central limit theorem, analogously to \cite[Section 4]{forien}.
  We are interested in $ P_t^N (\phi, \psi) $, the probability that two individuals sampled independently according to $ \phi, \psi $ from the $ N $-th rescaled SLFV model $ (\boldsymbol{\rho}_t^N, t \geq 0) $ are of the same type. The probability densities $\phi, \psi$ are assumed to be smooth and compactly supported. Recall from \Cref{subsec:WMF} that
\begin{equation*}
P_t^N (\phi , \psi) = \mathbb{E} \big[ \blangle \boldsymbol{\rho}_t^N \otimes \boldsymbol{\rho}_t^N , (\phi \otimes \psi) \mathds{1}_{\Delta} \brangle \big]
\end{equation*}
where $ \mathds{1}_\Delta : [0,1]^2 \rightarrow \mathbb{R} $ and $\Delta = \lbrace (k_1, k_2) \in [0,1]^2 : k_1 = k_2 \rbrace$. Since $ Z_t^N := \sqrt{N \eta_N } (\boldsymbol{\rho}_t^N - \lambda) $ and Lebesgue measure $ \lambda $ does not charge the diagonal $\Delta$, this is equal to
\begin{equation*}
P_t^N (\phi, \psi) = (N \eta_N)^{-1} \mathbb{E} \big[ \blangle Z_t^N \otimes Z_t^N , (\phi \otimes \psi ) \mathds{1}_\Delta \brangle \big].
\end{equation*}
As a first step, we need to justify that the exchange of the limit $N \to \infty$ with the expectation is possible.
\begin{lemma} 
For $t \geq 0$, two compactly supported and smooth densities $\phi, \psi$ and $Z_t$ the solution to \eqref{SPDE} we have
\begin{equation*}
	\lim_{N \to \infty} N \eta_N P_t^N (\phi , \psi) = \mathbb{E} \big[ \blangle Z_t \otimes Z_t , (\phi \otimes \psi) \mathds{1}_\Delta \brangle \big].
\end{equation*}
\end{lemma}
\begin{proof}
From \Cref{prop:fluctuationsrescaled} and \Cref{eq:stochasticintegralrepresentation} we know
\begin{equation*}
N \eta_N P_t^N (\phi , \psi) = \mathbb{E} \left[ \eta_N \int_0^t \blangle \Gamma^{\nu_N^\alpha} ( \boldsymbol{\rho}_s^N ) , \big( P_{t-s}^N \phi \otimes P_{t-s}^N \psi \big) \mathds{1}_\Delta \brangle ds \right].
\end{equation*}
Combining \Cref{lem:boundondominatingmeasures} and \Cref{lem:semigroupproperties} yields
\begin{equation*}
\begin{aligned}
&\Big\vert \blangle \eta_N \Gamma^{\nu_N^\alpha} ( \boldsymbol{\rho}_s^N ) , \big( P_{t-s}^N \phi \otimes P_{t-s}^N \psi \big) \mathds{1}_\Delta \brangle - \blangle \eta_N \Gamma^{\nu_N^\alpha} ( \boldsymbol{\rho}_s^N ) , \big( P_{t-s} \phi \otimes P_{t-s} \psi \big) \mathds{1}_\Delta \brangle \Big\vert\\
&\hspace{8cm} \leq C \vert t - s \vert^3 \delta_N^\theta e^{-2 \mu (t-s)}.
\end{aligned}
\end{equation*}
As in \Cref{lem:convergenceofquadraticvariation} one can then deduce
\begin{equation*}
\eta_N \int_0^t \blangle \Gamma^{\nu_N^\alpha} ( \boldsymbol{\rho}_s^N ) , \big( P_{t-s} \phi \otimes P_{t-s} \psi \big) \mathds{1}_\Delta \brangle ds \xrightarrow[N \to \infty]{} \int_0^t \blangle \mathcal{Q} , \big( P_s \phi \otimes P_s \psi \big) \mathds{1}_\Delta \brangle ds.
\end{equation*}
Again by \Cref{lem:boundondominatingmeasures} and \Cref{lem:semigroupproperties}, the dominated convergence theorems allows to exchange limit and expectation.
\end{proof}
\begin{proof}[Proof of \Cref{wmf}]
From the previous proof we know
\begin{equation*}
\mathbb{E} \big[ \blangle Z_t \otimes Z_t , (\phi \otimes \psi) \mathds{1}_\Delta \brangle \big] = \int_0^t \blangle \mathcal{Q}, \big( P_s \phi \otimes P_s \psi \big) \mathds{1}_{\Delta} \brangle ds.
\end{equation*}
The integral on the right-hand side stays finite by \Cref{lem:boundondominatingmeasures} and \Cref{lem:semigroupproperties} as $t\to \infty$, which shows
\begin{equation*}
\lim_{t \to \infty} \mathbb{E} \big[ \blangle Z_t \otimes Z_t , (\phi \otimes \psi) \mathds{1}_\Delta \brangle \big] = \int_0^\infty \blangle \mathcal{Q}, \big( P_s \phi \otimes P_s \psi \big) \mathds{1}_{\Delta} \brangle ds.
\end{equation*}
If we let $G_t^\alpha$ be the fundamental solution associated to $\partial_t - \mathcal{D}_\alpha$, we can rewrite the above equation with the representation
\begin{equation*}
P_t \phi (x,k) = e^{- \mu t} \int_{\mathbb{R}^{\dimension}} G_{ut}^\alpha (x-y) \phi (y,k) dy.
\end{equation*}
The expressions of \Cref{wmf} follow as $\mathds{1}_\Delta$ integrated with respect to the types part of $\mathcal{Q}$ found in \eqref{eq:equationQ} yields
\begin{equation*}
\int_{[0,1]^2} \mathds{1}_\Delta (k_1, k_2) (\delta_{k_1} (dk_2) dk_1 - dk_1 dk_2 ) = 1. \qedhere
\end{equation*}
\end{proof}

\section{Proofs of the lemmas} \label{sec:proof_lemmas}
This section collects all proofs requiring significant changes to generalize \cite{forien} to the new regimes.
\subsection{Bound on dominating measures} \label{subsec:bound_dominating_measures}
\begin{proof}[Proof of \Cref{lem:boundondominatingmeasures}]
The first estimates are the same as in \parencite[p.20f]{forien}, then we use new combinations of the bounds for the various cases. Recall
	\begin{equation*}
		\blangle D^N , \mathds{1}_{[s,t]} \phi \otimes \psi \brangle = \int_s^t \blangle \eta_N \vert \Gamma \vert^{\nu_N^\alpha} ( \boldsymbol{\rho}_v^N ) , \phi \otimes \psi \brangle d v.
	\end{equation*}
	Considering the four summands of $ \vert \Gamma \vert^{\nu_N^\alpha} $ separately, we can bound $ \blangle D^N , \mathds{1}_{[s,t]} \phi \otimes \psi \brangle $ by
	\begin{equation} \label{firstbounddominating}
 \begin{aligned}
		&4 \vert t - s \vert \int_{(0,1] \times (0,\infty)^2} \eta_N (N u)^2 \int_{(\mathbb{R}^{\dimension})^2} V_{\delta_N r_2} (x_1, x_2) \\
	&\hspace{1cm}\times \sup_{k \in [0,1]} \vert \phi (x_1,k) \vert \sup_{k \in [0,1]} \vert \psi (x_2, k) \vert dx_1 dx_2 \delta_N^{- (\dimension + \alpha)} \nu_N (du, dr_1, dr_2).
  \end{aligned}
	\end{equation}
We defined $ L^p $-norms in \eqref{equationlpnorm} to involve the supremum over $ k \in [0,1] $ and thus applying $ V_{\delta_N r_2} (x_1, x_2) \leq V_{\delta_N r_2} $ the inner integral is bounded by $ V_{\delta_N r_2} \Vert \phi \Vert_1 \Vert \psi \Vert_1 $. Alternatively, the Cauchy-Schwarz inequality results in
	\begin{equation} \label{boundonintegrals}
		\int_{(\mathbb{R}^{\dimension})^2} V_{\delta_N r_2} (x_1, x_2) \sup_{k \in [0,1]} \vert \phi (x_1, k) \vert \sup_{k \in [0,1]} \vert \psi (x_2, k) \vert dx_1 dx_2 
		\leq V_{\delta_N r_2}^2 \Vert \phi \Vert_2 \Vert \psi \Vert_2,
	\end{equation}
	since
	\begin{equation*}
		\int_{\mathbb{R}^{\dimension}} V_{\delta_N r_2 } (x_1, x_2) dx_2 = V_{\delta_N r_2}^2.
	\end{equation*}
We only demonstrate how the two estimates can be used to bound \eqref{firstbounddominating} for the one-tail regime. The two-tail regime follows in a similar fashion, as the double integrals can be written as a product of two integrals.

\paragraph{One-tail regime} If $\beta = \tail + \impact < d$, we split the integral in \eqref{firstbounddominating} and apply both bounds to obtain
	\begin{equation*}
		\begin{aligned}
			4 \vert t - s \vert \Bigg\{ &\int_1^{1/\delta_N} \delta_N^{\alpha - \tail - \impact} \frac{u_0^2}{r_2^{2 \impact}} V_{\delta_N r_2}^2 \delta_N^{- (\dimension + \alpha)} \frac{dr_2}{r_2^{1 + \dimension + \tail - \impact}} \Vert \phi \Vert_2 \Vert \psi \Vert_2\\
			& + \int_{1 /\delta_N}^\infty \delta_N^{\alpha - \tail - \impact} \frac{u_0^2}{r_2^{2 \impact}} V_{\delta_N r_2} \delta_N^{- (\dimension +\alpha)} \frac{dr_2}{r_2^{1 + \dimension + \tail - \impact}} \Vert \phi \Vert_1 \Vert \psi \Vert_1 \Bigg\}.
		\end{aligned}
	\end{equation*}
	The first integral can be bounded by a finite quantity independent of $N$
	\begin{align*}
	\int_1^{1/\delta_N} \delta_N^{\alpha - \tail - \impact} \frac{u_0^2}{r_2^{2 \impact}} V_{\delta_N r_2}^2 \delta_N^{- (\dimension + \alpha)} \frac{dr_2}{r_2^{1 + \dimension + \tail - \impact}}
	&= u_0^2 \int_1^{1/\delta_N}  V_{\delta_N r_2}^2 \delta_N \frac{dr_2}{(\delta_N r_2)^{1 + \dimension +\tail + \impact}} \\
    &= u_0^2 V_1^2 \int_{\delta_N}^{1}  \frac{dr_2}{r_2^{1 + \beta - \dimension}} \leq \frac{u_0^2 V_1^2}{\dimension - \beta},
	\end{align*}
since $\beta < \dimension$.
The second integral is equal to
\begin{align*}
	\int_{1 /\delta_N}^\infty \delta_N^{\alpha - \tail - \impact} \frac{u_0^2}{r_2^{2 \impact}} V_{\delta_N r_2} \delta_N^{- (\dimension + \alpha)} \frac{dr_2}{r_2^{1 + \dimension + \tail - \impact}} &= u_0^2 \int_{1/\delta_N}^\infty V_{\delta_N r_2} \delta_N \frac{dr_2}{(\delta_N r_2)^{1 + \dimension + \tail + \impact}} \\
 &= u_0^2 V_1 \int_{1}^\infty \frac{dr_2}{ r_2^{1 + \tail + \impact}} = \frac{u_0^2 V_1}{\beta}.
	\end{align*}
Similarly, if $\beta = \dimension$, the first integral is bounded by
\begin{align*}
\int_1^{1/\delta_N} \frac{\delta_N^{\alpha - \dimension}}{\log (1 /\delta_N )} \frac{1}{r_2^{2 \impact}} V_{\delta_N r_2}^2 \delta_N^{- (\dimension + \alpha )} \frac{dr_2}{r_2^{1 + \dimension + \tail - \impact}} \Vert \phi \Vert_2 \Vert \psi \Vert_2 &= \log (1 /\delta_N)^{-1} V_1^2 \int_1^{1/\delta_N} \frac{dr_2}{r_2} \Vert \phi \Vert_2 \Vert \psi \Vert_2 \\
&= V_1^2 \Vert \phi \Vert_2 \Vert \psi \Vert_2,
\end{align*}
and the second integral can be bounded by
\begin{align*}
\int_{1/\delta_N}^\infty \frac{\delta_N^{\alpha - \dimension}}{\log (1 /\delta_N )} \frac{1}{r_2^{2 \impact}} V_{\delta_N r_2} \delta_N^{- (\dimension + \alpha )} \frac{dr_2}{r_2^{1 + \dimension + \tail - \impact}} \Vert \phi \Vert_1 \Vert \psi \Vert_1 &= \frac{\delta_N^{- \dimension}}{\log (1 /\delta_N) } V_1 \int_{1/\delta_N}^\infty  \frac{dr_2}{r_2^{1 + \dimension}} \Vert \phi \Vert_1 \Vert \psi \Vert_1 \\
&= \frac{V_1 \Vert \phi \Vert_1 \Vert \psi \Vert_1}{\dimension \log (1/\delta_N)}.
\end{align*}
If $\beta > \dimension$, we can apply the Cauchy-Schwarz bound directly to \eqref{firstbounddominating}
\begin{align*}
\int_1^\infty \delta_N^{\alpha - \dimension} \frac{1}{r_2^{2 \impact}} V_{\delta_N r_2}^2 \delta_N^{- (\dimension + \alpha )} \frac{dr_2}{r_2^{1 + \dimension + \tail - \impact}} \Vert \phi \Vert_2 \Vert \psi \Vert_2 &= V_1^2 \int_1^\infty \frac{dr_2}{r_2^{1 + \tail + \impact - \dimension}} \Vert \phi \Vert_2 \Vert \psi \Vert_2 \\
&= \frac{V_1^2 \| \phi \|_2 \| \psi \|_2}{\beta - \dimension}.
\end{align*}
This concludes the proof of the lemma.
\end{proof}
In the proof of \Cref{lem:convergenceofquadraticvariation}, we will require a version of \Cref{lem:boundondominatingmeasures} for functions $\Psi \in \mathcal{S} ((\mathbb{R}^{\dimension})^2 \times [0,1]^2)$.
\begin{corollary} \label{lem:boundondominatingmeasures2}
	Let $ D^N $ be the dominating measure from \eqref{dominatingmeasuresdefinition}. Then in the one-tail and two-tail regimes we can find constants $ C > 0, q >1 $ such that for all $ N \geq 1, 0 \leq s  \leq t $ and $ \Psi \in \mathcal{S} ( (\mathbb{R}^{\dimension})^2 \times [0,1]^2) $
	\begin{equation} \label{boundondominatingmeasures2}
		\blangle D^N, \mathds{1}_{[s,t]} \Psi \brangle \leq C \vert t - s \vert \Bigg( \Vert \Psi \Vert_1 + \int_{\mathbb{R}^{\dimension}} \Big( \int_{\mathbb{R}^{\dimension}} \big\vert \Psi (x_1, x_2) \big\vert^q dx_2 \Big)^{1/q} dx_1 \Bigg),
	\end{equation}
 where we abbreviate $\sup_{k_1, k_2 \in [0,1]^2} \Psi (x_1, x_2, k_1, k_2 )$ by $\Psi (x_1, x_2)$.
\end{corollary}
\begin{proof}
To prove \eqref{boundondominatingmeasures2}, we need to bound the inner integral of \eqref{firstbounddominating} for functions $ \Psi \in \mathcal{S} ( (\mathbb{R}^{\dimension})^2 \times [0,1]^2) $. Using $V_{\delta_N r_2} (x_1, x_2 ) \leq V_{\delta_N r_2}$, we obtain the bound $V_{\delta_N r_2} \Vert \Psi \Vert_1$. Alternatively, we can use Hölder's inequality with $\frac{1}{p} + \frac{1}{q} = 1$ and set $\xi = 1/p \in [0,1)$
\begin{multline*}
\int_{(\mathbb{R}^{\dimension})^2} V_{\delta_N r_2} (x_1, x_2 )  \Psi (x_1, x_2 ) dx_1 dx_2\\
\begin{aligned}
    &\leq \int_{\R^{\dimension}} \left( \int_{\mathbb{R}^{\dimension}} V_{\delta_N r_2} (x_1, x_2)^{p \xi} dx_2 \right)^{1/p} \Bigg( \int_{\R^{\dimension}} \underbrace{V_{\delta_N r_2}^{q (1- \xi)} (x_1, x_2)}_{\leq V_{\delta_N r_2}} \big\vert \Psi (x_1, x_2 ) \big\vert^q dx_2 \Bigg)^{1/q} dx_1\\
&\leq V_{\delta_N r_2}^{\frac{2}{p} + \frac{1}{q}} \int_{\mathbb{R}^{\dimension}} \left( \int_{\mathbb{R}^{\dimension}} \Big\vert \Psi (x_1, x_2 ) \Big\vert^q dx_2 \right)^{1/q} dx_1.
\end{aligned}
\end{multline*}
We can choose $\frac{2}{p} + \frac{1}{q}$ arbitrarily close to $2$, which allows us to apply the same arguments as in \Cref{lem:boundondominatingmeasures}.
\end{proof}
\subsection{Convergence of generators} \label{subsec:convergence_of_generators}
In this section, we show the convergence of the generator $ \mathcal{L}^N $ from \eqref{equationrescaledgenerator} given by
\begin{align*}
	\mathcal{L}^N\phi (x)
	&= N \delta_N^{- \alpha} \int_{(0,N] \times (0, \infty)^2} u V_{r_2} \left(\overline{\overline{\phi}} (x , \delta_N r_1 , \delta_N r_2) - \phi (x) \right) \nu_N (du, dr_1, dr_2)
\end{align*}
to the fractional Laplacians $\mathcal{D}_\alpha $ defined in \eqref{alphastableprocessgenerator}.
Since the operators $\mathcal{L}^N$ only act on the spatial variable, we omit the dependence on the second variable $k$ in this whole section.
The main idea of the proof is that under our scaling one of the averaging radii is much smaller than the other, and in the limit only the largest averaging radius matters. 
The only case where the two radii are precisely the same, that is $ \power = 1 $ in the one-tail regime, was already covered by \parencite[Appendix A]{forien_central_2017}.
The proof will require the following lemmas.

\begin{lemma} \label{lpboundaverages}
Assume that $ \phi : \mathbb{R}^{\dimension} \rightarrow \mathbb{R} $ has continuous derivatives up to the second order and $ \max_{\vert \kappa \vert \leq 2} \Vert \partial_\kappa \phi \Vert_{q} $ is finite for $ 1 \leq q \leq \infty $. Then for all $ 1 \leq p \leq \infty $
\begin{equation} \label{equationaveragebound1}
\Big\Vert \overline{\phi} (r) - \phi \Big\Vert_p \leq \min \Big\{2 \Vert \phi \Vert_p , \frac{\dimension}{2} r^2 \max_{\vert \kappa \vert = 2} \Vert \partial_\kappa \phi \Vert_p \Big\}.
\end{equation}
\end{lemma}

\begin{proof}
The first estimate of \eqref{equationaveragebound1} stems from $ \Vert \overline{\phi} (r) \Vert_p \leq \Vert \phi \Vert_p $. The second estimate of \eqref{equationaveragebound1} is derived by Taylor's theorem in \parencite[Proposition A.1]{forien_central_2017}.
\end{proof}

\begin{lemma}
Assume that $\phi : \mathbb{R}^{\dimension} \rightarrow \mathbb{R}$ has continuous derivatives up to the fourth order and $ \max_{p \in [1, \infty]} \max_{\vert \kappa \vert = 4} \Vert \partial_\kappa \phi \Vert_p $ is finite. Then for all $1 \leq p \leq \infty$ and all $r,s \in [0,\infty)$
\begin{equation} \label{bound_double_average}
\Bigg\Vert \overline{\overline{\phi}} (x, r,s) - \phi (x) - \frac{r^2 + s^2}{2 (d+2)} \Delta \phi (x) \Bigg\Vert_p \leq \frac{d^3}{3!} \max_{\vert \kappa \vert = 4 } \Vert \partial_\kappa \phi \Vert_p (r^4 + s^4).
\end{equation}
\end{lemma}

\begin{proof}
This is a generalization of \parencite[Proposition A.1]{forien_central_2017} to two radii. 
\end{proof}
\begin{proof}[Proof of \Cref{lem:convergenceofgenerators}.] In the following proof, we focus on the more difficult second bound \eqref{equationboundondifferenceofgenerators}. The first bound \eqref{equationboundongenerator} can be shown using the same estimates.
\paragraph{One-tail regime} In the one-tail regime, substituting $\nu_N$ by its expression in \eqref{one_tail_regime}, we can write
\begin{equation} \label{expr_L_N}
\begin{aligned}
    \mathcal{L}^N\phi(x) &= \delta_N^{-\alpha} \int_1^\infty \frac{u_0}{r^c} V_r \left( \overline{\overline{\phi}}(x, \delta_N r^b, \delta_N r) - \phi(x) \right) \frac{dr}{r^{1+d+a-c}} \\
    &= u_0 V_1 \delta_N^{-\alpha} \int_1^\infty \left( \overline{\overline{\phi}}(x, \delta_N r^b, \delta_N r) - \phi(x) \right) \frac{dr}{r^{1+a}}. 
\end{aligned}
\end{equation}

\subparagraph*{The Brownian case}

Suppose that $a > 2$ and $\frac{a}{b} > 2$, then $\alpha = 2$.
Let $\sigma > 0$ be defined as
\begin{align*}
    \sigma^2 &= \frac{u_0 V_1}{d+2} \int_1^\infty (r^{2\power} + r^2) \frac{dr}{r^{1+\tail}} \\
    &= \frac{u_0 V_1}{d+2} \left( \frac{1}{\tail-2\power} + \frac{1}{\tail-2} \right).
\end{align*}
Let also $(A_N, N \geq 1)$ denote a sequence taking values in $(1, \infty]$ which will be chosen later.
	We can then bound, using the definition of $\sigma^2$,
	\begin{align*}
	    \Bigg\Vert \mathcal{L}^N \phi - \frac{\sigma^2}{2}& \Delta \phi \Bigg\Vert_p \leq u_0 V_1 \int_{1}^{A_N} \left\Vert  \overline{\overline{\phi}} (\cdot, \delta_N r_2^{\power}, \delta_N r) - \phi - \frac{ \delta_N^2 r^{2 \power} + \delta_N^2 r^2}{2 ( \dimension + 2)} \Delta \phi \right\Vert_p \delta_N^{-2} \frac{dr}{r^{1 + \tail}}\\
	    &+ u_0 V_1 \int_{A_N}^\infty \Big\Vert \overline{\overline{\phi}}(\cdot, \delta_N r^b, \delta_N r) - \phi \Big\Vert_p \delta_N^{-2} \frac{dr}{r^{1+\tail}} + \frac{u_0 V_1}{2(d+2)} \int_{A_N}^\infty (r^{2b} + r^2) \frac{dr}{r^{1+\tail}} \Vert \Delta \phi \Vert_p.
	\end{align*}
	Using \eqref{bound_double_average} in the first integral and \eqref{equationaveragebound1} in the second we obtain
	\begin{equation} \label{bound_Laplacian}
	    \left\Vert \mathcal{L}^N \phi - \frac{\sigma^2}{2} \Delta \phi \right\Vert_p \leq C \left( \delta_N^2 \int_1^{A_N} (r^{4b} + r^{4}) \frac{dr}{r^{1+\tail}} + 2 \int_{A_N}^\infty (r^{2b} + r^2) \frac{dr}{r^{1+\tail}} \right) \max_{0 \leq \vert \kappa \vert \leq 4 } \Vert \partial_\kappa \phi \Vert_p,
	\end{equation}
	for some constant $C > 0$.
	First note that, since $a > 2b$ and $a > 2$,
	\begin{align*}
	    \int_{A_N}^\infty (r^{2b} + r^2) \frac{dr}{r^{1+\tail}} = \frac{1}{(a-2b) A_N^{a-2b}} + \frac{1}{(a-2) A_N^{a-2}}.
	\end{align*}
	On the other hand
	\begin{align*}
	    \int_1^{A_N} r^{4} \frac{dr}{r^{1+\tail}} = \begin{cases}
	    \frac{A_N^{4-a}-1}{4-a} & a \neq 4, \\
	    \log(A_N) & a = 4,
	    \end{cases}, && \int_1^{A_N} r^{4b} \frac{dr}{r^{1+\tail}} = \begin{cases}
	    \frac{A_N^{4b-a}-1}{4b-a} & a \neq 4b, \\
	    \log(A_N) & a = 4b.
	    \end{cases}
	\end{align*}
	We can thus choose $A_N = \infty$ if $a > 4 \vee 4b$, and $A_N = \delta_N^{-\xi}$ with $0 < \xi < \frac{2}{4-a} $ if $a \leq 4$ and $0 < \xi < \frac{2}{4b-a}$ if $a \leq 4b$.
	In each case, the right hand side of \eqref{bound_Laplacian} is bounded by $C \delta_N^\theta \max_{0 \leq \vert \kappa \vert \leq 4 } \Vert \partial_\kappa \phi \Vert_p $ for some constant $C > 0$.

\subparagraph*{The stable case}

When $b<1$ and $\tail < 2$, $\alpha = \tail$ and, changing variables in the integral over $r$ in \eqref{expr_L_N}, we obtain
\begin{align} \label{expr_L_onetail}
    \mathcal{L}^N \phi(x) = u_0 V_{1} \int_{\delta_N}^\infty \left( \overline{\overline{\phi}}(x, \delta_N^{1-b} r^b, r) - \phi(x) \right) \frac{dr}{r^{1+a}}.
\end{align}
On the other hand, if $b > 1$ and $\frac{a}{b} < 2$, $\alpha = \frac{a}{b}$ and we instead write
\begin{align*} 
    \mathcal{L}^N\phi(x) = \frac{u_0 V_1}{b} \int_{\delta_N}^\infty  \left( \overline{\overline{\phi}}(x, r, \delta_N^{1-1/b} r^{1/b}) - \phi(x) \right) \frac{dr}{r^{1+\frac{a}{b}}}.
\end{align*}
Note that the two expressions above are of the same form, with $1/b$ in the second one playing the role of $b$ in the first one, so that the case $b > 1$ can be treated in exactly the same way as the case $b > 1$. The case $b = 1$ has already been covered in \cite{forien_central_2017}.

We start by detailing the case $b < 1$.
We see that, in \eqref{expr_L_onetail}, $\delta_N^{1-b} r^b \to 0$ as $N \to \infty$. Let us set
\begin{equation*} 
    \mathcal{D}^N \phi(x) := u_0 V_1 \int_{\delta_N}^\infty \left( \overline{\phi}(x, r) - \phi(x) \right) \frac{dr}{r^{1+\tail}}.
\end{equation*}
As a result,
\begin{align*}
    \| \mathcal{L}^N \phi - \mathcal{D}^N \phi \|_p \leq u_0 V_1 \int_{\delta_N}^{\infty} \| \overline{\overline{\phi}}(\cdot, \delta_N^{1-b} r^b, r) - \overline{\phi}(\cdot,r) \|_p \frac{dr}{r^{1+\tail}}.
\end{align*}
If $b = 0$, we can use \eqref{equationaveragebound1} to obtain
\begin{align*}
    \| \mathcal{L}^N \phi - \mathcal{D}^N \phi \|_p &\leq C \delta_N^2 \int_{\delta_N}^\infty \frac{dr}{r^{1+\tail}} \max_{\vert \kappa \vert = 2} \Vert \partial_\kappa \phi \Vert_p \\
    &\leq C \delta_N^{2-\tail} \max_{\vert \kappa \vert = 2} \Vert \partial_\kappa \phi \Vert_p.
\end{align*}
If $b > 0$, using \eqref{equationaveragebound1} for $r \leq \delta_N^{-\frac{1-b}{b}}$, we obtain, for $a \neq 2b$,
\begin{align*}
    \| \mathcal{L}^N \phi - \mathcal{D}^N \phi \|_p &\leq C \left( \delta_N^{2(1-b)} \max_{\vert \kappa \vert = 2} \Vert \partial_\kappa \phi \Vert_p \int_{\delta_N}^{\delta_N^{-\frac{1-b}{b}}} r^{2b-a-1} dr  + 2 \| \phi \|_p \int_{\delta_N^{-\frac{1-b}{b}}}^{\infty} \frac{dr}{r^{1+\tail}} \right) \\
    &\leq C \left( \delta_N^{2(1-b)} \frac{\delta_N^{-(2b-a)\frac{1-b}{b}} - \delta_N^{2b-a}}{2b-a} + \frac{2}{a} \delta_N^{\frac{a(1-b)}{b}} \right) \max_{0 \leq \vert \kappa \vert \leq 2} \Vert \partial_\kappa \phi \Vert_p \\
    &\leq C \left( \delta_N^{\frac{a(1-b)}{b}} + \delta_N^{2-a} \right) \max_{0 \leq \vert \kappa \vert \leq 2} \Vert \partial_\kappa \phi \Vert_p.
\end{align*}
If $a = 2b$, we instead obtain
\begin{align*}
    \| \mathcal{L}^N \phi - \mathcal{D}^N \phi \|_p &\leq C \left( \delta_N^{2(1-b)} |\log(\delta_N)| + \delta_N^{\frac{a(1-b)}{b}} \right) \max_{0 \leq \vert \kappa \vert \leq 2} \Vert \partial_\kappa \phi \Vert_p.
\end{align*}
Since $\tail < 2$ and $b < 1$, in every case the right hand side is of the form $C \delta_N^\theta \max_{0 \leq \vert \kappa \vert \leq 2 } \Vert \partial_\kappa \phi \Vert_p $ for some constant $C > 0$.

It thus remains to bound $\| \mathcal{D}^N \phi - \mathcal{D}_\alpha \phi \|_p$.
Using \eqref{equationaveragebound1},
\begin{align*}
    \| \mathcal{D}^N \phi - \mathcal{D}_\alpha \phi \|_p &\leq u_0 V_1 \int_{0}^{\delta_N} \| \overline{\phi}(\cdot,r) - \phi \|_p \frac{dr}{r^{1+\tail}} \\
    &\leq C \int_0^{\delta_N} r^{2-\tail-1} dr \max_{\vert \kappa \vert = 2} \Vert \partial_\kappa \phi \Vert_p \\
    &\leq C \delta_N^{2-\tail} \max_{\vert \kappa \vert = 2} \Vert \partial_\kappa \phi \Vert_p,
\end{align*}
for some constant $C > 0$.
This concludes the proof of \eqref{equationboundondifferenceofgenerators}.

\paragraph*{Two-tails regime}

In the two-tails regime, similar calculations as in \eqref{expr_L_N} lead to
\begin{align} \label{expr_L_twotails}
    \mathcal{L}^N \phi(x) = u_0 V_{1} \delta_N^{-\alpha} \int_1^\infty \int_1^\infty  \left( \overline{\overline{\phi}}(x, \delta_N r_1, \delta_N r_2) - \phi(x) \right) \frac{dr_1 dr_2}{r_1^{1+a_1} r_2^{1+a_2}}.
\end{align}

\subparagraph*{The Brownian case}

Suppose that $\tail_1 > 2$ and $\tail_2 > 2$.
Then $\alpha = 2$ and we set
\begin{align*}
    \sigma^2 &= \frac{u_0 V_1}{d+2} \int_1^\infty \int_1^\infty (r_1^2 + r_2^2) \frac{dr_1 dr_2}{r_1^{1+\tail_1} r_2^{1+\tail_2}} \\
    &= \frac{u_0 V_1}{d+2} \left( \frac{1}{(\tail_1-2) \tail_2} + \frac{1}{\tail_1(\tail_2-2)} \right).
\end{align*}
Then for $A_N \in (1,\infty]$,
\begin{equation*}
\begin{aligned}
    &\left\| \mathcal{L}^N \phi - \frac{\sigma^2}{2} \Delta \phi \right\|_p\\
    &\leq u_0 V_1 \int_1^{A_N} \int_1^{A_N} \left\| \overline{\overline{\phi}}(\cdot, \delta_N r_1, \delta_N r_2) - \phi - \frac{\delta_N^2 r_1^2 + \delta_N^2 r_2^2}{2(d+2)} \Delta \phi \right\|_p \delta_N^{-2} \frac{dr_1 dr_2}{r_1^{1+\tail_1} r_2^{1+\tail_2}} \\
    &\quad + u_0 V_1 \int_1^{A_N} \int_{A_N}^\infty \left( \left\| \overline{\overline{\phi}}(\cdot, \delta_N r_1, \delta_N r_2) - \phi \right\|_p \delta_N^{-2} + \frac{r_1^2 + r_2^2}{2(d+2)} \| \Delta \phi \|_p \right) \frac{dr_1 dr_2}{r_1^{1+\tail_1} r_2^{1+\tail_2}} \\
    &\quad + u_0 V_1 \int_{A_N}^{\infty} \int_1^\infty \left( \left\| \overline{\overline{\phi}}(\cdot, \delta_N r_1, \delta_N r_2) - \phi \right\|_p \delta_N^{-2} + \frac{r_1^2 + r_2^2}{2(d+2)} \| \Delta \phi \|_p \right) \frac{dr_1 dr_2}{r_1^{1+\tail_1} r_2^{1+\tail_2}}.
\end{aligned}
\end{equation*}
Using \eqref{bound_double_average} and \eqref{equationaveragebound1}, we obtain
\begin{equation*}
\begin{aligned}
    &\left\| \mathcal{L}^N \phi - \frac{\sigma^2}{2} \Delta \phi \right\|_p \\
    &\leq C \int_1^{A_N} \int_1^{A_N}  \delta_N^2 ( r_1^4 + r_2^4 ) \frac{dr_1 dr_2}{r_1^{1+\tail_1} r_2^{1+\tail_2}} \max_{\vert \kappa \vert = 4 } \Vert \partial_\kappa \phi \Vert_p \\
    &\quad + C \left( \int_1^{A_N} \int_{A_N}^\infty (r_1^2 + r_2^2) \frac{dr_1 dr_2}{r_1^{1+\tail_1} r_2^{1+\tail_2}} + \int_{A_N}^{\infty} \int_1^\infty (r_1^2 + r_2^2) \frac{dr_1 dr_2}{r_1^{1+\tail_1} r_2^{1+\tail_2}} \right) \max_{\vert \kappa \vert = 2 } \Vert \partial_\kappa \phi \Vert_p.
\end{aligned}
\end{equation*}
The bracketed term is equal to
\begin{align*}
    \frac{1-A_N^{-\tail_2}}{\tail_2 (\tail_1-2) A_N^{\tail_1-2}} + \frac{1-A_N^{-(\tail_2-2)}}{\tail_1(\tail_2-2) A_N^{\tail_1}} + \frac{1}{(\tail_1-2) \tail_2 A_N^{\tail_2}} + \frac{1}{\tail_1 (\tail_2-2) A_N^{\tail_2-2}} \leq C A_N^{-\theta},
\end{align*}
for some constants $C > 0$ and $\theta > 0$.
On the other hand, if $\tail_1 \neq 4$,
\begin{align*}
    \int_1^{A_N} \int_1^{A_N}  \delta_N^2 r_1^4 \frac{dr_1 dr_2}{r_1^{1+\tail_1} r_2^{1+\tail_2}} &= \delta_N^2 \frac{1-A_N^{-\tail_2}}{\tail_2} \times \frac{A_N^{4-\tail_1}-1}{4-\tail_1},
\end{align*}
while if $\tail_1 = 4$,
\begin{align*}
    \int_1^{A_N} \int_1^{A_N}  \delta_N^2 r_1^4 \frac{dr_1 dr_2}{r_1^{1+\tail_1} r_2^{1+\tail_2}} &= \delta_N^2 \frac{1-A_N^{-\tail_2}}{\tail_2} \log(A_N).
\end{align*}
The same integral with $r_2^4$ instead of $r_1^4$ is computed in the same way.
As a result, if both $\tail_1 > 4$ and $\tail_2 > 4$, we can take $A_N = \infty$, otherwise we chose $A_N = \delta_N^{-\xi} $ where $\xi$ is such that
\begin{align*}
    0 < \xi < \frac{2}{(4-\tail_1)^+} \wedge \frac{2}{(4-\tail_2)^+}.
\end{align*}
This yields
\begin{align*}
    \left\| \mathcal{L}^N \phi - \frac{\sigma^2}{2} \Delta \phi \right\|_p \leq C \delta_N^\theta \max_{0 \leq \vert \kappa \vert \leq 4} \Vert \partial_\kappa \phi \Vert_p,
\end{align*}
for some $C > 0$ and $\theta > 0$.

\subparagraph*{The stable case}

Let us now treat the case where $\tail_1 \wedge \tail_2 < 2$.
From \eqref{expr_L_twotails}, it is clear that $\tail_1$ and $\tail_2$ play symmetric roles, so let us assume that $\tail_1 < \tail_2$.
In this case, $\alpha = \tail_1$ and we change variables in the integral over $r_1$ to obtain
\begin{align*}
    \mathcal{L}^N \phi(x) = u_0 V_1 \int_{\delta_N}^\infty \int_1^\infty \left( \overline{\overline{\phi}}(x, r_1, \delta_N r_2) - \phi(x) \right) \frac{dr_1 dr_2}{r_1^{1+\tail_1} r_2^{1+\tail_2}}.
\end{align*}
Then, setting
\begin{align*}
    \mathcal{D}^N \phi(x) = \frac{u_0 V_1}{\tail_2} \int_{\delta_N}^\infty \left( \overline{\phi}(x,r) - \phi(x) \right) \frac{dr}{r^{1+\tail_1}},
\end{align*}
and using \eqref{equationaveragebound1}, we write
\begin{align*}
    \| \mathcal{L}^N \phi - \mathcal{D}^N \phi \|_p &\leq u_0 V_1 \int_{\delta_N}^\infty \int_1^\infty \| \overline{\overline{\phi}}(\cdot, r_1, \delta_N r_2) - \overline{\phi}(\cdot, r_1) \|_p \frac{dr_1 dr_2}{r_1^{1+\tail_1} r_2^{1+\tail_2}} \\
    &\leq C \int_{\delta_N}^\infty \int_1^\infty (1 \wedge (\delta_N r_2)^2 ) \frac{dr_1 dr_2}{r_1^{1+\tail_1} r_2^{1+\tail_2}} \max_{0 \leq \vert \kappa \vert \leq 2} \Vert \partial_\kappa \phi \Vert_p \\
    &\leq C \delta_N^{-\tail_1} \left( \delta_N^2 \int_1^{1/\delta_N} r^{1-\tail_2} dr + \int_{1/\delta_N}^\infty \frac{dr}{r^{1+\tail_2}} \right) \max_{0 \leq \vert \kappa \vert \leq 2} \Vert \partial_\kappa \phi \Vert_p.
\end{align*}
If $\tail_2 \neq 2$, this yields
\begin{align*}
    \| \mathcal{L}^N \phi - \mathcal{D}^N \phi \|_p &\leq C \delta_N^{-\tail_1} \left( \delta_N^2 \frac{\delta_N^{-(2-\tail_2)} - 1}{2-\tail_2} + \delta_N^{\tail_2} \right) \max_{0 \leq \vert \kappa \vert \leq 2} \Vert \partial_\kappa \phi \Vert_p \\
    &\leq C \delta_N^{-\tail_1} \left( \frac{\delta_N^{\tail_2} - \delta_N^2}{2-\tail_2} + \delta_N^{\tail_2} \right) \max_{0 \leq \vert \kappa \vert \leq 2} \Vert \partial_\kappa \phi \Vert_p.
\end{align*}
The right hand side above is of the form $C \delta_N^\theta \max_{0 \leq \vert \kappa \vert \leq 2} \Vert \partial_\kappa \phi \Vert_p$ for $C > 0$ and $\theta > 0$ since $\tail_1 < \tail_2$ and $\tail_1 < 2$.
If $\tail_2 = 2$, the above yields instead
\begin{align*}
    \| \mathcal{L}^N \phi - \mathcal{D}^N \phi \|_p &\leq C \delta_N^{-\tail_1} \left( \delta_N^2 |\log(\delta_N)| + \delta_N^{\tail_2} \right) \max_{0 \leq \vert \kappa \vert \leq 2} \Vert \partial_\kappa \phi \Vert_p,
\end{align*}
which can be bounded by $C \delta_N^\theta \max_{0 \leq \vert \kappa \vert \leq 2} \Vert \partial_\kappa \phi \Vert_p$ for some $C > 0$ and $\theta > 0$ for the same reasons.
\end{proof}
\subsection{Convergence of the martingales measures} \label{subsec:convergenceofmartingalemeasures}
\subsubsection{Martingale Jumps}
We now prove that the martingale measures become continuous in the limit.
\begin{proof}[Proof of \Cref{lem:martingalejumps}]
	We start as in \parencite[p.26f]{forien} and then adapt bounds to implement our regimes. Note
	\begin{align*}
		M_t^N (\phi) - M_{t^-}^N (\phi) &= \blangle Z_t^N - Z_{t-}^N, \phi \brangle - \int_{t-}^t \blangle Z_s^N, \mathcal{L}^N (\phi) - \mu \phi \brangle ds \\
		&= \sqrt{N \eta_N} \blangle \boldsymbol{\rho}_t^N - \boldsymbol{\rho}_{t-}^N, \phi \brangle\\
		&= \sqrt{N \eta_N} \Big( \blangle \rho_{Nt/\delta_N^\alpha}^N , \phi_N \brangle - \blangle \rho_{(Nt /\delta_N^\alpha)^-}^N , \phi_N \brangle \Big).
	\end{align*}
	If $ (t, x_0, u, r_1, r_2) \in \Pi $,
	\begin{align*}
		\left\vert \blangle \rho_t^N , \phi_N \brangle - \blangle \rho_{t^-}^N , \phi_N \brangle \right\vert &\leq u \sup_{k_0 \in [0,1]} \left\vert \int_{\mathbb{R}^{\dimension} \times [0,1]} \phi_N (x,k) \, \mathds{1}_{\vert x -x_0 \vert < r_2} \, ( \delta_{k   _0} (dk) - \rho_{t^-}^N (x, dk)) dx \right\vert \\
		&\leq 2 u \int_{\mathbb{R}^{\dimension}} \delta_N^{\dimension} \sup_{k \in [0,1]} \vert \phi (\delta_N x, k) \vert \mathds{1}_{\vert x- x_0 \vert < r_2} dx.
	\end{align*}
	Using $ \mathds{1}_{\vert x - x_0 \vert < r_2} \leq 1 $ on the one hand and $ |\phi(\delta_N x, k)| \leq \| \phi \|_\infty $ on the other we get
	\begin{equation*}
		\int_{\mathbb{R}^{\dimension}} \delta_N^{\dimension} \sup_{k \in [0,1]} \vert \phi (\delta_N x, k) \vert \mathds{1}_{\vert x - x_0 \vert \leq r_2} dx \leq \min \Big\{ \Vert \phi \Vert_1 , V_{r_2} \delta_N^d \| \phi \|_\infty \Big\}.
	\end{equation*}
	As a result, in the one-tail regime, there exists a constant $C > 0$, depending only on $\phi$ and $\dimension$, such that 
	\begin{align*}
		\sup_{t \geq 0 } \Big\vert M_t^N (\phi) - M_{t^-}^N (\phi) \Big\vert & \leq C \sqrt{\frac{\eta_N}{N}} \sup_{r_2 > 1} \left( \frac{1 \wedge (\delta_N r_2)^{\dimension}}{r_2^{\impact}} \right) \\
		& \leq C \sqrt{\frac{\eta_N}{N}} (\delta_N)^{{\dimension} \wedge {\impact}}.
	\end{align*}
	Here, we used the fact that, if $(t, x_0, u, r_1, r_2) \in \Pi$, then $u = u_0/(r_2^{\impact} N) $. The second inequality results from noting that
 \begin{equation*}
    \sup_{r > 1} \Big( \frac{1 \land (r \delta_N)^{\dimension} }{r^{\impact}} \Big) = \begin{dcases}
    \delta_N^d r^{\dimension-\impact} & \text{ if } 1 < r \leq \frac{1}{\delta_N}, \\
    r^{-\impact} & \text{ if } r > \frac{1}{\delta_N}.
    \end{dcases}
 \end{equation*}
 If $\dimension < \impact$, this function is decreasing and is bounded by $\delta_N^{\dimension}$. In case $\dimension > \impact$, the function assumes the maximum for $r = \frac{1}{\delta_N}$ in $\delta_N^{\impact}$.
    In the two-tails regime, there exists a constant $C > 0$ such that
	\begin{align*}
	    \sup_{t \geq 0 } \Big\vert M_t^N (\phi) - M_{t^-}^N (\phi) \Big\vert & \leq C \sqrt{\frac{\eta_N}{N}} \sup_{r_1 > 1, r_2 > 1} \left( \frac{1 \wedge (\delta_N r_2)^{\dimension}}{r_1^{\impact_1} r_2^{\impact_2}} \right) \\
		& \leq C \sqrt{\frac{\eta_N}{N}} (\delta_N)^{\dimension \wedge \impact_2}.
	\end{align*}
	In both regimes, the right-hand side converges to zero as $N \to \infty$ by \eqref{assumption_jumps_one_tail} and \eqref{assumption_jumps_two_tail}, which concludes the proof of the lemma.
\end{proof}

\subsubsection{Spatial Correlations} \label{sectionspatialcorrelations}
 In this section, we prove \Cref{lem:convergenceofquadraticvariation} to characterize the correlations of the limiting noises in \eqref{SPDE}. We first give a motivating calculation. We have already seen that
\begin{equation*}
	\blangle M^N (\phi) \brangle_t = \eta_N \int_0^t \blangle \Gamma^{\nu_\alpha^N} (\boldsymbol{\rho}_s^N) , \phi \otimes \phi \brangle ds,
\end{equation*}
where $\Gamma^\nu$ was defined in \eqref{eq:Gamma}, and claim
\begin{equation*}
	\eta_N \int_0^t \blangle \Gamma^{\nu_\alpha^N} (\boldsymbol{\rho}_s^N) , \phi \otimes \phi \brangle ds \xrightarrow[N \to \infty]{} t \blangle \mathcal{Q}, \phi \otimes \phi \brangle,
\end{equation*}
where $ \mathcal{Q} $ is defined as
\begin{equation} \label{equationspatialcorrelations}
	\mathcal{Q}(dx_1 dk_1 dx_2 dk_2) = \gamma K_\beta (dx_1, dx_2) (dk_1 \delta_{k_1} (dk_2) - dk_1 dk_2).
\end{equation}
The measure $ K_\beta(d x_1, d x_2) $ encodes the spatial correlations and coalescence rates in each regime and was defined in \eqref{eq:correlation}. The form of \eqref{equationspatialcorrelations} follows from $ \Gamma^{\nu_\alpha^N} $ by letting $ \boldsymbol{\rho}_t^N \rightarrow \lambda $. Recall that
\begingroup
\allowdisplaybreaks
\begin{multline} \label{equationcorrelationsmotivation}
	\eta_N \Gamma^{\nu_\alpha^N} (\boldsymbol{\rho}_s^N) (x_1, x_2, dk_1, dk_2)\\
		=\eta_N  \int_{(0,1] \times [0, \infty)^2} (N u)^2 \int_{\mathbb{R}^{\dimension}} \mathds{1}_{B(x_1, \delta_N r_2) \cap B(x_2, \delta_N r_2)} (x) \frac{1}{V_{\delta_N r_1}}\int_{B( x, \delta_N r_1 )}   \\
		 \Bigg[  \boldsymbol{\rho}_s^N (y, dk_1) \delta_{k_1} (dk_2)
		-  \boldsymbol{\rho}_s^N (y, dk_1) \boldsymbol{\rho}_s^N (x_2 , dk_2)\\
	-  \boldsymbol{\rho}_s^N (x_1 , dk_1) \boldsymbol{\rho}_s^N (y , dk_2)
 + \boldsymbol{\rho}_s^N (x_1, dk_1) \boldsymbol{\rho}_s^N ( x_2 , dk_2) \Bigg]\\ dy dx \frac{1}{\delta_N^{\dimension + \alpha} } \nu_N (du, dr_1, dr_2).
\end{multline}
If we replace $ \rho_{s}^N $ by $ \lambda $, the right-hand side of \eqref{equationcorrelationsmotivation} should converge to
\begin{multline*}
		\eta_N \int_{(0,1] \times [0, \infty)^2} (N u)^2 \int_{\mathbb{R}^{\dimension}} \mathds{1}_{B(x_1, \delta_N r_2) \cap B(x_2, \delta_N r_2)} (x) \frac{1}{V_{\delta_N r_1}}\int_{B( x, \delta_N r_1 )}   \\
		 \hspace{0.5cm} \Big[  dk_1 \delta_{k_1} (dk_2)-  dk_1 dk_2 \Big] dy dx \frac{1}{\delta_N^{\dimension + \alpha} } \nu_N (du, dr_1, dr_2)\\
		= \underbrace{\eta_N \int_{(0,1] \times (0, \infty)^2} (N u)^2 V_{\delta_N r_2} (x_1, x_2)  \frac{1}{\delta_N^{\dimension + \alpha} } \nu_N (du, dr_1, dr_2)} \Big[  dk_1 \delta_{k_1} (dk_2)-  dk_1 dk_2 \Big].
\end{multline*}
\endgroup
The term in brackets contains all spatial information and will equal $\gamma K_\beta(dx_1, dx_2) $. For example, in the one-tail regime, if $\beta = \tail + \impact < \dimension$, we obtain
\begin{align*}
	\eta_N \int_{(0,1]\times (0, \infty)^2}& (N u)^2 V_{\delta_N r_2} (x_1, x_2) \delta_N^{-(\dimension + \alpha)} \nu_N (du, dr_1, dr_2) \\
 &\begin{aligned}
 &= \delta_N^{\alpha - \tail - \impact} u_0^2 \int_1^\infty \frac{1}{r_2^{2 \impact}} V_{\delta_N r_2} (x_1, x_2) \delta_N^{-(d + \alpha)} \frac{dr_2}{r_2^{1 + \dimension + \tail - \impact}}\\
	&= u_0^2 \int_1^\infty V_{\delta_N r_2} (x_1, x_2) \delta_N \frac{dr_2}{(\delta_N r_2)^{1 + \dimension + \tail + \impact}}\\
	&= u_0^2 \int_{\delta_N}^\infty V_{ r_2} (x_1, x_2) \frac{dr_2}{r_2^{1 + \dimension +\tail + \impact}}.
 \end{aligned}
 \end{align*}
Taking the limit $\delta_N \to 0$ and abbreviating $\frac{x_1 -x_2}{\Vert x_1 - x_2 \Vert} = e_1$, this converges to
\begin{align*}
&u_0^2 \int_{0}^\infty V_{ r_2} (x_1, x_2) \frac{dr_2}{r_2^{1 + \dimension +\tail + \impact}} \\ 
&= u_0^2 \int_{\mathbb{R}^{\dimension}} \int_0^\infty \mathds{1}_{\{ \vert x_1 -x_2 - y \vert < r_2 \} } \mathds{1}_{\{ \vert 0 -y \vert < r_2\}} \frac{dr_2}{r_2^{1 + \dimension + \beta}} dy\\
&= u_0^2 \int_{\mathbb{R}^{\dimension}} \int_0^\infty \mathds{1}_{\{ \vert e_1 - y \vert < r_2 /\Vert x_2 - x_1 \Vert \} } \mathds{1}_{\{ \vert y \vert < r_2 / \Vert x_2 - x_1 \Vert \}} \frac{dr_2}{r_2^{1 + \dimension + \beta}} \Vert x_1 - x_2 \Vert^{\dimension} dy\\
&= u_0^2 \int_{\mathbb{R}^{\dimension}} \int_0^\infty \mathds{1}_{\{ \vert e_1 - y \vert < r_2 \} } \mathds{1}_{\{ \vert y \vert < r_2 \}} \frac{1}{\Vert x_1 - x_2 \Vert^\beta} \frac{dr_2}{r_2^{1 + \dimension + \beta}} dy\\
 &= u_0^2 C_{d, \beta}^{(2)} K_{\beta} (x_1, x_2). 
\end{align*}
If $\beta > \dimension$, we instead arrive at
\begin{align*}
	\eta_N \int_{(0,1] \times (0, \infty)^2}& u^2  V_{r_2} (x_1,x_2) \nu_N^\alpha (dr_1, dr_2)\\
	&= u_0^2 \eta_N \int_{(0,1] \times (0, \infty)^2} (N u)^2 V_{\delta_N r_2} (x_1,x_2) \delta_N^{- (\dimension + \alpha)} \nu_N (du, dr_1, dr_2) \\
	&= u_0^2 \delta_N^{\alpha - \dimension} \int_1^\infty V_{\delta_N r_2} (x_1,x_2) \delta_N^{- (\dimension + \alpha)} \frac{dr_2}{r_2^{1 + \dimension + \alpha + \impact}} \\
	&= u_0^2 V_1^2 \int_1^\infty \frac{V_{\delta_N r_2} (x_1,x_2)}{V_{\delta_N r_2}^2 }  \frac{dr_2}{r_2^{1 + \alpha + \impact - \dimension}},
\end{align*}
and as this quantity is integrated against test functions
\begin{equation*}
\begin{aligned}
u_0^2 V_1^2 \int_1^\infty	\int_{(\mathbb{R}^{\dimension})^2} \phi (x_1, k_1) \phi (x_2&, k_2) \frac{V_{\delta_N r_2 } (x_1, x_2 )}{V_{\delta_N r_2}^2}  dx_1 dx_2 \frac{dr_2}{r_2^{1 + \alpha + \impact - \dimension}}\\
&\xrightarrow[N \to \infty]{} u_0^2 V_1^2 \int_1^\infty  \frac{dr_2}{r_2^{1 + \alpha + \impact - \dimension}} \int_{\mathbb{R}^{\dimension}}  \phi (x, k_1) \phi (x, k_2) dx.
\end{aligned}
\end{equation*}
Noting $\int_1^\infty r_2^{-(1 + \alpha + \impact - \dimension)} dr_2 = (\alpha + \impact-\dimension)^{-1}$ the above expression matches the corresponding result in \Cref{table:parameters}. The spatial correlations vanish and the noise becomes white in space. A similar calculation holds in the case $\beta = \dimension$. The next proof will demonstrate that all these calculations are rigorously justified.
\begin{proof}[Proof of \Cref{lem:convergenceofquadraticvariation}]
	In this proof, we will focus on the one-tail regime as the two-tail regime can be treated with similar arguments. First, we will consider the case $\beta = \tail + \impact < \dimension$ extending \parencite[Subsection 3.5, Proof of Lemma 3.10 (3.32), stable case]{forien}. 
	Second, we will adapt \parencite[Subsection 3.5, Proof of Lemma 3.10 (3.32), fixed case]{forien} to analyse the case $\beta > \dimension$, where the spatial correlations vanish. Afterwards, we will treat the case $\beta = \dimension$. The covariation process can be written as
	\begin{align}\label{eqspatialcorrelations}
		\blangle M^N (\phi) \brangle_t 
		&= \eta_N \int_0^t \int_{(0,1] \times [0, \infty)^2} (N u)^2 \blangle \Gamma^{(\delta_N r_1, \delta_N r_2)} (\boldsymbol{\rho}_s^N) , \phi \otimes \phi \brangle \delta_N^{- (\alpha +d)} \nu_N (du, dr_1, dr_2) ds, 
	\end{align}
	where $ \Gamma^{( \delta_N r_1, \delta_N r_2)} $ is short for evaluation of $ \Gamma^{\nu} $ at the Dirac measure at $ (1, \delta_N r_1, \delta_N r_2) $. 

    \paragraph{Case $\boldsymbol{\beta < \dimension}$}
Substituting integrals, \eqref{eqspatialcorrelations} translates to:
	\begin{equation} \label{martingaleconvergenconetailregimeeq1}
		\blangle M^N (\phi) \brangle_t = \int_0^t u_0^2 \int_{\delta_N}^\infty \blangle \Gamma^{(\delta_N^{1 - \power} r_2^{\power}, r_2)} (\boldsymbol{\rho}_s^N) , \phi \otimes \phi \brangle \frac{dr_2}{r_2^{1 +\dimension + \tail + \impact}} ds.
	\end{equation}
	In \eqref{martingaleconvergenconetailregimeeq1} we can change the limit of integration of the inner integral from $ [\delta_N, \infty) $ to $ [0, \infty) $, as the difference converges to zero:
	\begin{align*}
		&\int_0^t u_0^2 \int_{0}^{\delta_N} \blangle \Gamma^{(\delta_N^{1 - \power} r_2^{\power}, r_2)} (\boldsymbol{\rho}_s^N), \phi \otimes \phi \brangle \frac{dr_2}{r_2^{1 + \dimension + \tail + \impact}} ds \\
		& \hspace{1cm} \leq 4 u_0 t \Vert \phi \Vert_2^2 \int_0^{\delta_N} V_{r_2}^2  \frac{dr_2}{r_2^{1 + \dimension + \tail + \impact}} \xrightarrow[N \to \infty]{}  0,
	\end{align*}
	where we used again \eqref{boundonintegrals}. 
	Defining
	\begin{equation*}
	 \mathfrak{Q}_{r_2} (x_1, x_2, dk_1, dk_2) := V_{r_2} (x_1, x_2) (dk_1 \delta_{k_1} (dk_2) - dk_1 dk_2)
	\end{equation*} 
	it remains to show
	\begin{equation*}
		\begin{aligned}
			&\int_0^t u_0^2 \int_0^\infty \blangle \Gamma^{(\delta_N^{1 - \power} r_2^{\power}, r_2)} (\boldsymbol{\rho}_s^N), \phi \otimes \phi \brangle \frac{dr_2}{r_2^{1 + \dimension + \tail + \impact}} ds\\
			& \hspace{1cm} \xrightarrow[N \to \infty]{} \int_0^t u_0^2 \int_0^\infty \blangle \mathfrak{Q}_{r_2}, \phi \otimes \phi \brangle \frac{dr_2}{r_2^{1 + \dimension + \tail + \impact}} ds.
		\end{aligned}
	\end{equation*}
We proceed by dominated convergence. Calculations in \Cref{lem:boundondominatingmeasures} provide us with the necessary bound in $ (r_1, r_2) $. To see pointwise convergence of the integrand, note we can write
\begin{equation*} 
		\blangle \Gamma^{(r_1, r_2)} (\boldsymbol{\rho}) , \phi \otimes \phi \brangle = V_{r_2}^2 \blangle \boldsymbol{\rho} , \Psi_{(r_1, r_2)}^{(1)} \brangle + \blangle \boldsymbol{\rho} \otimes \boldsymbol{\rho} , \Psi_{(r_1, r_2)}^{(2)} \brangle
	\end{equation*}
	for two functions $\Psi_{(r_1, r_2)}^{(1)} (x_1, k_1), \Psi_{(r_1, r_2)}^{(2)} (x_1, x_2 , k_1, k_2)$. The first function is defined as
	\begin{equation*}
		\Psi_{(r_1, r_2)}^{(1)} (x_1, k_1) := \frac{1}{V_{r_1}} \int_{B(x_1, r_1)} \overline{\phi} (x, k_1, r_2)^2 dx
	\end{equation*}
	and the second function as
	\begin{equation*}
		\begin{aligned}
			&\Psi_{(r_1, r_2)}^{(2)} (x_1, x_2 , k_1, k_2) \\
			&:= V_{r_2} (x_1, x_2) \phi (x_1 , k_1 ) \phi (x_2, k_2) - \frac{V_{r_2}}{V_{r_1}}\phi (x_1, k_1) \int_{B(x_1,r_1) \cap B(x_2, r_2)} \overline{\phi} (y , k_2 , r_2) dy\\
			& \hspace{1cm} - \frac{V_{r_2}}{V_{r_1}} \phi (x_2, k_2) \int_{B(x_1, r_1) \cap B(x_2,r_2)} \overline{\phi} (y, k_1 , r_2) dy.
		\end{aligned}
	\end{equation*}
Since $\Psi^{(1)} \in \mathcal{S} (\mathbb{R} \times [0,1] )$, we know from \eqref{equationconvergencetolebesgueinequality} that
\begin{equation} \label{boundclt2}
\mathbb{E} \Bigg[ \sup_{t \in [0,T]} \Big\vert \blangle Z_t^N, \Psi^{(1)} \brangle \Big\vert^2 \Bigg]^{1/2} \leq C \max_{q \in \{ 1,2 \}} \max_{\vert \kappa \vert \leq 2 } \Vert \partial_\kappa \Psi^{(1)} \Vert_p.
\end{equation}
The derivatives on the right-hand side are
\begin{equation*}
\frac{\partial}{\partial x_1} V_{r_2}^2 \Psi_{(r_1, r_2)}^{(1)} (x_1, k_1) = \frac{2 V_{r_2}^2}{V_{r_1}} \int_{B(x_1,r_1)} \overline{\phi} (y, k_1, r_2) \overline{\phi_{x}} (y, k_1, r_2) dy
\end{equation*} 
and
\begin{equation*}
\frac{\partial^2}{\partial x_1^2} V_{r_2}^2 \Psi_{(r_1, r_2)}^{(1)} (x_1, k_1) = \frac{2 V_{r_2}^2}{V_{r_1}} \int_{B (x_1,r_1)} \overline{\phi_{x}} (y,k_1, r_2)^2 + \overline{\phi} (y, k_1, r_2) \overline{\phi_{x x}} (y, k_1, r_2) dy.
\end{equation*}
The Schwartz space $\mathcal{S} (\mathbb{R}^{\dimension} \times [0,1])$ is closed under multiplication and averages. This, together with \Cref{lpboundaverages}, allows us to bound \eqref{boundclt2} for $\Psi^{(1)}$ independently of $r_1$ and show pointwise convergence for the first summand. 
To see the pointwise convergence of the second summand note
\begin{equation*}
\mathbb{E} \Big[ \blangle Z_t^N \otimes Z_t^N , \Psi^{(2)} \brangle \Big] \leq \mathbb{E} \Big[ \blangle D^N, \mathds{1}_{[0,t]} \Psi^{(2)} \brangle \Big].
\end{equation*}
\Cref{lem:boundondominatingmeasures2} provides us with a uniform bound on the right-hand side. This implies the pointwise convergence of $\blangle (\boldsymbol{\rho}_t^N - \lambda ) \otimes ( \boldsymbol{\rho}_t^N - \lambda ) , \Psi^{(2)} \brangle$. The identity
\begin{equation*}
\boldsymbol{\rho}_t^N \otimes \boldsymbol{\rho}_t^N - \lambda \otimes \lambda = \lambda \otimes (\boldsymbol{\rho}_t^N - \lambda) + (\boldsymbol{\rho}_t^N - \lambda ) \otimes \lambda + (\boldsymbol{\rho}_t^N - \lambda) \otimes ( \boldsymbol{\rho}_t^N - \lambda)
\end{equation*}
allows us to make the transition to
\begin{equation} \label{equationtransition}
    \blangle \boldsymbol{\rho}_t^N \otimes \boldsymbol{\rho}_t^N - \lambda \otimes \lambda, \Psi^{(2)} \brangle = (N \eta_N)^{-1/2} \blangle Z_t^N, \tilde{\Psi} \brangle + (N \eta_N)^{-1} \blangle Z_t^N \otimes Z_t^N, \Psi^{(2)} \brangle.
\end{equation}
Here, $\tilde{\Psi} (x,k) = \blangle \lambda, \Psi^{(2)} (\cdot ,x , k ) \brangle + \blangle \lambda , \Psi^{(2)} (x,\cdot, k ) \brangle$ and it remains to prove the pointwise convergence of the middle term. We can represent $\blangle \lambda, \Psi^{(2)} (\cdot, x, k) \brangle$ through averaged functions
\begin{align*}
&\blangle \lambda, \Psi^{(2)} (\cdot, x_2, k_2) \brangle \\
&\begin{aligned}
&= \int_{\mathbb{R}^{\dimension}} V_{r_2} (x_1, x_2) \int_0^1 \phi (x_1 , k_1) dk_1 dx_1 \phi (x_2 ,k_2)\\
& \hspace{1cm} - \frac{V_{r_2}}{V_{r_1}} \int_{\mathbb{R}^{\dimension}} \int_0^1 \phi (x_1, k_1 ) dk_1 \int_{B(x_1, r_1)} \mathds{1}_{B(x_2, r_2)} \overline{\phi} (y,k_2, r_2) dy dx_1 \\
& \hspace{1cm} - \frac{V_{r_2}}{V_{r_1}} \phi (x_2, k_2 ) \int_{\mathbb{R}^{\dimension}} \int_{B(x_1, r_1)} \mathds{1}_{B(x_2, r_2)} (y) \int_0^1 \overline{\phi} (y, k_1, r_2) dk_1 dy dx_1
\end{aligned}\\
&\begin{aligned}
&= \overline{\overline{\phi}} (x_2, r_2) \phi (x_2,k_2) V_{r_2}^2 - V_{r_2}^2 \overline{\overline{\phi} (\cdot, r_1) \overline{\phi} (\cdot, k_2, r_2)} (x_2, r_2) - V_{r_2}^2 \overline{\overline{\phi}} (x_2, r_2) \phi (x_2,k_2).
\end{aligned}
\end{align*}
Again, this time for $\blangle \lambda, \Psi^{(2)} (\cdot, x, k) \brangle$, we can bound the right-hand side of the inequality \eqref{boundclt2} independently of $r_1$ and the missing term of \eqref{equationtransition} converges.
	
	\paragraph{Case $\boldsymbol{\beta > \dimension}$}
	
	We use a separate strategy, which extends \parencite{forien} to different radii $r_1, r_2$ for the impacted area and range from where to choose parents. Our strategy is to approximate the four summands of $ \blangle \Gamma^{(r_1, r_2)} (\rho), \phi \otimes \phi \brangle $ by more tractable terms. The first summand equals
	\begin{equation*}
		V_{r_2}^2 \blangle \overline{\rho}_{r_1} , \overline{\phi} (\cdot , r_2)^2 \brangle,
	\end{equation*}
	where we define for $ r > 0 $
	\begin{equation*}
		\overline{\rho}_r (x, dk) := \frac{1}{V_{r}} \int_{B(x,r)} \rho (z, dk) dz.
	\end{equation*}
	The second term is of the form
	\begin{equation} \label{equationspatialcorvanish}
		\int_{\mathbb{R}^{\dimension}} \int_{[0,1]} \int_{B(x, r_2)} \phi (x_1, k_1) dx_1 \overline{\rho}_{r_1}(x, dk_1) \int_{B(x, r_2)} \int_{[0,1]}  \phi (x_2, k_2) \rho (x_2, dk_2) dx_2 dx.
	\end{equation}
	If we replace $ \phi (x_1, k_1) $ by $ \phi (x, k_1) $ and $ \phi (x_2, k_2) $ by $ \phi (x, k_2) $ in \eqref{equationspatialcorvanish}, the distance between the new expression
	\begin{equation*}
		V_{r_2}^2 \blangle \overline{\rho}_{r_1} \cdot \overline{\rho}_{r_2} , \phi \otimes \phi \brangle = V_{r_2}^2 \int_{\mathbb{R}^{\dimension} \times [0,1]^2} \phi (x, k_1) \phi (x, k_2) \overline{\rho}_{r_1} (x, dk_1) \overline{\rho}_{r_2} (x, dk_2) dx
	\end{equation*}
	and \eqref{equationspatialcorvanish} is smaller than
	\begin{equation*}
		\min \Big\{ 2 r_2 V_{r_2}^2 \Vert \phi \Vert_1 \max_{\vert \kappa \vert = 1 } \Vert \partial_\kappa \phi \Vert_\infty,
		 4 V_{r_2}^2 \Vert \phi \Vert_1 \Vert \phi \Vert_\infty \Big\}.
	\end{equation*}
	This can be seen by applying
	\begin{equation*}
		\Bigg\vert \int_{B(x, r) \times [0,1]} (\phi (x_1, k_1) - \phi (x, k_1)) \rho (x_1, dk_1) dx_1 \Bigg\vert \leq \min \Big\{ r V_r \max_{\vert \kappa \vert = 1 } \Vert \partial_\kappa \phi \Vert_1,
		2 V_r \Vert \phi \Vert_\infty \Big\} .
	\end{equation*}
	The last two terms of $ \blangle \Gamma^{(r_1, r_2)}, \phi \otimes \phi \brangle $ can be approximated in the same way to obtain
	\begin{equation*}
		\begin{aligned}
			&\Big\vert \blangle \Gamma^{(r_1, r_2)} (\rho) , \phi \otimes \phi \brangle - V_{r_2}^2 \big( \blangle \overline{\rho}_{r_1} , \overline{\phi} (\cdot , r_2)^2 \brangle - \blangle \overline{\rho}_{r_1} \cdot \overline{\rho}_{r_2} , \phi \otimes \phi \brangle \big) \Big\vert\\
			& \hspace{1cm} \leq \min \Big\{
			8 r_2 V_{r_2}^2 \max_{\vert \kappa \vert = 1 } \Vert \partial_\kappa \phi \Vert_\infty \Vert  \phi \Vert_1,
			16 V_{r_2}^2 \Vert \phi \Vert_\infty \Vert \phi \Vert_1 \Big\}.
		\end{aligned}
	\end{equation*}
	In total, if $\tail + \impact > \dimension$ this leads to 
	\begin{align*}
			&\Bigg\vert \langle M^N (\phi) \brangle_t - u_0^2 \delta_N^{-2d} \int_0^t \int_1^\infty V_{\delta_N r_2}^2 \Big( \blangle \overline{(\boldsymbol{\rho}_s^N)}_{\delta_N r_2^{\power}} , \overline{\phi} (\cdot , \delta_N r_2)^2 \brangle \\
			& \hspace{1cm}- \blangle \overline{(\boldsymbol{\rho}_s^N)}_{\delta_N r_2^{\power}} \cdot \overline{(\boldsymbol{\rho}_s^N)}_{\delta_N r_2}, \phi \otimes \phi \brangle \Big) \frac{dr_2}{r_2^{1 + \dimension+ \tail + \impact}} ds \Bigg\vert\\ 
			&  \leq u_0^2 t \delta_N^{-2\dimension}  8 \Vert \phi \Vert_1 \int_1^\infty \min \Big\{ \delta_N r_2 V_{\delta_N r_2}^2 \max_{\vert \kappa \vert = 1} \Vert \partial_\kappa \phi \Vert_\infty,
			2 V_{\delta_N r_2}^2 \Vert \phi \Vert_\infty
			\Big\} \frac{dr_2}{r_2^{1 + \dimension + \tail + \impact}} \\
			&  \leq u_0^2 t \delta_N^{\tail + \impact -\dimension}  8 \Vert \phi \Vert_1 \int_{\delta_N}^\infty \min \Big\{
		    r_2 V_{r_2}^2 \max_{\vert \kappa \vert = 1} \Vert \partial_\kappa \phi \Vert_\infty,
			2 V_{r_2}^2 \Vert \phi \Vert_\infty \Big\} \frac{dr_2}{r_2^{1 + \dimension + \tail + \impact}} \\
			&  \begin{aligned}
			&\leq u_0^2 t \delta_N^{\tail + \impact -\dimension}  8 \Vert \phi \Vert_1 \max_{\vert \kappa \vert = 1} \Vert \partial_\kappa \phi \Vert_\infty \int_{\delta_N}^1 r_2 V_{r_2}^2  \frac{dr_2}{r_2^{1 + \dimension + \tail + \impact}}\\
			& \hspace{1cm} + u_0^2 t \delta_N^{\tail + \impact -\dimension}  16 \Vert \phi \Vert_1 \Vert \phi \Vert_\infty \int_{1}^\infty V_{r_2}^2  \frac{dr_2}{r_2^{1 + \dimension + \tail + \impact}}. \\
			\end{aligned}
	\end{align*}
	This right-hand side converges to zero as $N \to \infty$. Again from \Cref{theo:convergencetolebesgue} we know $ \boldsymbol{\rho}_t^N \to \lambda $ in probability, which implies $ \overline{(\boldsymbol{\rho}_s^N)}_{\delta_N r_2^{\power}} \to \lambda $ and $ \overline{(\boldsymbol{\rho}_s^N)}_{\delta_N r_2} \to \lambda $ in the same sense. This shows
	\begin{equation*}
		\blangle M_t^N (\phi) \brangle_t \xrightarrow[N \to \infty]{} u_0^2 V_1^2 \int_0^t \int_1^\infty \Big( \blangle \lambda, \phi^2 \brangle - \blangle \lambda \cdot \lambda , \phi \otimes \phi \brangle \Big) \frac{dr_1}{r_1^{1 + \tail + \impact - \dimension}} ds = \blangle \mathcal{Q}, \phi \otimes \phi \brangle.
	\end{equation*}
	
	\paragraph{Case $\boldsymbol{\beta = \dimension}$}

    Similar to the argument following \eqref{firstbounddominating}, we can show the bound
    \begin{equation*}
        \Big\vert \blangle \Gamma^{(r_1, r_2)} (\rho) , \phi \otimes \phi \brangle \Big\vert \leq 4 V_{r_2} \Vert \phi \Vert_1^2.
    \end{equation*}
    The following distance vanishes as $N \to \infty$:
    \begingroup
    \allowdisplaybreaks
    \begin{align*}
        &\begin{aligned}
            &\Bigg\vert \langle M^N (\phi) \brangle_t - u_0^2 \frac{\delta_N^{-2d}}{\log (1/\delta_N)} \int_0^t \int_1^{1/\delta_N} V_{\delta_N r_2}^2 \Big( \blangle \overline{(\boldsymbol{\rho}_s^N)}_{\delta_N r_2^{\power}} , \overline{\phi} (\cdot , \delta_N r_2)^2 \brangle \\
			& \hspace{1cm}- \blangle \overline{(\boldsymbol{\rho}_s^N)}_{\delta_N r_2^{\power}} \cdot \overline{(\boldsymbol{\rho}_s^N)}_{\delta_N r_2}, \phi \otimes \phi \brangle \Big) \frac{dr_2}{r_2^{1 + \dimension+ \tail + \impact}} ds \Bigg\vert\\ 
        \end{aligned}\\
		&\begin{aligned}  
		&\leq u_0^2 t \frac{\delta_N^{-2\dimension}}{\log (1/\delta_N)}  8 \Vert \phi \Vert_1 \int_1^{1/\delta_N}
			\delta_N r_2 V_{\delta_N r_2}^2 \max_{\vert \kappa \vert = 1} \Vert \partial_\kappa \phi \Vert_\infty  \frac{dr_2}{r_2^{1 + \dimension + \tail + \impact}} \\
		& \hspace{1cm} + u_0^2 t \frac{\delta_N^{-2 \dimension}}{\log (1/\delta_N)} \int_{1/\delta_N}^\infty \Big\vert \blangle \Gamma^{(\delta_N r_1, \delta_N r_2)} (\boldsymbol{\rho}_s^N) , \phi \otimes \phi \brangle \Big\vert \frac{dr_2}{r_2^{1 + \dimension + \tail + \impact}}
		\end{aligned}\\
		&\begin{aligned}  
		&\leq u_0^2 t \frac{\delta_N^{-2\dimension}}{\log (1/\delta_N)}  8 \Vert \phi \Vert_1 \int_1^{1/\delta_N} \delta_N r_2 V_{\delta_N r_2}^2 \max_{\vert \kappa \vert = 1} \Vert \partial_\kappa \phi \Vert_\infty  \frac{dr_2}{r_2^{1 + \dimension + \tail + \impact}} \\
		& \hspace{1cm} + u_0^2 t \frac{\delta_N^{-2 \dimension }}{\log (1/\delta_N)} 4 \Vert \phi \Vert_1^2 \int_{1/\delta_N}^\infty V_{\delta_N r_2} \frac{dr_2}{r_2^{1 + \dimension + \tail + \impact}} \to 0.
		\end{aligned}
    \end{align*}
    Substituting $u := \frac{\log(r_2)}{\log (1/\delta_N)}$ we obtain
\begin{align*}
&\begin{aligned}
    &u_0^2 \frac{\delta_N^{-2d}}{\log (1/\delta_N)} \int_0^t \int_1^{1/\delta_N} V_{\delta_N r_2}^2 \Big( \blangle \overline{(\boldsymbol{\rho}_s^N)}_{\delta_N r_2^{\power}} , \overline{\phi} (\cdot , \delta_N r_2)^2 \brangle \\
	& \hspace{1cm}- \blangle \overline{(\boldsymbol{\rho}_s^N)}_{\delta_N r_2^{\power}} \cdot \overline{(\boldsymbol{\rho}_s^N)}_{\delta_N r_2}, \phi \otimes \phi \brangle \Big) \frac{dr_2}{r_2^{1 + \dimension+ \tail + \impact}} ds
\end{aligned}\\
&\begin{aligned}
&= u_0^2 V_1^2 \int_0^t \int_1^{1/\delta_N} \Big( \blangle \overline{(\boldsymbol{\rho}_s^N)}_{\delta_N r_2^{\power}} , \overline{\phi} (\cdot , \delta_N r_2)^2 \brangle \\
	& \hspace{1cm}- \blangle \overline{(\boldsymbol{\rho}_s^N)}_{\delta_N r_2^{\power}} \cdot \overline{(\boldsymbol{\rho}_s^N)}_{\delta_N r_2}, \phi \otimes \phi \brangle \Big) \frac{dr_2}{\log (1 /\delta_N) r_2} ds
\end{aligned}\\
&\begin{aligned}
&= u_0^2 V_1^2 \int_0^t \int_0^{1} \Big( \blangle \overline{(\boldsymbol{\rho}_s^N)}_{\delta_N^{1- \power u}} , \overline{\phi} (\cdot , \delta_N^{1- u})^2 \brangle \\
	& \hspace{1cm}- \blangle \overline{(\boldsymbol{\rho}_s^N)}_{\delta_N^{1 - \power u}} \cdot \overline{(\boldsymbol{\rho}_s^N)}_{\delta_N^{1- u}}, \phi \otimes \phi \brangle \Big) du ds.
\end{aligned}
\end{align*}
Following a similar argument based on dominated convergence as in the case $\beta > \dimension$, we can conclude the proof.
\end{proof}
\endgroup


\appendix

\section{Existence and uniqueness of the SLFV with mutations} \label{sec:proof_existence}

The aim of this section is to prove \cref{thm:existence}.
We adapt the arguments used to prove Theorem~1.2 in \cite{etheridge_rescaling_2020} to the setting of the SLFV with different dispersal and parent search radius.
We start by proving uniqueness of the solution to the martingale problem associated to $(\mathcal{G}, D(\mathcal{G}))$, using a duality argument.

\subsection{The dual of the SLFV with mutations} \label{subsec:duality}

We define an ancestral process $(\mathcal{A}_t, t \geq 0)$ taking values in
\begin{equation*}
\begin{aligned}
&\Bigg\{ \big((\xi^1, S^1), (\xi^2, S^2),\dots, (\xi^l, S^l ) \big) \in  \bigcup_{n = 1}^{N_0} \Big( \mathbb{R}^{\dimension} \times \{ \partial \} \times \mathcal{P} ( \{ 1, \dots , N_0 \}) \Big)^n \\
&\hspace{1cm} : \bigcup_{i = 0}^l S_i = \{ 0,1,2,\dots, N_0 \} \text{ and } \forall 0 \leq i \neq j \leq l \leq N_0: S_i \cap S_j = \emptyset \Bigg\}.
\end{aligned}
\end{equation*}
At each time $t$, $\mathcal{A}_t$ consists of a finite collection of lineages in $\R^{\dimension} \times \lbrace \partial \rbrace$. Their locations are denoted by $(\xi^j_t, 1 \leq j \leq N_t)$, where $N_t$ is the number of lineages in $\mathcal{A}_t$ at time $t$. By $ \partial$ we denote a cemetery state for the lineages. Forwards-in-time, if interpreted as an infinite-population limit, the mutation mechanism corresponds to individuals mutating at rate $\mu$. Backwards-in-time, which corresponds to tracing samples in the population, we therefore see mutation along the ancestral lineages at the same rate $\mu$. However, as we hit a mutation event backwards-in-time, we cannot know the type of the lineage prior to the event and therefore allocate the lineage the type $\partial$. We do not need to know the type of the ancestor prior to the mutation in order to know the distribution of the type of the individuals at the leaves of the coalescent, as we know the form of the mutation kernel. Each of the lineages $\xi^j$ comes along with an element $S^j$ of $\mathcal{P} (\{ 1, \dots, N_0 \} )$ indicating the labels of initial samples descended from the lineage. 

\begin{definition} \label{def:ancestral_process}
Let $\Pi$ be a Poisson point process on $\R_+ \times \R^{\dimension} \times (0,1] \times (0,\infty)^2$ with intensity measure $dt \otimes dx \otimes \nu(du, dr_1, dr_2)$. The initial condition is given by \[\mathcal{A}_0 = \big( (\xi^1, \{ 1 \}) , ( \xi^2 , \{ 2 \}) , \ldots, ( \xi^{N_0} , \{ N_0 \} ) \big),\] where $\xi^j \in \mathbb{R}^{\dimension}$ and $N_0 \in \mathbb{N}$.
At each event $(t,x,u,r_1,r_2) \in \Pi$,
\begin{enumerate}
    \item pick a location $y$ uniformly at random from $B(x,r_1)$, the parent search area,
    \item for each lineage which finds itself in the replacement area just before the event, i.e. such that $\xi^j_{t^-} \in B(x,r_2)$, mark this lineage independently of the others with probability $u$, and let $J \subset \{ 1, \ldots , N_0 \}$ denote the set of marked lineages,
    \item all marked lineages coalesce at time $t$, so the entries $\{ (\xi_{t-}^j , S_{t-}^j) , j \in J \}$ are replaced by the single entry $ (y, \cup_{j \in J} S_{t-}^j)$.
\end{enumerate}
Moreover, each lineage is killed at rate $\mu$ independently of each other, and killed lineages instantly jump to $\partial$, where they are no longer affected by reproduction events.
\end{definition}

Notice that each lineage jumps when it finds itself in the replacement area of an event and when it is marked, which happens at rate
\begin{align*}
    \int_{(0,1] \times (0,\infty)^2} u\,  |B(0,r_2)| \, \nu(du, dr_1, dr_2).
\end{align*}
By \eqref{existence_condition}, this rate is finite, and since we consider only finitely many lineages, the jump rate of $(\mathcal{A}_t, t \geq 0)$ is almost surely bounded by a deterministic constant (the number of lineages cannot increase).
Hence the above definition gives rise to a well-defined process.
The following proposition states the duality relation between the SLFV with mutations and the ancestral process defined above.
Its proof is omitted as it is very similar to the proof of the corresponding result in \cite{etheridge_rescaling_2020} (Proposition~2.2), see also \cite[Theorem~1]{veber_spatial_2015}.

\begin{proposition} \label{prop:duality}
Any $\Xi$-valued Markov process $(\rho_t, t \geq 0)$ solution to the martingale problem associated to $(\mathcal{G}, D(\mathcal{G}))$ is dual to the ancestral process $(\mathcal{A}_t, t \geq 0)$ in the sense that, for any $N_0 \in \mathbb N$, for any $\psi : \R^{\dimension} \to \R$ continuous and compactly supported and any $\chi : [0,1] \to \R$ continuous, for $\psi \boxtimes \chi (x,k) := \psi (x) \chi (k)$,
\begin{align*}
    \mathbb{E}_{\rho_0} \Big[ \langle \rho_t, \psi \boxtimes \chi \rangle^{N_0} \Big] = \int_{(\R^{\dimension})^{N_0}} \mathbb{E}_{(x_1, \ldots, x_{N_0} )} \Bigg[ \prod_{i=1}^{N_t} \int_{[0,1]} (\chi(k))^{\vert S^i_t \vert} \rho_0(\xi^i_t, dk) \Bigg] \prod_{j=1}^{N_0} \psi(x_j) dx_1 \ldots dx_{N_0},
\end{align*}
with the convention that $\rho_0(\partial,dk) = \lambda(dk)$.
Here, $\mathbb{E}_{\rho_0} [\cdot]$ denotes the expectation with respect to the law of $(\rho_t, t \geq 0)$ with initial condition $\rho_0$. The cardinality of the set $S$ is denoted by $\vert S \vert$ and $\mathbb{E}_{(x_1, \ldots, x_{N_0})}[\cdot]$ is the expectation with respect to the law of the ancestral process started from $ \mathcal{A}_0 = \big((x_1, \{ 1 \}, \ldots, (x_{N_0}, \{ N_0 \}) \big)$.
\end{proposition}

Since the ancestral process is well defined and the expectation of the set of functions of the form $\rho \mapsto \langle \rho, \psi \boxtimes \chi \rangle^{N_0}$ characterises the distribution of $\rho$ (it is dense in the set of continuous real-valued functions on $\Xi$, see Lemma~2.1 in \cite{veber_spatial_2015}), \cref{prop:duality} ensures that the solution to the martingale problem associated to $(\mathcal{G}, D(\mathcal{G}))$ is unique (if it exists).

\subsection{Proof of Theorem~\ref{thm:existence}} \label{subsec:proof_existence}

To prove the existence of a solution to the martingale problem associated to $(\mathcal G, D(\mathcal{G}))$, we follow the argument of \cite[Theorem 1.2]{etheridge_rescaling_2020}, i.e. we construct a sequence of processes $(\rho^{(n)}, n \geq 0)$ which converges in distribution to a solution of the martingale problem.
The processes $(\rho^{(n)}, n \geq 0)$ will be constructed as solutions to a sequence of approximating martingale problems, $(\mathcal{G}^{(n)}, D(\mathcal{G}^{(n)}))$, such that the overall jump rate of $(\rho^{(n)}_t, t \geq 0)$ is bounded almost surely for each $n \geq 0$.

To achieve this, consider an increasing sequence of compact subsets of $\R^{\dimension}$, $(E_n, n \geq 0)$, such that $\cup_{n \geq 0} E_n = \R^{\dimension}$, and an increasing sequence $(\nu_n, n \geq 0)$ of finite measures on $(0,1] \times (0,\infty)^2$ which converges to $\nu$ as $n \to \infty$.
We then set, for any $F = F_{f,\phi} \in D(\mathcal{G}^{(n)}) := D(\mathcal{G})$,
\begin{align*}
    \mathcal{G}^{(n)} F_{f,\phi}(\rho) &:= \mu f'(\langle \rho, \phi \rangle) \langle \lambda - \rho, \phi \rangle \\
    &\begin{aligned}
     \hspace{1cm}+ \int_{E_n} \int_{(0,1] \times (0,\infty)^2} \frac{1}{|B(x,r_1)|} \int_{B(x,r_1)} \int_{[0,1]} &\Big( F_{f,\phi}(\Theta_{x,r_2,u,k}(\rho)) - F_{f,\phi}(\rho) \Big) \\
    &\rho(y,dk) dy\, \nu_n(du, dr_1, dr_2)\, dx.
    \end{aligned}
\end{align*}
We then define $(\rho^{(n)}_t, t \geq 0)$ as the unique $\Xi$-valued Markov process which solves the martingale problem associated to $(\mathcal{G}^{(n)}, D(\mathcal{G}^{(n)}))$.
This martingale problem is well posed because the overall jump rate of $(\rho^{(n)}_t, t \geq 0)$ is almost surely bounded (note that the definition of $(\rho^{(n)}_t, t \geq 0)$ amounts to replacing $\Pi$ by a Poisson point process with intensity $\mathds{1}_{x \in E_n} dt \otimes dx \otimes \nu_n(du, dr_1, dr_2)$ - which is a finite measure - in \cref{SLFVdefinition}).

To prove \cref{thm:existence}, it remains to show that $\rho^{(n)}$ converges in distribution as $n \to \infty$ to a solution to the martingale problem associated to $(\mathcal{G}, D(\mathcal{G}))$.

\begin{lemma} \label{lemma:tightness_rho_n}
    For any $T > 0$, the sequence of $\Xi$-valued Markov processes $(\rho^{(n)}, n \geq 0)$ is tight in $D([0,T], \Xi)$.
\end{lemma}

\begin{proof}
    First recall that, by Lemma~2.1 in \cite{veber_spatial_2015}, $\Xi$ is a compact space.
    Therefore, by the Aldous-Rebolledo criterion (see e.g. \cite[Theorem~1.17]{etheridge_introduction_2000}), it is enough to show that both the finite variation parts and the quadratic variations of the real-valued processes $(F_{f,\phi}(\rho^{(n)}_t), t \geq 0)$ are tight, for any $f \in C^1(\R)$ and $\phi \in C_c(\R^{\dimension})$.
    
    By the definition of $(\rho^{(n)}_t, t \geq 0)$, the finite variation part of $(F_{f,\phi}(\rho^{(n)}_t), t \geq 0)$ is given by
    \begin{align} \label{finite_var_part}
        \int_0^t \mathcal{G}^{(n)} F_{f,\phi}(\rho^{(n)}_s) ds.
    \end{align}
    First note that, for any $\rho \in \Xi$,
    \begin{align} \label{bound_mutations}
        \left| f'(\langle \rho, \phi \rangle) \langle \lambda - \rho, \phi \rangle \right| \leq 2 \| \phi \|_1 \sup_{|x| \leq \| \phi \|_1} |f'(x)|.
    \end{align}
    In addition,
    \begin{align*}
        \left| F_{f,\phi}(\Theta_{x,r,u,k}(\rho)) - F_{f,\phi}(\rho) \right| &\leq u \sup_{|x| \leq \| \phi \|_1} |f'(x)| \left| \langle \mathds{1}_{B(x,r)} (\delta_k - \rho), \phi \rangle \right| \\
        &\leq 2 u \sup_{|x| \leq \| \phi \|_1} |f'(x)| \int_{B(x,r)} \sup_{k \in [0,1]} |\phi(y,k)| dy.
    \end{align*}
    Since $\phi$ is continuous and compactly supported, there exists a constant $C > 0$ (depending only on $f$ and $\phi$) such that, for any $\rho \in \Xi$,
    \begin{align} \label{bound_jump}
        \left| F_{f,\phi}(\Theta_{x,r,u,k}(\rho)) - F_{f,\phi}(\rho) \right| \leq C u (r^d \wedge 1) \mathds{1}_{B(x,r) \cap supp(\phi) \neq \emptyset}.
    \end{align}
    Using \eqref{bound_mutations}, \eqref{bound_jump} and the fact that there exists a constant $C > 0$ such that
    \begin{align} \label{bound_support}
        \int_{\R^{\dimension}} \mathds{1}_{B(x,r) \cap supp(\phi) \neq \emptyset} dx \leq C (r^d \vee 1),
    \end{align}
    we see that
    \begin{align*}
        \left| \mathcal{G}^{(n)} F_{f,\phi}(\rho) \right| \leq C \int_{(0,1] \times (0,\infty)^2} u (r_2)^d \nu_n(du, dr_1, dr_2) + C',
    \end{align*}
    for some constants $C, C' > 0$.
    As a result, under condition \eqref{existence_condition}, the modulus of continuity of \eqref{finite_var_part} is almost surely uniformly bounded over all $n \in \mathbb N$, hence the finite-variation parts of $(F_{f,\phi}(\rho^{(n)}_t), t \geq 0)$ are tight.
    
    The quadratic variation of $(F_{f,\phi}(\rho^{(n)}_t), t \geq 0)$ is given by
    \begin{align*}
        \int_0^t \Gamma^{(n)} F_{f,\phi}(\rho^{(n)}_s) ds,
    \end{align*}
    where
    \begin{multline*}
        \Gamma^{(n)} F_{f,\phi}(\rho) := \int_{E_n} \int_{(0,1] \times (0,\infty)^2} \frac{1}{|B(x,r_1)|} \int_{B(x,r_1)} \int_{[0,1]} \Big( F_{f,\phi}(\Theta_{x,r_2,u,k}(\rho)) - F_{f,\phi}(\rho) \Big)^2\\
        \rho(y,dk) dy\, \nu_n(du, dr_1, dr_2)\, dx.
    \end{multline*}
    But, using \eqref{bound_jump} and the fact that $u \in (0,1]$, there exists a constant $C > 0$ such that, for all $\rho \in \Xi$,
    \begin{align*}
        \left| F_{f,\phi}(\Theta_{x,r,u,k}(\rho)) - F_{f,\phi}(\rho) \right|^2 \leq C u (r^d \wedge 1) \mathds{1}_{B(x,r) \cap supp(\phi) \neq \emptyset}.
    \end{align*}
    As a result there exists a constant $ C > 0 $ such that, for all $\rho \in \Xi$,
    \begin{align*}
        \left| \Gamma^{(n)} F_{f,\phi}(\rho) \right| \leq C,
    \end{align*}
    and this shows that the quadratic variation processes of $(F_{f,\phi}(\rho^{(n)}_t), t \geq 0)$ are tight.
    
    By the Aldous-Rebolledo criterion, this concludes the proof of the lemma.
\end{proof}

We can now conclude the proof of \cref{thm:existence}.

\begin{proof}[Proof of \cref{thm:existence}]
    Uniqueness of the solution to the martingale problem associated to $(\mathcal{G}, D(\mathcal{G}))$ follows by the duality relation of \cref{prop:duality}.
    To prove existence, by \cref{lemma:tightness_rho_n}, there exists a converging subsequence of the sequence $(\rho^{(n)}, n \geq 0)$, and it is enough to show that the limit of any such subsequence must solve the martingale problem associated to $(\mathcal{G}, D(\mathcal{G}))$ (we will also obtain the convergence in distribution of the whole sequence $(\rho^{(n)}, n \geq 0)$).
    By Theorem~4.8.10 in \cite{ethier_markov_1986}, it suffices to prove that, for any $F_{f,\phi} \in D(\mathcal{G})$,
    \begin{align*}
        \mathcal{G}^{(n)} F_{f,\phi} (\rho) \to \mathcal{G} F_{f,\phi}(\rho) \text{ as } n \to \infty, \text{ uniformly in $\rho \in \Xi$.}
    \end{align*}
    But, using \eqref{bound_jump} and the fact that $\nu_n \leq \nu$, for any $\rho \in \Xi$,
    \begin{multline*}
        \left| \mathcal{G}^{(n)} F_{f,\phi}(\rho) - \mathcal{G} F_{f,\phi}(\rho) \right| \leq C \int_{E_n^c} \int_{(0,1] \times (0,\infty)^2} u (r_2^d \wedge 1) \mathds{1}_{B(x,r_2) \cap supp(\phi) \neq \emptyset} \nu(du, dr_1, dr_2) dx \\ + C \int_{\R^{\dimension}} \int_{(0,1] \times (0,\infty)^2} u (r_2^d \wedge 1) \mathds{1}_{B(x,r_2) \cap supp(\phi) \neq \emptyset} (\nu-\nu_n)(du, dr_1, dr_2).
    \end{multline*}
    The first term on the right hand side tends to zero as $n \to \infty$ by Lebesgue's dominated convergence theorem, using \eqref{existence_condition}, and so does the second term, using \eqref{bound_support}.
    Since the right hand side does not depend on $\rho$, the convergence is uniform on $\Xi$, and the result is proved.
\end{proof}
\section{Definition of worthy of martingale measures}
\begin{definition}[{\parencite[Chapter 2]{walsh_introduction_1986}, \parencite[Definition 1.1]{Cho1995}}] \label{def:martingalemeasure}
Let $ (\Omega, \mathcal{F}, \mathcal{F}_t , \mathbb{P}) $ be a filtered probability space and $ \mathcal{B} (\mathbb{R}^{\dimension}) $ be the Borel-$ \sigma $-algebra on $ \mathbb{R}^{\dimension} $. A random set function $ M : \mathbb{R}_+ \times \mathbb{R}^{\dimension} \rightarrow \mathbb{R} $ is a martingale measure if the following three properties are satisfied:
\begin{enumerate}
\item
For disjoint $ A, B \in \mathcal{B}(\mathbb{R}^{\dimension}) , A \cap B = \emptyset$ the set function $ M $ is additive, i.e. $ M(t, A \cup B) = M(t, A) + M(t,B) $ almost surely.
\item
For $ A \in \mathcal{B} (\mathbb{R}^{\dimension}) $, $ \big(M(t, A), t \geq 0 \big) $ is a martingale.
\item
There exists an increasing sequence $ (E_n, n \geq 1) $ where $ E_n \in \mathcal{B}(\mathbb{R}^{\dimension}) $ whose union is $ \mathbb{R}^{\dimension} $ such that
\begin{enumerate}[label=(\roman*)]
\item
for every $ t > 0 $, if $ \mathcal{E}_n = \mathcal{B} (E_n) $ then $ \sup_{A \in \mathcal{E}_n} \mathbb{E} \big(M(t, A)^2 \big) < \infty $;
\item
and for every $ t > 0 $, if $ (A_k, k \geq 1) $ is a sequence of $ \mathcal{E}_n $ decreasing to $ \emptyset $, then
\begin{equation*}
\mathbb{E} \big( M(t, A_k)^2 \big) \rightarrow 0.
\end{equation*}
\end{enumerate}
\end{enumerate}
\end{definition}
\begin{definition}[{\parencite[p.291]{walsh_introduction_1986}}] \label{def:worthy}
A martingale measure $ M $ is worthy if there exists a random $ \sigma $-finite measure $ D $ on $ \mathcal{B} (\mathbb{R}_+ \times \mathbb{R}^{\dimension} \times \mathbb{R}^{\dimension}) $ with the following properties:
\begin{enumerate}
\item
$ D $ is symmetric with respect to the two spatial marginals and for fixed $ A, B \in \mathcal{B} (\mathbb{R}^{\dimension}) $ the process $ \big( D ([0,t] \times A \times B) , t \geq 0 \big) $ is predictable.
\item
For any $ (s,t] \times A \times B \subset \mathbb{R}_+ \times \mathbb{R}^{\dimension} \times \mathbb{R}^{\dimension} $ the covariance functional is bounded by $ D $, i.e.
\begin{equation*}
\Big\vert \blangle M(A), M(B) \brangle_t - \blangle M(A), M(B) \brangle_s \Big\vert \leq D \big((s,t] \times A \times B \big).
\end{equation*}
\end{enumerate}
\end{definition}

\section{Convergence Theorems}

\begin{theorem}[{\parencite[Theorem 3.7]{forien}}] \label{theo:mainconvergencetheorem}
Let $ (M^N, N \geq 1) $ be a sequence of worthy martingale measures on $ \mathbb{R}^{\dimension} \times [0,1] $ with corresponding dominating measures $ (D^N, N \geq 1) $. Let \begin{equation} \psi^N : \big\{ (s,t) : 0 \leq s \leq t \big\} \times \mathbb{R}^{\dimension} \times [0,1] \rightarrow \mathbb{R}
\end{equation} be a sequence of functions and $ U_t^N $ be the sequence of stochastic integrals
\begin{equation*}
	U_t^N := \int_{[0,t] \times \mathbb{R}^{\dimension} \times [0,1]} \psi_{s,t}^N (x,k) M^N (ds, dx, dk).
\end{equation*}
If $ (\psi^N, N \geq 1) $ and $ (D^N, N \geq 1) $ satisfy the conditions below with constants $ C_1, C_2, C_3 > 0$ and $k \in \mathbb{N}  $, the sequence $ U_t^N $ is tight in $ D (\mathbb{R}_+, \mathbb{R}) $ and for all $ T > 0 $ there exists a constant $ C(T) $ such that for all $ N \geq 1 $
\begin{equation*}
	\mathbb{E} \Bigg( \sup_{t \in [0, T]} \vert U_t^N \vert^2 \Bigg) \leq C(T) C_1 k (C_2^2 + C_3^2).
\end{equation*}
The conditions are:
\begin{enumerate}
\item
There exists $ C_1 >0 $ and $ k \in \mathbb{N} $ such that for any $ 0 \leq s \leq t$, for all $ N \geq 1 $ and any Schwartz function $ \phi \in \mathcal{S} (\mathbb{R}^{\dimension} \times [0,1]) $
\begin{equation} \label{equationtheoremdominatingmeasures}
	\int_{[s,t] \times ( \mathbb{R}^{\dimension} \times [0,1])^2} \phi (x_1,k_1) \phi (x_2, k_2) D^N (ds, dx_1, dk_1, dx_2, dk_2) \leq C_1 \vert t - s \vert \sum_{p = 1}^k \Vert \phi \Vert_p^2.
\end{equation}
\item
For any $ 0 \leq s \leq t $ and for all $ N \geq 1 $ the function $ \psi_{s,t}^N $ is in $ \mathcal{S} (\mathbb{R}^{\dimension} \times [0,1]) $ and the functions $ t \mapsto \psi_{s,t}^N $ and $ s \mapsto \psi_{s,t}^N $ are both continuous.
\item
There exist $ C_2, C_3, \mu > 0 $ such that for all $ N  \geq 1, 0 \leq s \leq t $ and $ p \in [1,k] $
\begin{equation*}
	\big\Vert \psi_{s,t}^N \big\Vert_p \leq C_2 e^{- \mu (t- s)},
\end{equation*}
and for all $ N \geq 1, 0 \leq s' \leq s \leq t \leq t', p \in [1,k] $
\begin{equation*}
	\big\Vert \psi_{s,t'}^N - \psi_{s,t}^N \big\Vert_p \leq C_3 \vert t'-  t\vert e^{- \mu (t-s)}, \hspace{1cm} \big\Vert \psi_{s', t}^N - \psi_{s,t}^N \big\Vert_p \leq C_3 \vert s' - s \vert e^{- \mu (t-s)}.
\end{equation*}
\end{enumerate}
\end{theorem}

\begin{theorem}[{\parencite[Theorem VIII 3.11]{JJ2003}}] \label{theo:martingaleconvergence}
	Let $ ((X_t^N)_{t \geq 0}, N \geq 1) $ be a sequence of c\`adl\`ag, locally square-integrable, d-dimensional martingales and let $ (X_t, t \geq 0) $ be a continuous, d-dimensional Gaussian martingale. We write $ X_t^N = \big(X_t^{N,1}, ..., X_t^{N,d} \big) $ and $ X_t = \big(X_t^1,..., X_t^d \big) $ for the corresponding vectors. Suppose that
	\begin{enumerate}
		\item
		$ \sup_{t \geq 0} \big\vert X_t^N - X_{t^-}^N \big\vert $ is uniformly bounded for $ N \geq 1 $ and converges to zero in probability,
		\item
		for each $ t \in \mathcal{Q} $ and $ 1\leq i,j \leq d $, the pairwise predictable quadratic variations
		\begin{equation*} \blangle X^{N,i}, X^{N,j} \brangle_t \xrightarrow[N \to \infty]{} \blangle X^i, X^j \brangle_t 
		\end{equation*} 
		converge in probability.
	\end{enumerate}
	Then $ X^N $ converges to $ X $ in distribution in $ D(\mathbb{R}_+, \mathbb{R}^{\dimension}) $.
\end{theorem}

 \begin{theorem}[{\parencite[Proposition E.1]{forien}}]
 \label{theo:convergenceoffinitedimensionaldistributions}
 Consider a sequence of worthy martingale measures $(M^N, N \geq 1)$ on $\mathbb{R}^{\dimension} \times [0,1]$ with corresponding dominating measures ${(D^N, N \geq 1)}$. Suppose that $(M^N, N \geq 1)$ converges in distribution to a martingale measure $M$ in $D(\mathbb{R}_+ , \mathcal{S}' ( \mathbb{R}^{\dimension} \times [0,1] ))$ and that the dominating measures are bounded in the sense of \eqref{equationtheoremdominatingmeasures} for a $k\geq 1$. If a family of deterministic real-valued functions $(f_i^N, N \geq 1, i \leq n)$ on $[0,T] \times \mathbb{R}^{\dimension} \times [0,1]$ satisfies
 \begin{enumerate}
     \item 
     for all $N \geq 1, 1 \leq i \leq n$ and $s \in [0,T]$ the function $f_i^N (s, \cdot)$ is contained in $\mathcal{S}(\mathbb{R}^{\dimension} \times [0,1] )$,
     \item
     the functions $f_i^N$ are uniformly bounded in $L^p$ for any $p \in [1,k]$ and $1 \leq i \leq n$,
     \begin{equation*}
         \sup_{N \geq 1} \sup_{s \in [0,T]} \big\Vert f_i^N (s, \cdot) \big\Vert_p < + \infty,
     \end{equation*}
     \item
     there exists limiting functions $f_i: [0, T] \times \mathbb{R}^{\dimension} \times [0,1] \rightarrow \mathbb{R}, \enspace 1 \leq i \leq n$ such that the convergence is uniform in $L^p$ for $p \in [1,k]$,
     \begin{equation*}
         \lim_{N \to \infty} \sup_{s \in [0,T]} \big\Vert f_i^N ( s, \cdot) - f_i (s, \cdot) \big\Vert_p = 0,
     \end{equation*}
 \end{enumerate}
 then the finite-dimensional distributions 
 \begin{equation*}
     \Big( M_t^N (f_1^N) ,..., M_t^N (f_n^N) \Big)_{t \in [0,T]} \xrightarrow[N \to \infty]{} \Big( M_t (f_1) ,..., M_t (f_n) \Big)_{t \in [0,T]}
 \end{equation*}
 converge in distribution in $D([0,T], \mathbb{R}^n)$. 
 \end{theorem}

\printbibliography




\end{document}